\documentclass[preprint,12pt]{elsarticle}

\usepackage{amssymb}
\usepackage{amsmath}

\usepackage{multirow}
\usepackage{arydshln} 
\usepackage{bm}
\usepackage{mathtools}
\usepackage{graphicx}
\usepackage[colorlinks=true, allcolors=blue]{hyperref}
\usepackage{hyperref}
\usepackage{subfigure}
\usepackage{booktabs} 
\usepackage{caption}  
\usepackage{algorithm}
\usepackage{algorithmicx}
\usepackage[algo2e,linesnumbered,lined,ruled]{algorithm2e} %

\usepackage{xcolor} 

\newcommand{\new}[1]{{\color{black} #1}}

\newtheorem{remark}{Remark}

\begin{document}

\begin{frontmatter}

\title{An Adaptive Subdomain Coupling Approach in Domain Decomposition for Multiphase Porous Media Flow} 

\author[2]{Shizhe Li}
\author[1,2]{Li Zhao\corref{cor1}}\ead{lizhao@lsec.cc.ac.cn}
\author[2]{Chen-Song Zhang}

\affiliation[1]{organization={College of Computer Science and Technology, National University of Defense Technology},
            addressline={Changsha}, 
            postcode={410073}, 
            country={P. R. China}}

\affiliation[2]{organization={SKLMS, Academy of Mathematics and Systems Science, Chinese Academy of Sciences, and School of Mathematical Sciences, University of Chinese Academy of Sciences},
            addressline={Beijing}, 
            postcode={100190}, 
            country={P. R. China}}

\cortext[cor1]{Corresponding author.}

\begin{abstract}
\vskip -0.2cm
The numerical simulation of large-scale multiphase flow in porous media is of considerable importance across various application fields, particularly in the petroleum industry. The fully implicit method is preferred in reservoir simulations owing to its superior numerical stability and more relaxed time step constraints.
However, this method requires solving a large nonlinear system, which becomes highly nonlinear in complex heterogeneous media with small grid scales, emphasizing the need for efficient and convergent numerical methods to accelerate nonlinear solvers on parallel computing systems. In this paper, we present an adaptively coupled subdomain framework based on domain decomposition methods. This framework effectively handles strong local nonlinearities in global problems by solving subproblems within the coupled regions. Furthermore, we propose several adaptive coupling strategies and present a novel method for calculating initial guesses, aimed at improving the convergence and scalability of nonlinear solvers. A series of numerical experiments validate the effectiveness and robustness of the proposed framework. Additionally, large-scale reservoir simulations demonstrate that the proposed method achieves competitive parallel performance.
\footnotetext{The work was partially supported by the National Science and Technology Major Project of China (No. 2025ZD1406105) and the National Natural Science Foundation of China (No. 12571445).}
\end{abstract}


%
%
%


\end{frontmatter}



\section{Introduction}
Predicting multiphase flow in porous media is a key area of research in petroleum reservoir simulation~\citep{Aziz1979,chen_computational_2006,ZHANG2022}.
As reservoir development becomes more complex and the need for improved resource utilization efficiency grows, traditional coarse-grid simulation methods are no longer sufficient to capture the intricate geological features and fluid dynamics within reservoirs.
Modern reservoirs, characterized by high heterogeneity, complex fault systems, and fracture networks, require refined simulations to accurately predict fluid flow and optimize extraction strategies. 
However, the high grid resolution required for  these refined simulations considerably increases computational costs.
At the same time, advancements in computer hardware, particularly in parallel computing architectures, have significantly boosted computational power, enabling the simulation of larger and more complex models within practical time frames.

In high-resolution reservoir simulation, the fully implicit method (FIM)~\\\citep{DouglasFIM1959} is one of the most robust approaches, offering unconditional stability and allowing the relaxation of the Courant-Friedrichs-Lewy (CFL) condition~\citep{CoatsCFL2003}.
Using a fully implicit method requires solving a large nonlinear system at each time step. 
The standard approach to addressing this nonlinearity involves variations of Newton iterations, where a system of linear equations must be solved during each iteration~\citep{YANG2018}.
This process incurs significant computational costs, making it the primary expense in the simulation~\citep{Chensong2022}.
Currently, extensive research is focused on accelerating the solution of linear systems that arise from nonlinear equations in reservoir simulations~\citep{2014A2,YANG2014417,Feng2024,LiZheng2017,WANG2018443,YANG20192,Zhao2022,Zhao2023}.

Nonlinear preconditioning techniques provide an alternative approach by targeting the elimination of imbalanced nonlinearities within the system. This improvement enhances the global convergence properties of nonlinear methods (such as the Newton method), thereby reducing the number of global linear iterations required.
Imbalances in nonlinearity typically arise from factors such as discontinuities in permeability coefficients, wide variations in fluid properties, strong capillary effects with limited spatial extent, complex source terms, and singularities at corners, faults, or voids. In such cases, the Newton method may experience poor convergence, potentially leading to stagnation or divergence~\citep{luo2021nonlinear}.
Similar to linear preconditioning, nonlinear preconditioning can be applied to either the left or right side of nonlinear functions.
Left nonlinear preconditioners, such as additive Schwarz preconditioned inexact Newton (ASPIN) method~\citep{doi:10.1137/S106482750037620X,hwang2007class}, multiplicative Schwarz preconditioned inexact Newton (MSPIN) method~\citep{doi:10.1137/140970379,liu2016convergence,Liu2023OverlappingMS}, and restricted additive Schwarz preconditioned exact Newton (RASPEN) method~\citep{dolean:hal-01171167}, solve local nonlinear problems to provide preconditioning for the global nonlinear problem, thereby improving its convergence.
In contrast, right nonlinear preconditioners, such as the nonlinear elimination (NE) method~\citep{doi:10.1137/080736272,yang2016active,yang2018adaptive}, can be viewed as an inner correction step before the global Newton iterations aimed at preconditioning areas with strong nonlinearities in the solution.

ASPIN was introduced by~\citet{skogestad2013domain} for solving multiphase flow problems in porous media, demonstrating its potential in addressing challenging problems. \citet{klemetsdal2022numerical} evaluated the robustness of the ASPIN method across various complex scenarios, particularly in fractured reservoirs and three-phase compositional models. 
They also investigated the method's sensitivity to the pattern of domain decomposition.
Additionally, \citet{luo2021nonlinear} proposed and compared several different NE strategies for two-phase flow problems discretized using the fully implicit discontinuous Galerkin (DG) finite element method. The results demonstrated the superiority of the proposed methods over the classical Newton approach.
Furthermore, \citet{liu2018note} proposed an adaptive nonlinear preconditioning framework based on convergence monitors, allowing nonlinear preconditioning to be turned off during outer Newton iterations when it is not needed, thereby reducing computational costs while maintaining robustness.

To the best of our knowledge, there is limited work focused on efficient nonlinear algorithms for multiphase flow in porous media in large-scale parallel computing. In parallel computing, the ASPIN method is closely tied to domain decomposition. In the classical ASPIN approach, subproblems are defined independently within subdomains, with one or more processes assigned to solve them. As a result, the quality of the preconditioner is dependent on the domain decomposition pattern.
As the number of processes and subdomains increases, the convergence performance of traditional single-level additive Schwarz methods (ASM) inevitably deteriorates~\citep{xu1998some}. Although a two-level strategy can improve convergence, constructing an efficient coarse-grid problem remains challenging.
Furthermore, in large-scale, refined reservoir simulations, the coarse-grid problem itself can become sufficiently large to pose challenges similar to those of the original problem. Additionally, efficiently implementing multi-level algorithms requires substantial modifications to existing code.

To address these issues, a potential solution is to dynamically merge the original subdomains into larger subdomains during the simulation and define subproblems within these newly formed larger subdomains, to be collectively solved by all the processes originally assigned to the individual subdomains. This approach not only preserves the important couplings between subdomains but also enhances the ability to capture strong local nonlinearities. 
An appropriate coupling pattern of subdomains is expected to improve the convergence performance of single-level ASM and accelerate the solution process.

In this paper, we propose a novel adaptively coupled domain decomposition method (ADDM) designed for large-scale multiphase flow in porous media. The main contributions of this work are:

\begin{itemize}

\item We develop an efficient subdomain coupling framework based on domain decomposition methods to improve both the convergence and parallel performance of nonlinear solvers, overcoming the limitations of classical domain decomposition techniques.

\item We introduce several physics-based adaptive coupling strategies and utilize subproblem solutions defined on the coupling regions as initial guesses to accelerate the Newton iterations for the global problem.

\item The proposed methods are implemented into the open-source simulator OpenCAEPoro~\citep{OCP2024,ShizhePHD2024} to simulate multiphase flow in porous media. The efficiency and robustness of these methods are validated through complex heterogeneous media cases and extremely large-scale simulations.

\end{itemize}

The structure of the paper is outlined as follows. Section \ref{sec:Mathematical Model and Discretization Method} presents the mathematical model for multiphase and multicomponent flows in porous media, followed by the corresponding fully implicit discretization. In Section \ref{sec:ADDM}, we provide a comprehensive introduction to subdomain adaptively coupled decomposition methods, which is central to the proposed nonlinear solver. In Section \ref{sec:Numerical Experiments}, we evaluate the effectiveness and parallel performance of the proposed methods through numerical experiments. Finally, we conclude the paper by summarizing the work presented in Section \ref{sec:Conclusions}.

\section{Mathematical model and discretization method}\label{sec:Mathematical Model and Discretization Method}
This section reviews the governing equations and discretization methods employed in the simulation of multiphase and multicomponent flow within porous media.

\subsection{Mathematical model}\label{sec:model}
We consider an isothermal multicomponent model that includes $n_{c}$ components and $n_{p}$ phases~\citep{chen_computational_2006}. \new{Let $\Omega$ denote the physical domain, and let $t$ be the time variable defined on the temporal interval $\mathcal{T}$. The mass conservation equation for each component $i~(i=1, \ldots, n_c)$ is given by}
\begin{equation} \label{eq:mass-conservation}
		\frac\partial{\partial t}{\new{\Big(\phi \sum_{j=1}^{n_p} x_{ij}\,\xi_j\,S_j\Big)}}+\nabla\cdot\sum_{j=1}^{n_p}\Big(x_{ij}\xi_{j}\mathbf{u}_j-\xi_{j}\mathbf{D}_{ij}\nabla x_{ij} \Big)=Q_i,\quad \text{in}~\Omega \times \mathcal{T}, 
\end{equation}
where \new{$\phi$ is the porosity of the porous medium, $S_{j}$ is the saturation of phase $j$,} $x_{ij}$ is mole fraction of component $i$ in phase $j$, $\xi_{j}$ is molar density of phase $j$, $\mathbf{u}_{j}$ is volumetric flow rate of phase $j$, and $\mathbf{D}_{ij}$ is diffusion coefficient tensor of component $i$ in phase $j$. $Q_{i}$ is volumetric molar injection or production rate for component $i$. Wells are described using a standard Peaceman well model~\citep{Peaceman_1978}. 

Based on Darcy's Law, we have
\begin{equation} \label{eq:darcy-law}
		\mathbf{u}_j=-\frac{\boldsymbol{\kappa}\kappa_{rj}}{\mu_j}(\nabla P_j-\rho_{j}g\nabla z),\quad j=1, \ldots, n_p,
\end{equation}
where $\boldsymbol{\kappa}$ is effective permeability of rock, $\kappa_{rj}$ is relative permeability of phase $j$, $\mu_{j}$ is viscosity coefficient of phase $j$, $P_{j}$ is pressure of phase $j$, $\rho_{j}$ is mass density of phase $j$, $g$ is gravity acceleration, and $z$ is depth.

Additionally, certain constraints must be imposed on these physical quantities.

\begin{itemize}
    \item Saturation constraint equation:
	\begin{equation}\label{eq:saturation}
		\sum_{j=1}^{n_p} S_{j} = 1.
	\end{equation}
    Alternatively, the equivalent relationship:
    \begin{equation*}
        V_{f}-V_{p}=0,
    \end{equation*}
    where $V_{f}$ is fluid volume, and $V_{p}$ is pore volume.
	\item Molar fraction constraint equation:
	\begin{equation}\label{eq:molar-fraction}
	\sum_{i=1}^{n_c} x_{ij} = 1, \quad j=1, \ldots, n_p.
	\end{equation}
	
	\item Capillary pressure equation:
	\begin{equation}\label{eq:capillary}
		 P_j = P - P_{cj}, \quad j=1, \ldots, n_p,
	\end{equation}
	where $ P $ is pressure of reference phase, and $ P_{cj} $ is capillary pressure between the reference phase and phase $ j $.
\end{itemize}

\new{
Note that the system involves $n_c n_p + 3 n_p + 1$ unknown variables, namely $x_{ij}$, $\mathbf{u}_{j}$, $P_{j}$, $P$, and $S_{j}$, with $i = 1, \ldots, n_c$ and $j = 1, \ldots, n_p$. Therefore, $n_c n_p + 3 n_p + 1$ independent relations are required to uniquely determine the solution of the system. Equations~\eqref{eq:mass-conservation}--\eqref{eq:capillary} provide $n_c + 3 n_p + 1$ independent relations, either differential or algebraic. The remaining $n_c (n_p - 1)$ relations are supplied by the equilibrium conditions~\eqref{eq:fugacity}, which relate the numbers of moles of each component among different phases. Moreover, the above system, together with the initial and boundary conditions~\eqref{eq:Comp-initial-Sj}--\eqref{eq:Comp-boundary-mass}, is well posed.
\begin{itemize}
	\item Initial conditions:
	\begin{equation}\label{eq:Comp-initial-Sj}
	\begin{split}
		S_{j}(\mathbf{x}, 0) &= S_{j}^{0}(\mathbf{x}), \quad j = 1,\ldots,n_p-1, \\
		P(\mathbf{x}, 0) &= P^{0}(\mathbf{x}),
	\end{split}
	\end{equation}
	where $S_{j}^{0}(\mathbf{x})$ and $P^{0}(\mathbf{x})$ are given known functions for $\mathbf{x} \in \Omega$. 

	\item Boundary conditions:
	\begin{equation}\label{eq:Comp-boundary-mass}
		\frac{\boldsymbol{\kappa} \kappa_{r j}}{\mu_{j}}
		\Big(\nabla P_{j} - \rho_{j} g \nabla z \Big)
		\cdot \bm{n} \big|_{\partial \Omega} = 0,
		\quad j = 1,\ldots,n_p.
	\end{equation}
\end{itemize}
}

\new{
\subsection{Thermodynamic equilibrium equations}
In multiphase multicomponent flow, phase-equilibrium calculations are used to determine the phase state, phase fractions, and phase compositions at specified pressure $P$, temperature $T$, and overall composition $z_i$. In this work, isothermal conditions are assumed, i.e., the temperature $T$ is constant.

Assuming local thermodynamic equilibrium, the fugacity of each component is equal in all coexisting phases, i.e.,
\begin{equation}\label{eq:fugacity}
\ln f_{ij} = \ln f_{i,n_p}, \quad i = 1, \ldots, n_c,~j = 1, \ldots, n_p-1.
\end{equation}
Here, $f_{ij} = x_{ij}\varphi_{ij}P$ denotes the fugacity of component $i$ in phase $j$, $\varphi_{ij} $ is the fugacity coefficient of component $i$ in phase $j$.

A phase-stability analysis based on Gibbs free-energy minimization is first conducted to determine whether the mixture remains stable as a single phase or splits into multiple phases. Typically, this analysis requires checking only the heaviest phase. If instability is detected, a flash calculation is then performed for a prescribed number of phases. For a multiphase system, the material balance is given by
\begin{equation*}
\begin{split}
x_{ij} &= \frac{ z_i K_{ij} }{1 + \sum_{\ell=1}^{n_p-1} \beta_\ell (K_{i\ell}-1)},\quad i = 1, \ldots, n_c,~j = 1, \ldots, n_p,\\
\end{split}
\end{equation*}
where $K_{ij} = x_{ij} / x_{i,n_p}$, and $\beta_j$ denotes the fraction of phase \( j \), determined from the Rachford-Rice equation,
\begin{equation}\label{eq:RRE}
\sum_{i=1}^{n_c} \frac{z_i(K_{ij}-1)}{1 + \sum_{\ell=1}^{n_p-1} \beta_\ell (K_{i\ell}-1)} = 0,\quad j = 1, \ldots, n_p-1.
\end{equation}

Generally, equation \eqref{eq:RRE} is solved in the successive substitution method to update $K_{ij}$ and obtain initial estimates of $x_{ij}$, followed by the Newton--Raphson method to accelerate the convergence of the nonlinear fugacity-equilibrium system~\eqref{eq:fugacity}. The resulting phase compositions and thermophysical properties are then incorporated into the governing equations of porous-media flow. For more details, the reader is referred to~\citep{michelsen1982isothermal-psa,michelsen1982isothermal-psc}.
}

\subsection{Discretization method}\label{Discretization method}
The equations are discretized in space using a finite-volume method with a two-point flux approximation and upwind weighting~\citep{Aavatsmark2002,Aavatsmark2008}, and in time using an implicit (backward) Euler method.
In each grid cell, we have
\begin{align}
	F_{v} \coloneqq V_{p}^{(n+1)}-V_{f}^{(n+1)}&=0, \label{FIM01}\\
	F_{N,i} \coloneqq \frac{N_{i}^{(n+1)}-N_{i}^{(n)}}{{\delta t}^{(n+1)}}+\sum\limits_{s \in \mathcal{I}}\sum_{j=1}^{n_{p}} T_{s,ij}^{(n+1)}\left[\Delta P_{j}^{(n+1)}-\rho_{j}^{(n+1)}g\Delta z\right]_{s}+Q_{i}^{(n+1)}&=0,  \label{FIM02} 
\end{align}
where \new{$N_{i} \coloneqq \phi \sum_{j=1}^{n_p} x_{ij}\,\xi_j\,S_j$ is molar concentration of component $i~(i=1, \ldots, n_{c})$,} $\mathcal{I}$ represents the set of all interfaces between a grid cell and its neighboring cells, and the operator $\Delta$ denotes the difference in the corresponding variable between the two grid cells on either side of the interface. $[\cdot]_s$ represents the value of the physical quantities at
the interface $s$, approximated using the upstream-weighted method. $T_{s,ij}$ represents the transmissibility of component $i$ in phase $j$ across interface $s$. \new{For a detailed discussion of the discretization, readers are referred to~\citep{chen_computational_2006,Qiao2015,ShizhePHD2024}.}

The above nonlinear system is then solved simultaneously using a fully implicit formulation and a Newton-type method.
In our approach, the primary variables are chosen as $P, N_{1}, \dots, N_{n_{c}}$. The Newton search direction is
\begin{equation}
    \boldsymbol{d} = -\boldsymbol{J}^{-1}\boldsymbol{F},
\end{equation}
where $\boldsymbol{F}$ represents the residual and $\boldsymbol{J}$ is the Jacobian matrix.

\section{Adaptively coupled domain decomposition method}\label{sec:ADDM}
In this section, we first introduce the adaptively coupled domain decomposition method, presenting the motivation for its development and describing three types of adaptive coupling strategies for subdomains.
Next, we discuss the method for setting boundary conditions on the subdomains. 
Finally, we present a nonlinear solution framework based on the adaptively coupled domain decomposition method, which provides efficient and robust initial values for the fully implicit method.

\subsection{Subdomain adaptive coupling strategy}\label{sec:ADDM-subsec1}
Domain decomposition methods (DDM) are naturally well-suited for parallel computing, as they divide the computational domain into multiple subdomains, each of which is managed by one or more processes. The basic idea behind the adaptively coupled domain decomposition method (ADDM) is to dynamically merge existing subdomains into larger subdomains during the simulation process, defining subproblems within these newly formed, larger subdomains. These subproblems are then solved collaboratively by all the processes originally responsible for the smaller subdomains. The strategy of dynamically coupling subdomains during reservoir simulation is based on three fundamental observations:
\begin{enumerate}
    \item Effectively addressing strong local nonlinearities in the global problem is crucial for improving convergence performance.
    \item Regions with strong nonlinearities (e.g., the advancing fluid fronts near wells when fluid injection begins) are typically localized, occupying relatively small areas within the entire computational domain, and they evolve dynamically throughout the simulation.
    \item Most linear and nonlinear solving algorithms achieve significantly higher parallel efficiency with a smaller number of processes, such as tens to hundreds, compared to when the number of processes increases to thousands or tens of thousands, where efficiency tends to decrease rapidly.
\end{enumerate}
This method of dynamically coupling subdomains in a nonlinear manner preserves the essential coupling relationships between them. 
As a result, it is anticipated that an appropriate coupling pattern will effectively capture the nonlinear characteristics of the global problem. 
Consequently, treating the problems defined on these newly formed subdomains as initial values or preconditioners can significantly improve the convergence performance of the global problem, in comparison to the classical single-level ASM.
Furthermore, since regions with strong nonlinearities are typically localized, the number of subdomains requiring coupled solving is generally small. Therefore, solving the problems within these coupled regions tends to achieve high parallel efficiency.

Waterflooding or gas injection is commonly employed to displace one phase, such as oil, through the reservoir. During this process, the fluid dynamics near the advancing front exhibit pronounced nonlinear behavior, particularly when the fluids flow through heterogeneous or complex media. Proper treatment of the nonlinear regions around the advancing front is expected to enhance the convergence performance of the global problem. In the framework of the adaptively coupled domain decomposition method, a natural approach is to dynamically couple the subdomains nonlinearly around the advancing front. This approach involves three key steps: identifying the location of the advancing front, constructing a coupling strength graph for the subdomains, and partitioning the graph to determine the optimal coupling patterns between the subdomains.

\begin{figure}[htpb]
    \centering
    \includegraphics[width=0.3\textwidth]{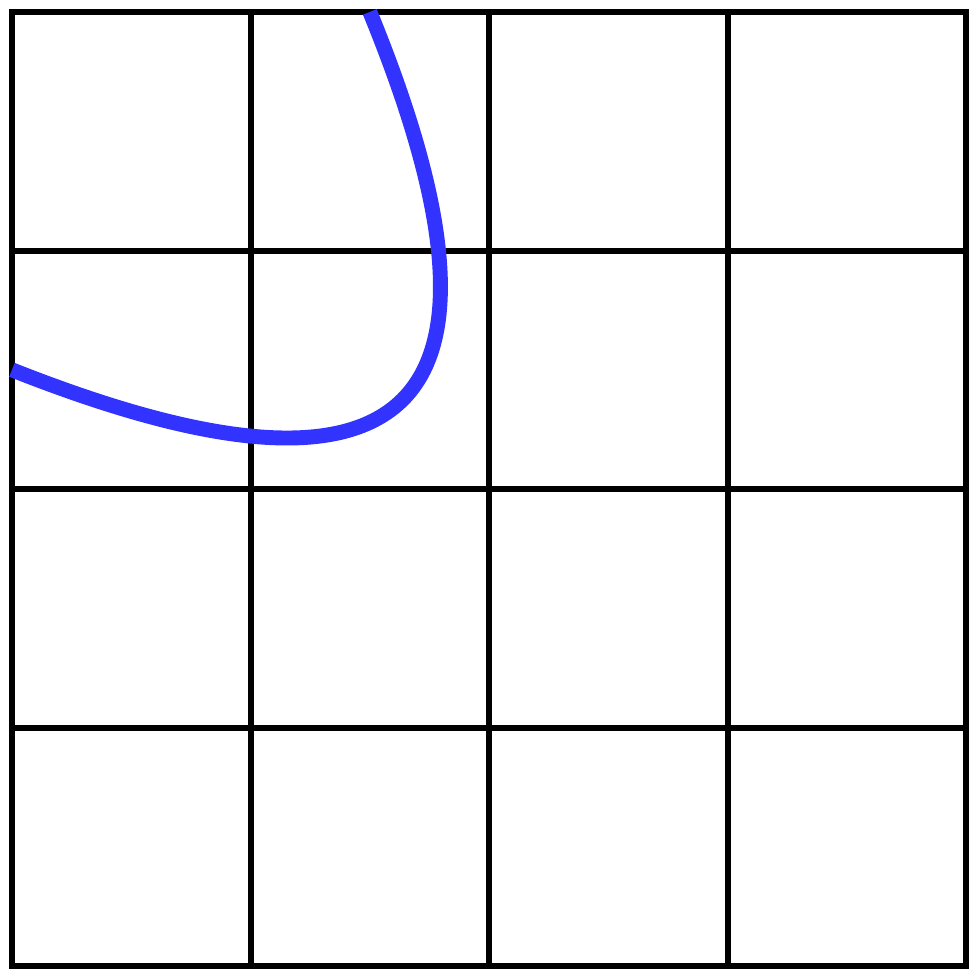}
    \caption{Subdomain partitioning and displacement front.}
    \label{fig:ADDM-curve}
\end{figure}

Assume that the computational domain, denoted as $\Omega$, is subdivided into $N$ subdomains, denoted as $\Omega_i$, such that $\Omega = \bigcup_{i=1}^{N} \Omega_i$ and $\Omega_i \cap \Omega_j = \varnothing$ for all $i \neq j$.
We begin by identifying the location of the advancing front within each subdomain.
At this stage, the subdomains are extended to determine whether the advancing front crosses multiple subdomains, which is then used to assess the coupling strength between them.
The location of the advancing front is determined by the threshold $c_{_S}$ for saturation change in grid cells between adjacent time steps:
\[
M_{k,l} = \left\{ i \mid \exists \, j \; \text{s.t.} \; |\Delta S_{j,i}^{(n)}| > c_{_S}, \, j = 1, \dots, n_{p}, \, i \in \Omega_{k,l} \right\},
\]
where \( |\Delta S_{j,i}^{(n)}| \) represents the saturation change in the \( j \)-th phase within grid cell \( i \) between the \( n \)-th and \( (n-1) \)-th time steps.
 \(\Omega_{k,l}\) represents the set of subdomains obtained by expanding \(\Omega_{k}\) outward by \(l\) layers.
We have \(\Omega_{k} = \Omega_{k,0}\), and let \(M_{k} = M_{k,0}\). Next, we construct the coupling strength graph $G(\Omega)$ between subdomains using $M_{k,l}$. Figure~\ref{fig:ADDM-curve} illustrates an example in which the computational domain is divided into \(4 \times 4\) subdomains. 
The blue curve marks the grid cells identified by the saturation change threshold, and this curve spans multiple subdomains.

We present three different strategies for adaptive coupling of subdomains:
\begin{enumerate}
    \item \textbf{Based on saturation changes at subdomain boundaries (Strategy A)}: For adjacent subdomains \(\Omega_{i}\) and \(\Omega_{j}\), if grid cells near their boundary exhibit saturation changes exceeding $c_{_S}$, then \(\Omega_{i}\) and \(\Omega_{j}\) are coupled:
   \[
   E(\Omega_{i}, \Omega_{j}) \in G({\Omega}), \quad \text{if and only if} \quad M_{i,l} \cap M_{j,l} \neq \varnothing.
   \]
   The coupling pattern is determined by calculating the connected components of the graph \(G({\Omega})\). Figure~\ref{fig:ch3:ADDM01} shows an example for \(l=1\), where green subdomains are coupled for solving and gray subdomains are solved independently. This strategy focuses on identifying critical couplings between subdomains to accurately capture the advancing front with a minimal number of subdomains.
   
    \begin{figure}[htpb]
    \centering
    \includegraphics[width=0.85\textwidth]{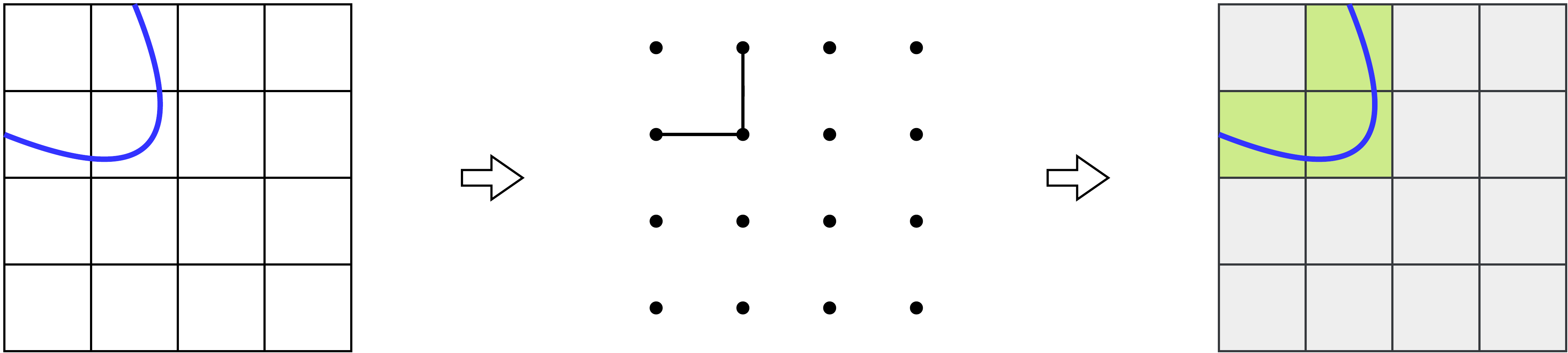}
    \caption{Illustration of the first subdomain coupling strategy (Strategy A).}
    \label{fig:ch3:ADDM01}
    \end{figure}
    
    \item \textbf{Based on active subdomains (Strategy B):} A subdomain \( \Omega_i \) is considered active if there are grid cells within it where the saturation change exceeds \( c_{_S} \). In this case, \( \Omega_i \) is coupled with all its neighboring subdomains:
   \[
   E(\Omega_{i}, \Omega_{j}) \in G({\Omega}), \; \text{for all neighboring subdomains } \,\Omega_{j}, \quad \text{if } M_{i} \neq \varnothing.
   \]
   The coupling pattern is determined by calculating the connected components of the graph \(G({\Omega})\). As shown in Figure~\ref{fig:ch3:ADDM02}, the green subdomains are coupled for solving, while the gray subdomains are solved independently. This approach accounts for both significant couplings between subdomains and the movement of the advancing front, thereby enhancing robustness, although it results in a larger number of coupled subdomains.

    \begin{figure}[!htbp]
    \centering
    \includegraphics[width=0.85\textwidth]{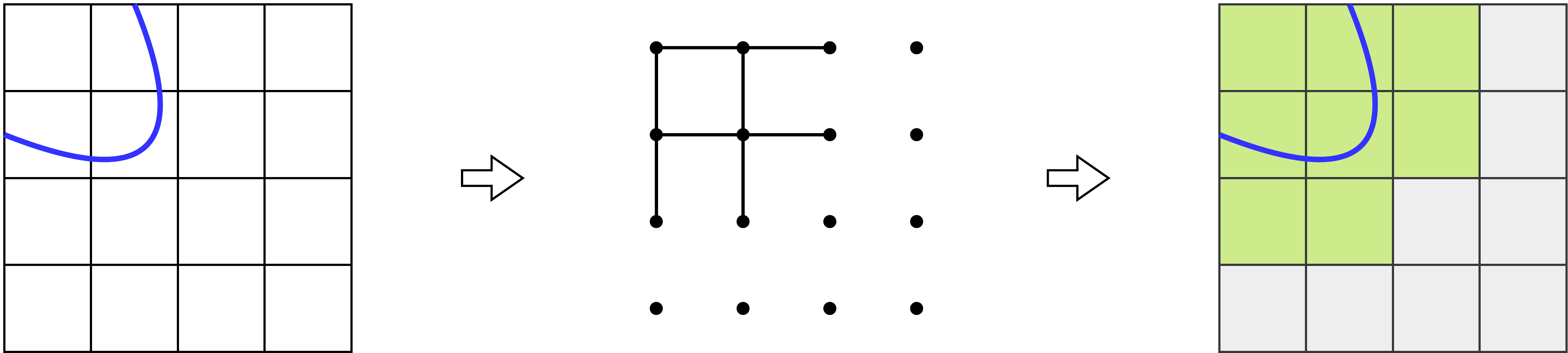}
    \caption{Illustration of the second subdomain coupling strategy (Strategy B).}
    \label{fig:ch3:ADDM02}
    \end{figure}
    
    \item \textbf{Based on weighted graph of subdomains (Strategy C)}: For adjacent subdomains \(\Omega_{i}\) and \(\Omega_{j}\), we define their connection weight as
   \[
   W(\Omega_{i}, \Omega_{j}) = \left(\sum_{k \in M_{i}} + \sum_{k \in M_{j}} + \sum_{k \in M_{i,l} \cap M_{j,l}}\right) \new{\sum_{\iota=1}^{n_{p}} \left|\Delta S_{\iota,k}\right|},\text{for}\, M_{i,l} \cap M_{j,l} \neq \varnothing,
   \]
   which takes into account both its own activity and the coupling strength with its neighboring subdomains.
   If \(W(\Omega_{i}, \Omega_{j}) = 0\), a very small value is assigned to represent the connection between the two subdomains. The weighted graph \(G(\Omega)\) is then partitioned into a specified number of blocks, as shown in Figure~\ref{fig:ch3:ADDM03}. Subdomains marked with the same color are coupled for solving. 
   This strategy provides a comprehensive consideration of couplings between subdomains, effectively utilizing parallel computing resources by coupling as many subdomains as possible.  
   Additionally, it helps prevent the formation of excessively large coupled subdomains by maintaining a controlled number of partitions.
   However, it is less flexible in capturing irregular shapes and sacrifices some recognition capability in critical regions. The number of subgraphs is determined by the predefined maximum number of coupled subdomains, \( N_{\text{cs}} \).

    \begin{figure}[!htbp]
    \centering
    \includegraphics[width=0.85\textwidth]{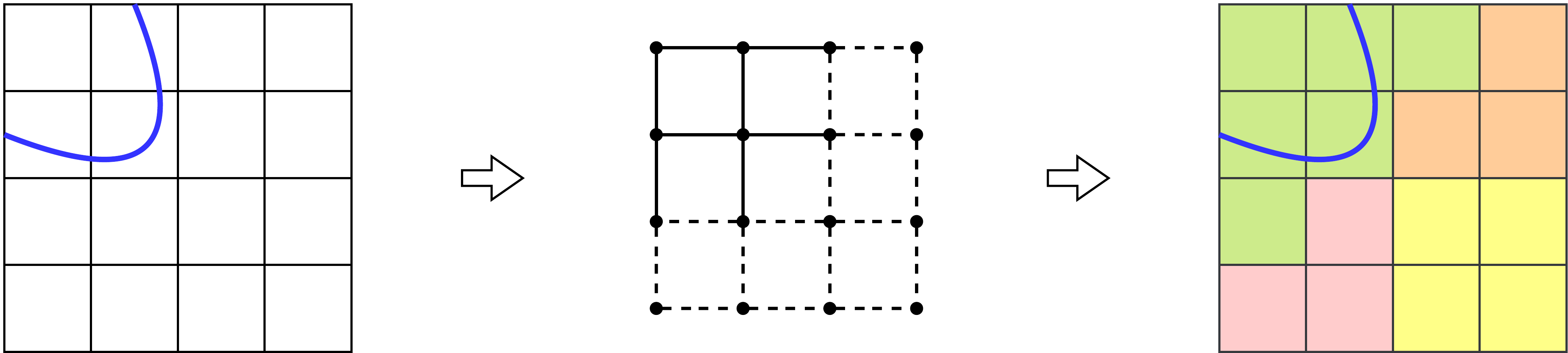}
    \caption{Illustration of the third subdomain coupling strategy (Strategy C).}
    \label{fig:ch3:ADDM03}
    \end{figure}
\end{enumerate}

The Boost Graph Library (BGL) toolset~\citep{siek_boost_2002} can be used to compute the connected components of a graph for Strategies A and B. 
BGL is a highly efficient C++ library tailored for graph-related problems, offering a comprehensive range of data structures and algorithms, including graph traversal, shortest path computation, minimum spanning tree construction, and connected components identification. 
For partitioning weighted graphs in strategy C, the Metis toolset~\citep{karypis_metis_2009} is highly suited. Metis is a robust and efficient library specifically designed for large-scale graphs partitioning, offering advanced algorithms that enable rapid multilevel graph partitioning. Additionally, It is important to note that the computational cost of subdomain coupling partitioning is typically very low, as its computational complexity is related to the initial number of subdomains, \(N\), which is usually small.

\subsection{Subdomain boundary conditions}\label{sec:ADDM-subsec2}
Following the coupling of the subdomains, the system retains the structure of the domain decomposition, with some subdomains being formed by the coupling of multiple initial subdomains, which are then solved in parallel by the corresponding processes.
In this context, the boundary conditions for these subdomains must be specified, analogous to conventional domain decomposition methods.
Typically, the boundary conditions between subdomains are specified using Dirichlet conditions for constant pressure or constant flow velocity, with these values generally taken from the previous time iteration or nonlinear iteration step.

\subsection{The solution framework of ADDM}
In this section, we propose an adaptively coupled domain decomposition method solution framework. 
Without loss of generality, we assume that at time $t$, the domain decomposition after the coupling of initial subdomains is \({\Omega}=\bigcup_{i=1}^{{N}^t}{\Omega}_{i}^t\), where \({\Omega}_{i}^t \cap {\Omega}_{j}^t = \varnothing\) (for \(\forall \, i \neq j \)). 
Here, ${N}^t$ denotes the number of subdomains after coupling, \({\Omega}_{i}^t = \bigcup_{j \in \mathcal{N}_{i}^t }\Omega_{j}\), where \(\mathcal{N}_{i}^t\) is the index set of the subdomains that form \({\Omega}_{i}^t\), satisfying \(\bigcup_{i=1}^{N^t} \mathcal{N}_{i}^t = \left\{1, 2, \dots, N\right\}\) and \(\mathcal{N}_{i}^t \cap \mathcal{N}_{j}^t = \varnothing\) (for \(\forall \, i \neq j \)). 
At time \( t \), if all initial subdomains are coupled and solved together, the method results in a fully coupled algorithm with \({N}^t = 1\) and \({\Omega}_{1}^t=\Omega\). On the other hand, if each subdomain is solved independently, the method reduces to the classical domain decomposition method (denoted as CDDM in our context) with \({N}^t=N\) and \({\Omega}_{i}^t=\Omega_{i}\).

The key feature of this solution framework lies in its ability to dynamically identify and address significant coupling relationships between subdomains during the simulation, thereby providing a foundation for accelerating the solution.
After the subdomains are coupled, the system still maintains a domain decomposition structure. 
Therefore, strategies from classical domain decomposition methods, such as overlapping subdomains and multilevel domain decomposition~\citep{hwang2007class}, are also be applied within this framework. 
Methods that already incorporate domain decomposition strategies can also replace their domain decomposition components with adaptively coupled domain decomposition to handle local strong nonlinearities, thereby enhancing convergence performance.
For example, in the classical ASPIN algorithm, the original domain decomposition strategies \(\left\{\Omega_{i}\right\}_{i=1}^{N}\) can be replaced by \(\left\{{\Omega}_{i}^t\right\}_{i=1}^{{N}^t}\) at the beginning of each time step, while the remaining steps remain unchanged.

We propose an efficient algorithm within the adaptively coupled domain decomposition solution framework in the context of parallel computing, where the solution to the local nonlinear problem defined on \({\Omega}_{i}^t\) is used as the initial guess for the global nonlinear problem. As an example, in this algorithm, \new{the standard Newton--Krylov method} is employed to solve both the global and local nonlinear problems. For ease of comparison, we integrate workflows of multiple algorithms into a single process, as detailed in Algorithm~\ref{alg:ADDM-InitialGuess}.

\begin{algorithm}[htpb]
  \caption{Initial Guess with ADDM for Nonlinear Problems}\label{alg:ADDM-InitialGuess}
  \setcounter{AlgoLine}{0}
  \LinesNumbered
  
  \KwIn{$\mathbf{X}^{(n)}$, $\{{\Omega}_i\}_{i=1}^{N}$, \emph{asm\_flag},  $\epsilon_{\text{global-init}} $, $\epsilon_{\text{local}}$, $MaxIt$;}

	Determine \(\left\{{\Omega}_{i}^t\right\}_{i=1}^{{N}^t}\) using adaptive coupling strategies (see Section~\ref{sec:ADDM-subsec1});
    
	Set $\mathbf{X}^{1}=\mathbf{X}^{(n)}$, \(\mathbf{X}_i^{1}=\left.\mathbf{X}^{(n)}\right|_{{\Omega}_{i}^t}\), for $\forall~i=1,\ldots,N^t$;

	\For(\emph{\small // Solve local problem $\mathcal{F}_i (\mathbf{X}) = 0$ with initial guess $\mathbf{X}_i^{1}$})
	{$k=1,\ldots,MaxIt$}
	{
		\textbf{Parallel} \For(\emph{\small // Solve subproblems in parallel}){$i=1,\ldots,N^t$}
		{
			Compute local Jacobian matrix $\mathbf{J}_i = \mathcal{F}'_i (\mathbf{X}_i^{k})$;
			
			Solve local linear system $\mathbf{J}_i \mathbf{d}_i = - \mathbf{F}_i $;
			
			Compute local step size $\alpha_i$ obtained by backtracking line search;
			
			Update $\mathbf{X}_i^{k+1} = \mathbf{X}_i^{k} + \alpha_i \mathbf{d}_i$;
			
			\uIf(\emph{\small // ASM\_ADDM}){asm\_flag $==$ true}
			{
				Exchange information of boundary elements between subdomains;

				\lIf{$\| \mathcal{F} (\mathbf{X}^{k+1}) \| \leq \epsilon_{\text{global-init}} \| \mathcal{F} (\mathbf{X}^{1}) \|$}{\textbf{break}}
			}
			\Else(\emph{\small // ADDM}){

				\lIf{$\| \mathcal{F}_i (\mathbf{X}_i^{k+1}) \| \leq \epsilon_{\text{local}} \| \mathcal{F}_i (\mathbf{X}_i^{1}) \|$}{\textbf{break}}
			}
		}

	}

    Set \(\left.\mathbf{\bar{X}}\right|_{{\Omega}_{i}^t}=\mathbf{X}_{i}^{k+1}\), for $\forall~i=1,\ldots,N^t$;

   \KwOut{$\mathbf{\bar{X}}$.}
\end{algorithm}

\begin{algorithm}[htpb]
  \caption{Solution for Nonlinear Problems with an Initial Guess}\label{alg:ADDM-FIM}
  \setcounter{AlgoLine}{0}
  \LinesNumbered
  
  \KwIn{$\mathbf{X}^{(n)}$, $\{{\Omega}_i\}_{i=1}^{N}$,  \emph{asm\_flag}, $\epsilon_{\text{local}}$, $\epsilon_{\text{global-init}} $,  $\epsilon_{\text{global}}$, $\epsilon$, $MaxIt$;}

   Compute initial guess $\mathbf{\bar{X}}^{1}$ by calling Algorithm~\ref{alg:ADDM-InitialGuess};

	\For(\emph{\small //Solve global problem $\mathcal{F} (\mathbf{X}) = 0$ with initial guess $\mathbf{\bar{X}}^{1}$})
	{$k=1,\ldots,MaxIt$}
	{
			Compute global Jacobian matrix $\mathbf{J} = \mathcal{F}' (\mathbf{\bar{X}}^{k})$;
			
			Solve global linear system $\mathbf{J} \mathbf{d} = - \mathbf{F} $;
			
			Compute global step size $\alpha$ obtained by backtracking line search;
			
			Update $\mathbf{\bar{X}}^{k+1} = \mathbf{\bar{X}}^{k} + \alpha \mathbf{d}$;
			
			\lIf{$\| \mathcal{F} (\mathbf{\bar{X}}^{k+1}) \| \leq \epsilon_{\text{global}} \| \mathcal{F} (\mathbf{\bar{X}}^{1}) \|$ {\rm{or}} $\| \mathbf{\bar{X}}^{k+1} - \mathbf{\bar{X}}^{k}  \| \leq \epsilon $}{\textbf{break}}
	}

	Set \(\mathbf{X}^{(n+1)}=\mathbf{\bar{X}}^{k+1}\);

   \KwOut{$\mathbf{X}^{(n+1)}$.}
\end{algorithm}

\new{The algorithmic framework (see Algorithm~\ref{alg:ADDM-FIM}) follows the Newton--Krylov paradigm~\citep{KNOLL2004357}. Within this framework, a nonlinear iteration is applied at the outer level, while each resulting linearized system is solved at the inner level by a Krylov subspace method, usually equipped with an efficient preconditioner.} Algorithm~\ref{alg:ADDM-FIM} provides flexibility in solving the nonlinear problem through different domain decomposition and coupling strategies, allowing for the implementation of the following four methods:

\begin{itemize} 
    \item \new{\textbf{Standard:} This approach directly addresses the global nonlinear problem using a standard Newton--Krylov method. Specifically, Algorithm~\ref{alg:ADDM-FIM} uses the solution vector from the previous time step, $\mathbf{X}^{(n)}$, as the initial guess, i.e., $\mathbf{\bar{X}}^{1} = \mathbf{X}^{(n)}$.} 

    \item \textbf{ASM\_CDDM:} Algorithm~\ref{alg:ADDM-FIM} utilizes the initial domain decomposition pattern with ASM to provide the initial guess. Specifically, the subdomain adaptive coupling strategy is not applied, and we have $\Tilde{N}(t) \equiv N$ and $\Tilde{\Omega}_{i}(t) \equiv \Omega_{i}$, for $i = 1, 2, \dots, N$. The \emph{asm\_flag} in Algorithm~\ref{alg:ADDM-InitialGuess} is set to be true.

    \item \textbf{ASM\_ADDM:} Algorithm~\ref{alg:ADDM-FIM} utilizes adaptive domain decomposition methods with ASM to provide the initial guess.
    Specifically, the subdomain adaptive coupling strategy is applied, and the \emph{asm\_flag} in Algorithm~\ref{alg:ADDM-InitialGuess} is set to be true.

    \item \textbf{ADDM:} Algorithm~\ref{alg:ADDM-FIM} utilizes adaptive domain decomposition methods without using ASM to provide the initial guess.
    Specifically, the subdomain adaptive coupling strategy is applied, and the \emph{asm\_flag} in Algorithm~\ref{alg:ADDM-InitialGuess} is set to be false.
\end{itemize}

The key difference between ADDM and ASM\_ADDM is that the former does not exchange information between boundary cells during the solution process of a single time step. 
The boundary cell values are fixed at those from the previous time step.
This approach is grounded in one of the core objectives of the subdomain adaptive coupling strategy: identifying the significant coupling relationships between subdomains.
This means that the coupling strength between subdomains in the newly formed domain decomposition is weak, assuming the coupling pattern is ideal.
Therefore, using fixed boundary values from the previous time step is reasonable and potentially enhances the convergence of the local solution process.

\new{
\begin{remark}
Subdomain selection and updating during the adaptive coupling in Algorithm~\ref{alg:ADDM-InitialGuess} are described as follows. 
The initial domain decomposition, obtained using the graph partitioning software ParMETIS~\citep{ParMETIS2020}, is given by $\Omega = \bigcup_{i=1}^{N} \Omega_i$ and $\Omega_i \cap \Omega_j = \varnothing~(i \neq j)$.
During the adaptive coupling process, the subdomains are dynamically updated. At time $t$, the domain decomposition becomes \({\Omega}=\bigcup_{i=1}^{{N}^t}{\Omega}_{i}^t\) and \({\Omega}_{i}^t \cap {\Omega}_{j}^t = \varnothing~(i \neq j)\). Each updated subdomain $\Omega_i^t$ is formed by merging a group of initial subdomains, i.e.,
\({\Omega}_{i}^t = \bigcup_{j \in \mathcal{N}_{i}^t }\Omega_{j}\).
Here, $\mathcal{N}_i^t$ denotes the index set of subdomains associated with $\Omega_i^t$, which is determined by the adaptive coupling strategy described in Section~\ref{sec:ADDM-subsec1} using the BGL or METIS toolsets.
\end{remark}
}

\begin{remark}
The convergence criterion in Algorithm~\ref{alg:ADDM-FIM}, \( \| \mathbf{\bar{X}}^{k+1} - \mathbf{\bar{X}}^{k} \| \leq \epsilon \), specifically refers to the change in the pressure and saturation variables between two consecutive iterations, i.e.,  
\[
\max_{\tau} | P^{k+1}_{\tau} - P^{k}_{\tau} | \leq \epsilon_{P} \quad \text{and} \quad \max_{j,\tau} | S^{k+1}_{j,\tau} - S^{k}_{j,\tau} | \leq \epsilon_{S}, \quad j = 1, \dots, n_p,
\]
where \( \tau \) is the grid cell index.
\end{remark}

\section{Numerical experiments}\label{sec:Numerical Experiments}
In this section, we present numerical results from a series of experiments designed to evaluate the convergence and parallel performance of the newly proposed adaptively coupled domain decomposition method (ADDM). The test cases include simulations in both simple homogeneous and complex heterogeneous media, as well as parallel strong scalability tests involving hundreds of millions of grid points and parallel weak scalability tests with grid sizes refined to the half-billion level.
The choice of time steps is typically crucial to the performance of solution methods. However, in practice, determining an optimal strategy can be challenging. To facilitate method comparison and improve the significance of our experimental results, we adopted the following approach: multiple experiments were conducted to identify the optimal time steps for \new{the standard Newton--Krylov method}, and this configuration was subsequently applied to all other methods.

The proposed method is implemented in our open-source parallel reservoir simulator, OpenCAEPoro\footnote{\url{https://github.com/OpenCAEPlus/OpenCAEPoroX}}~\citep{OCP2024}. 
For linear solvers, we used the constrained pressure residual (CPR) preconditioned iterative method~\citep{Wallis1983,Wallis1985}, implemented with the portable, extensible toolkit for scientific computation (PETSc)~\citep{petsc_web_page2024} and Hypre~\citep{Hypre2002} libraries. 
In the CPR preconditioner, the first stage employs Hypre's Boomer-AMG to solve the pressure subsystem, while the second stage uses PETSc's Block-Jacobi with BILU(0) to solve the overall system. 
Additionally, the iterative method employed is the flexible generalized minimal residual method (FGMRES)~\citep{Saad2003}.

The parameters and experimental setup for the numerical experiments are as follows. The convergence criteria are set to \( \epsilon_{\text{local}} = 10^{-2} \) (or \( \epsilon_{\text{global-init}} = 10^{-2} \)), \( \epsilon_{\text{global}} = 10^{-4} \), \( \epsilon_{P} = 1 \), and \( \epsilon_{S} = 1 \times 10^{-3} \). The maximum number of Newton iterations for both local and global nonlinear problems is limited to 10, while the maximum number of linear iterations for both local and global linear problems is capped at 50. During the Newton–Raphson iteration, the choice of the Newton step size critically influences the convergence behavior of nonlinear problems. Instabilities—and even convergence failures—can arise when substantial phase-state changes occur within elements. In this study, we adopt the Appleyard chopping strategy~\citep{Appleyard1983}, a widely used step-size control method that performs well across most scenarios and is implemented in many simulators~\citep{eclipse_manual_2021,rasmussen2021open}. Specifically, within each Newton iteration, we cap the maximum change in saturation at 0.2 to enhance stability and convergence.

Numerical experiments are conducted on a supercomputer, with each compute node equipped with two Intel 6458Q CPUs, each CPU having 32 cores running at 3.1 GHz, and 256 GB of memory.

\subsection{Case 1}\label{subsec:Case 1} 
This case is an extended version of the SPE1 benchmark~\citep{odeh_comparison_1981}. 
Specifically, the original model is refined with a higher grid resolution and extended in the horizontal direction. 
The final grid size is \(1232 \times 1232 \times 10\), with each grid cell measuring \(20 \text{ ft} \times 20 \text{ ft} \times 10 \text{ ft}\). 
The well layout is modified to a five-spot pattern, with four injection wells located at the corners and one production well at the center. 
The injection wells are perforated in the top two layers, with a target injection rate of 20,000 thousand standard cubic feet per day (Mscf/day). The production well is perforated in the bottom five layers, with a target oil production rate of 20,000 standard barrels per day (stb/day).
\begin{figure}[htpb]
\centering
\includegraphics[width=0.45\linewidth]{Partition.pdf}
\caption{The domain partitioning by 784 processes using ParMetis, with colors corresponding to the rank of each process.}
\label{fig:case1:partition}
\end{figure}

The total simulation time is 3000 days, with the injection process divided into three phases: pre-gas breakthrough at the production well, gas breakthrough at the production well, and post-gas breakthrough at the production well.
The second phase is relatively short (approximately from day 2850 to day 2900) but presents the greatest computational challenges. 
As such, this case provides an excellent basis for comparing different methods. 
In the following tests, 784 MPI processes are used. 
Figure~\ref{fig:case1:partition} presents the domain partitioning pattern. Notably, due to the significantly larger number of grid cells in the horizontal direction compared to the vertical direction, the grid partitioning achieved with ParMetis is effectively two-dimensional.

\subsubsection{Correctness verification and performance comparison}
To validate the correctness and performance of the proposed methods, we compare four approaches: Standard, ASM\_CDDM, ASM\_ADDM, and ADDM. In both ASM\_ADDM and ADDM, Strategy B is chosen for subdomain adaptive coupling, with \( c_{s} \) set to \( 5 \times 10^{-3} \). Additionally, a constant pressure condition is applied as the boundary condition in ASM\_CDDM, ASM\_ADDM, and ADDM (For ASM\_CDDM and ASM\_ADDM, the pressure values at the neighboring boundary grids are taken from the previous nonlinear iteration step. For ADDM, the pressure values at the neighboring boundary grids are taken from the previous time step).

\begin{figure}[htpb]
	\centering
	\subfigure[\texttt{FPR}]{\label{fig:ADDM:case1:3000Day:FPR}  
		\begin{minipage}[t]{0.31\linewidth}
			\centering
			\includegraphics[trim=0.1cm 0.3cm 0.1cm 0.1cm,clip,width=4.9cm]{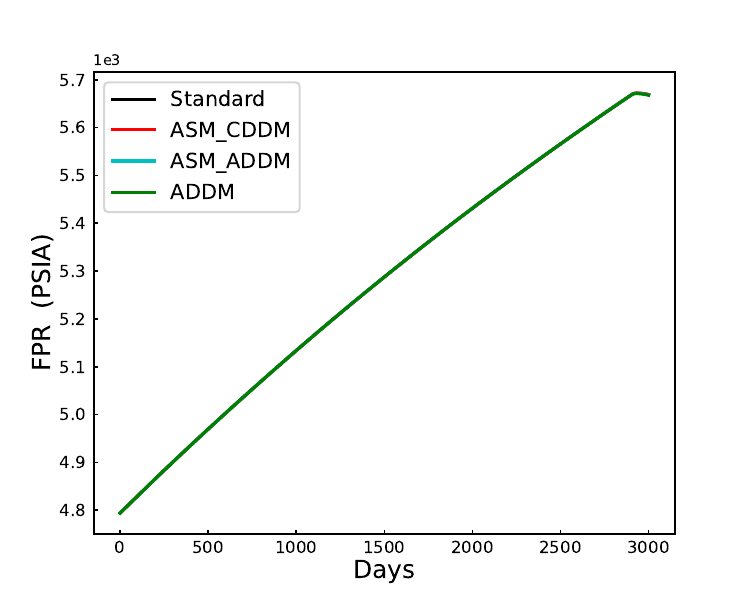}  
		\end{minipage}
	}%
	\subfigure[\texttt{FGPR}]{\label{fig:ADDM:case1:3000Day:FGPR}  
		\begin{minipage}[t]{0.31\linewidth}
			\centering
			\includegraphics[trim=0.1cm 0.3cm 0.1cm 0.1cm,clip,width=4.9cm]{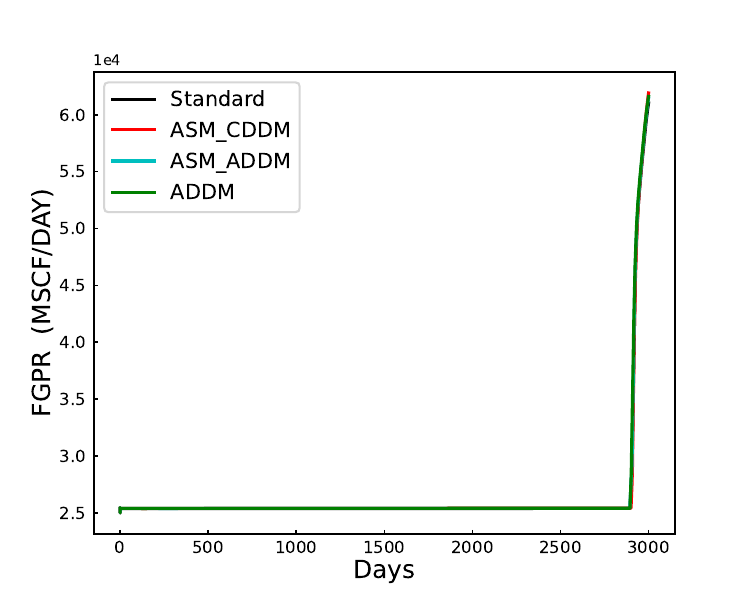}  
		\end{minipage}
	}%
	\subfigure[\texttt{FWPR}]{\label{fig:ADDM:case1:3000Day:FWPR} 
		\begin{minipage}[t]{0.31\linewidth}
			\centering
			\includegraphics[trim=0.1cm 0.3cm 0.1cm 0.1cm,clip,width=4.9cm]{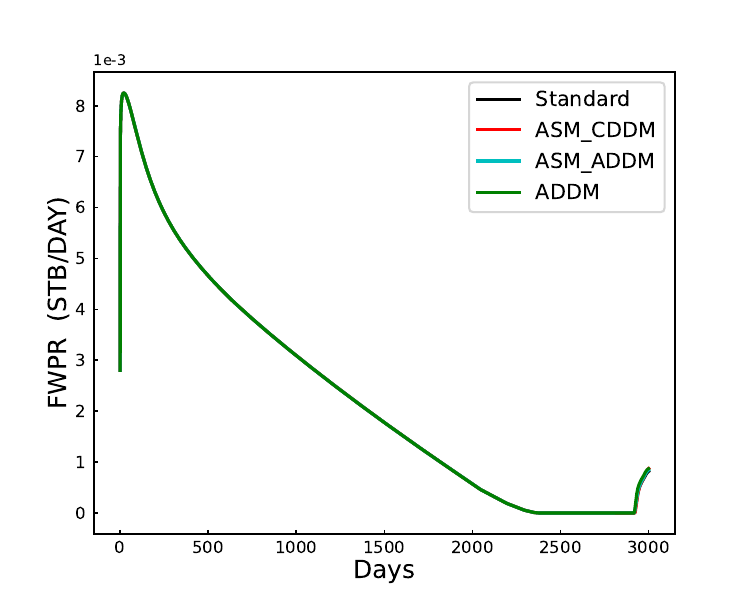}  
		\end{minipage}
	}%
    \caption{Comparisons of \texttt{FPR}, \texttt{FGPR}, and \texttt{FWPR} among four methods for Case 1.}
    \label{fig:case1:3000Day-1}  
\end{figure}

\begin{figure}[htpb]
	\centering
	\subfigure[\texttt{NRiter}]{\label{fig:ADDM:case1:3000Day:NRiter}
		\begin{minipage}[t]{0.31\linewidth}
			\centering
			\includegraphics[trim=0.1cm 0.3cm 0.1cm 0.1cm,clip,width=4.9cm]{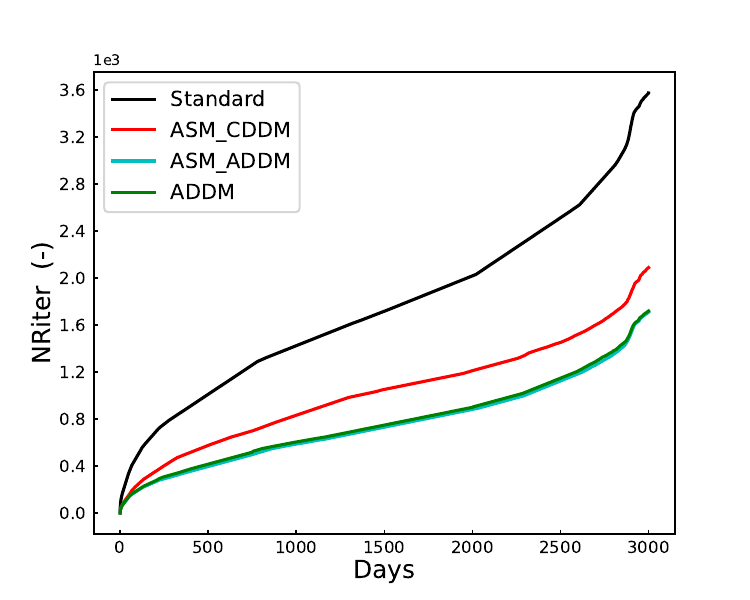}  
		\end{minipage}
	}%
	\subfigure[\texttt{LSiter}]{\label{fig:ADDM:case1:3000Day:LSiter}  
		\begin{minipage}[t]{0.31\linewidth}
			\centering
			\includegraphics[trim=0.1cm 0.3cm 0.1cm 0.1cm,clip,width=4.9cm]{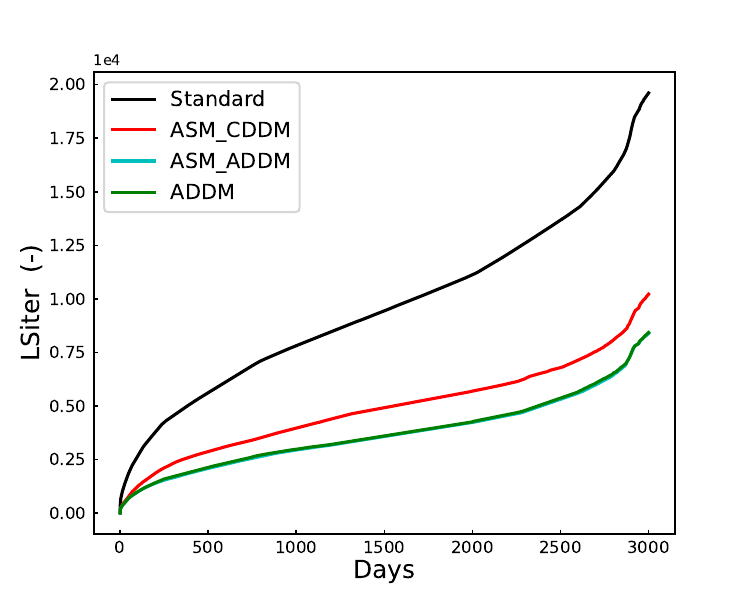}  
		\end{minipage}
	}%
	\subfigure[\texttt{Runtime}]{\label{fig:ADDM:case1:3000Day:Runtime}  
		\begin{minipage}[t]{0.31\linewidth}
			\centering
			\includegraphics[trim=0.1cm 0.3cm 0.1cm 0.1cm,clip,width=4.9cm]{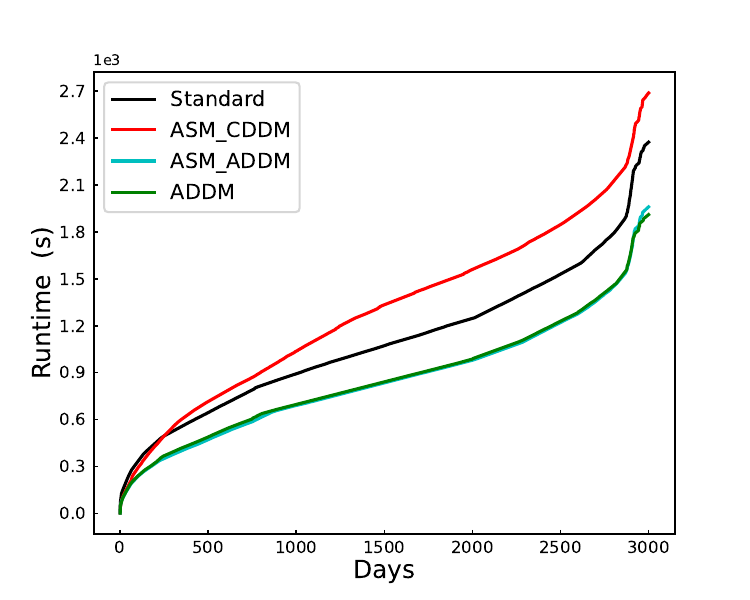} 
		\end{minipage}
	}%
    \caption{Comparisons of \texttt{NRiter}, \texttt{LSiter}, and \texttt{Runtime} among four methods for Case 1.}
    \label{fig:case1:3000Day-1-2}  
\end{figure}

Figure~\ref{fig:case1:3000Day-1} shows the field average pressure (\texttt{FPR}), field gas production rate (\texttt{FGPR}), and field water production rate (\texttt{FWPR}), validating the consistency of the computational results across the different solution methods. Figure~\ref{fig:case1:3000Day-1-2} presents the cumulative global Newton--Raphson iterations (\texttt{NRiter}), cumulative global linear iterations (\texttt{LSiter}), and total simulation runtime (\texttt{Runtime}). Table~\ref{tab:case1:3000day-1} provides detailed performance comparisons, including the number of time steps (\texttt{Timestep}) and the number of local Newton--Raphson iterations required in the initial solution process (\texttt{NRiter(DDM)}).

From Figure~\ref{fig:case1:3000Day-1-2}, it can be seen that, compared to Standard, the other three methods significantly reduce the number of global Newton--Raphson iterations, thereby decreasing the number of global linear iterations. However, in terms of total runtime, ASM\_CDDM, which does not employ the subdomain adaptive coupling strategy, is slower than Standard. 
This is because the reduction in global Newton--Raphson iterations is insufficient to compensate for the additional cost of computing initial solutions. In contrast, ASM\_ADDM and ADDM, which incorporate the subdomain adaptive coupling strategy, demonstrate runtime advantages by providing higher-quality initial solution approximations at a significant lower computational cost.

\begin{table} [htpb]
    \caption{Performance comparison results of four methods for Case 1. The numbers in parentheses represent the iterations wasted due to solver failures. Bold indicates the best performance results.}  
    \label{tab:case1:3000day-1}  
    \centering
    \small  
    \setlength{\tabcolsep} {2.5pt}  
    \renewcommand{\arraystretch} {1.2}  
	\begin{tabular} {cccccccc}  
        \hline
        \texttt{Method} & \texttt{Timestep} & \texttt{NRiter} & \texttt{LSiter} &
        \texttt{NRiter(DDM)} &
        \texttt{Runtime}(s) \\
        \hline
        Standard       & 725 & 3573(+315) & 19608(+2939) & 0 & 2373 \\
        ASM\_CDDM & 855 & 2088(+160) & 10221(+1064) & 3941(+505)  & 2688 \\
        ASM\_ADDM & 739 & 1708(+170) & 8386(+1213) & 2494(+152)  & 1960 \\
        ADDM      & 732 & 1719(+112) & 8431(+894) & 2461(+135)  & \textbf{1910} \\
        \hline
    \end{tabular} 
\end{table}

As shown in Table~\ref{tab:case1:3000day-1}, compared to Standard, ASM\_CDDM reduces global Newton--Raphson iterations by 42.2\% and global linear iterations by 50.0\%; however, it increases runtime by 13.3\%. ASM\_ADDM achieves a 51.7\% reduction in global Newton--Raphson iterations and a 57.4\% reduction in global linear iterations, along with a 17.4\% decrease in runtime. ADDM results in a 52.9\% reduction in global Newton--Raphson iterations and a 58.6\% reduction in global linear iterations, leading to a 19.5\% decrease in runtime. 
Additionally, as indicated in Table~\ref{tab:case1:3000day-1}, the use of initial value techniques effectively identifies potential error-prone steps early, thereby reducing computational costs. Furthermore, the application of subdomain adaptive coupling techniques significantly accelerates the convergence of the initial value problem. For instance, compared to ASM\_CDDM, ASM\_ADDM reduces the number of Newton--Raphson iterations by 40.5\% during the initial solution process.

\subsubsection{Coupling strategies and boundary conditions}  
\begin{figure}[htpb]
\subfigtopskip = -1pt
	\centering
	\subfigure[750d]{\label{fig:ADDM:case1:750d-sags}  
		\begin{minipage}[t]{0.23\linewidth}
			\centering
			\includegraphics[width=3.3cm]{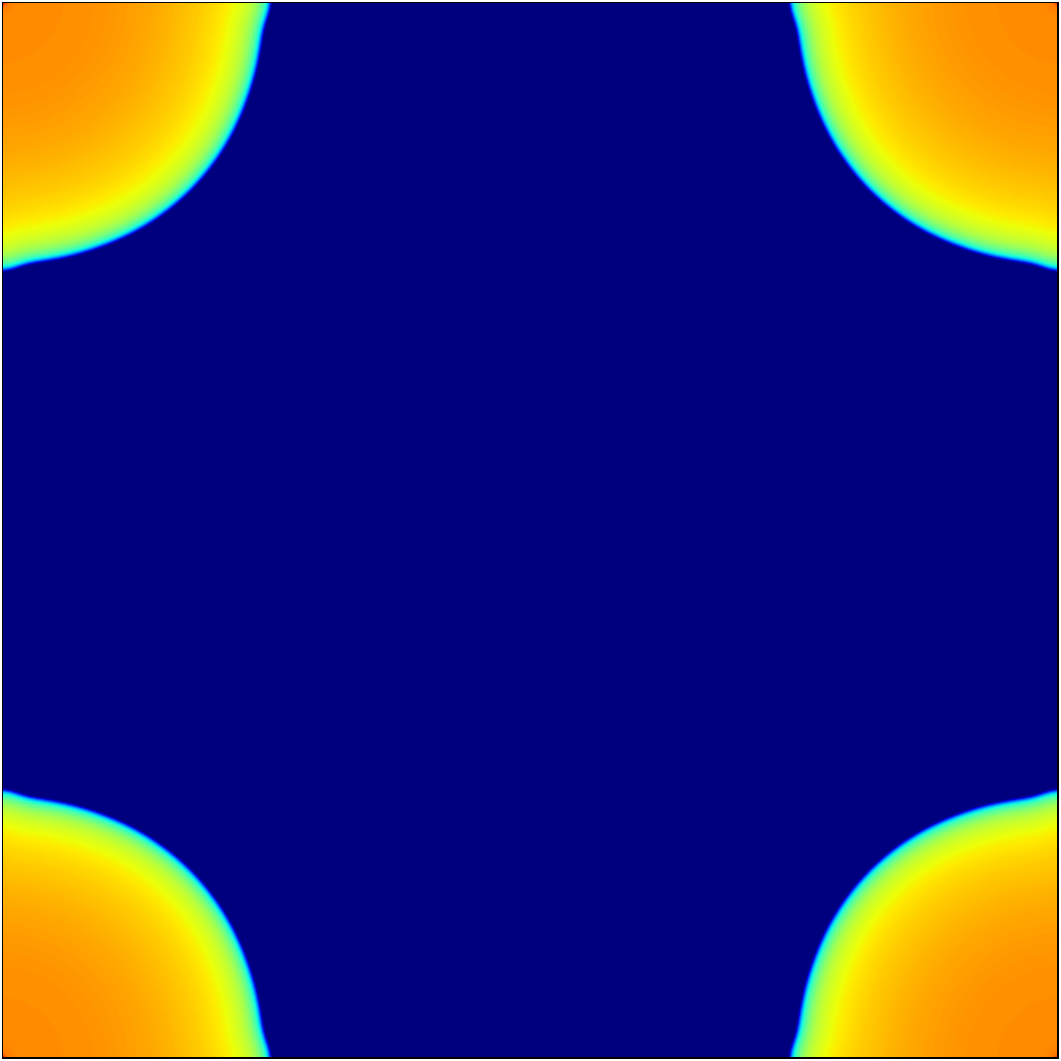}  
		\end{minipage}
	}%
	\subfigure[2000d]{\label{fig:ADDM:case1:2000d-sags} 
		\begin{minipage}[t]{0.23\linewidth}
			\centering
			\includegraphics[width=3.3cm]{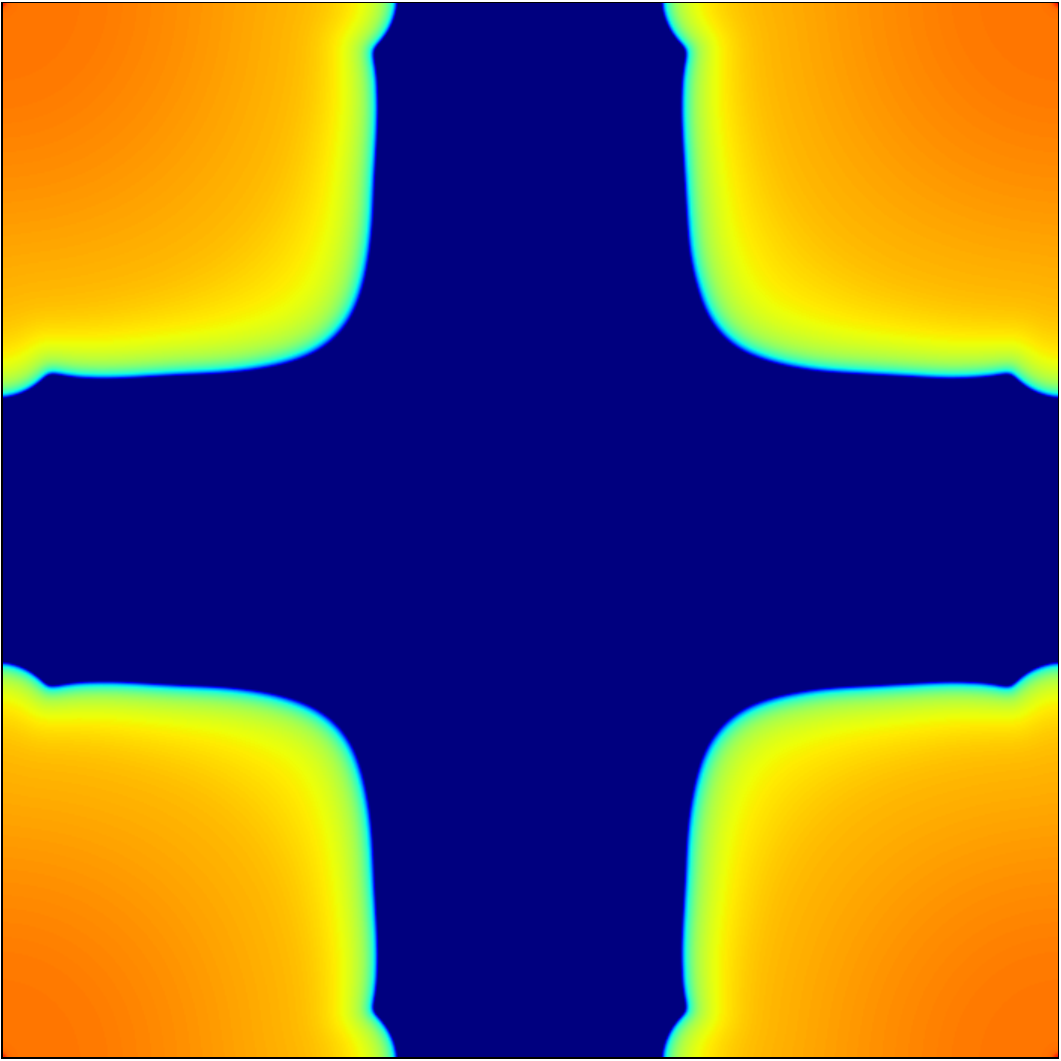}  
		\end{minipage}
	}%
	\subfigure[2820d]{\label{fig:ADDM:case1:2820d-sags} 
		\begin{minipage}[t]{0.23\linewidth}
			\centering
			\includegraphics[width=3.3cm]{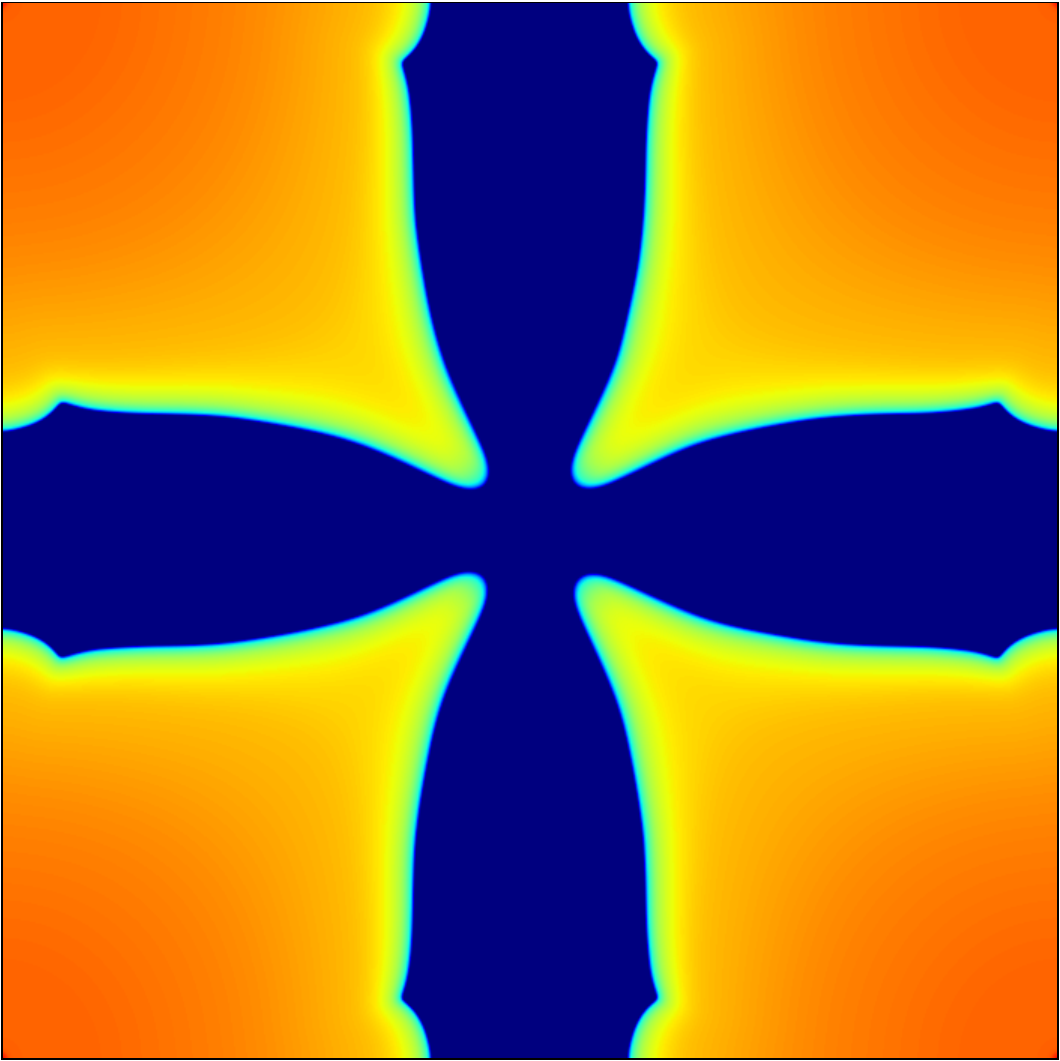}  
		\end{minipage}
	}%
	\subfigure[3000d]{\label{fig:ADDM:case1:3000d-sags} 
		\begin{minipage}[t]{0.23\linewidth}
			\centering
			\includegraphics[width=3.3cm]{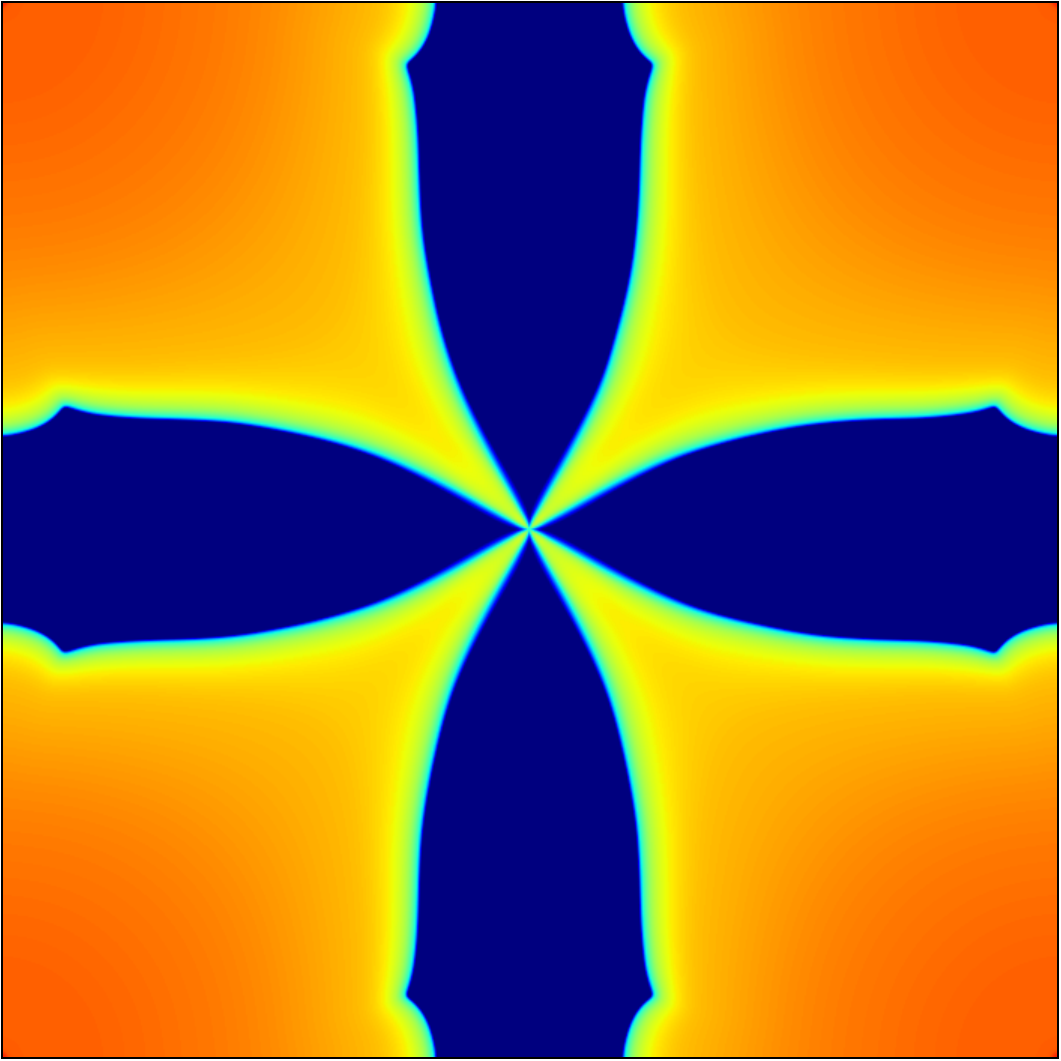}  
		\end{minipage}
	}%

	\subfigure[750d]{\label{fig:ADDM:case1:750d-cs01}  
		\begin{minipage}[t]{0.23\linewidth}
			\centering
			\includegraphics[width=3.3cm]{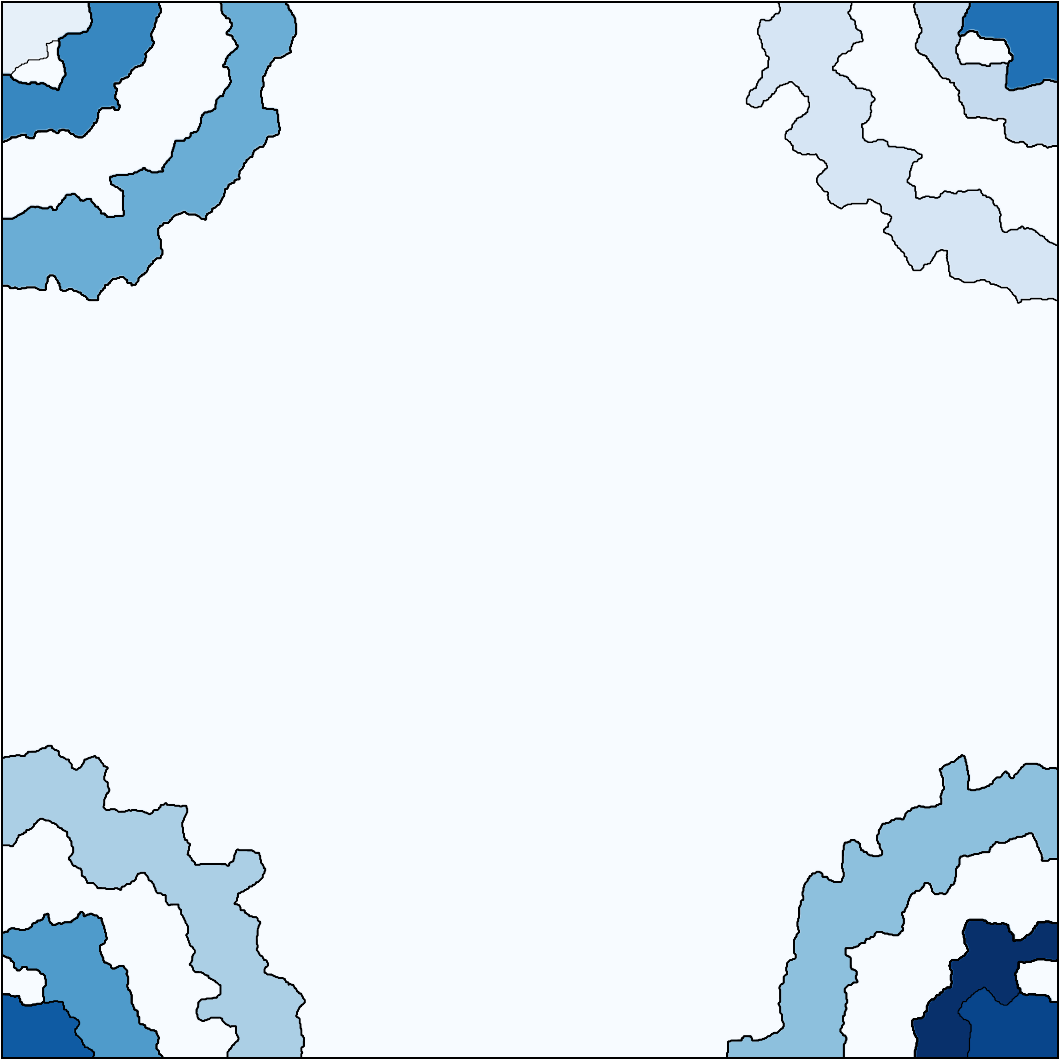}  
		\end{minipage}
	}%
	\subfigure[2000d]{\label{fig:ADDM:case1:2000d-cs01} 
		\begin{minipage}[t]{0.23\linewidth}
			\centering
			\includegraphics[width=3.3cm]{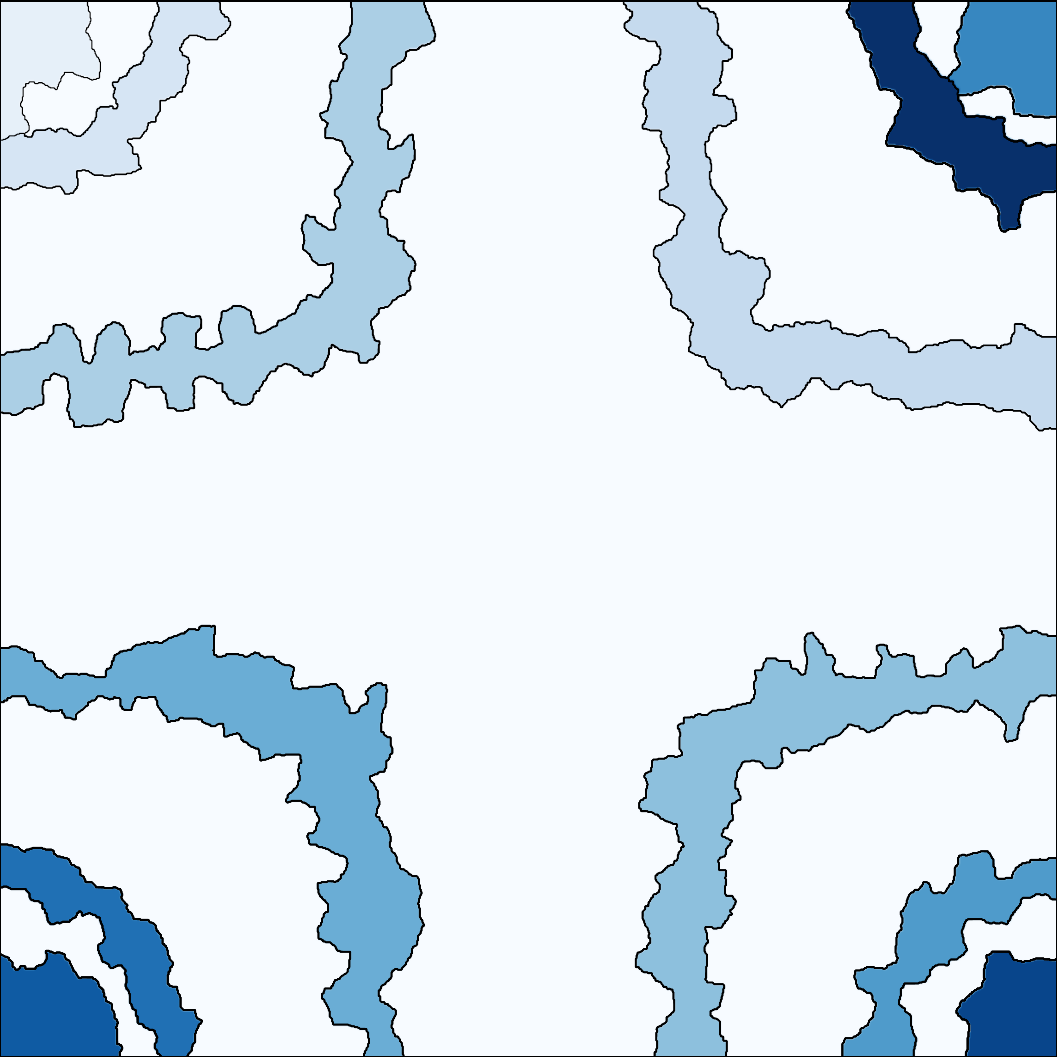}  
		\end{minipage}
	}%
	\subfigure[2820d]{\label{fig:ADDM:case1:2820d-cs01} 
		\begin{minipage}[t]{0.23\linewidth}
			\centering
			\includegraphics[width=3.3cm]{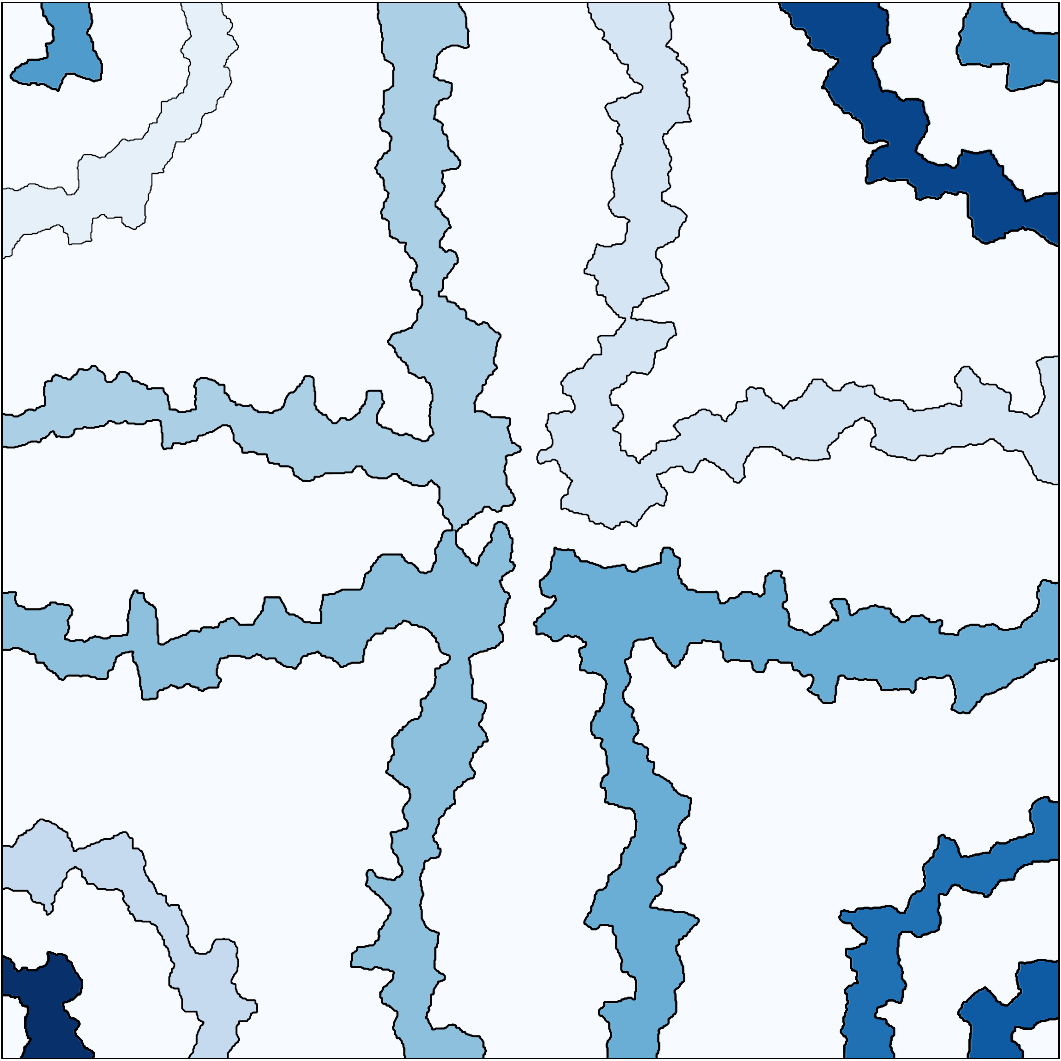}  
		\end{minipage}
	}%
	\subfigure[3000d]{\label{fig:ADDM:case1:3000d-cs01} 
		\begin{minipage}[t]{0.23\linewidth}
			\centering
			\includegraphics[width=3.3cm]{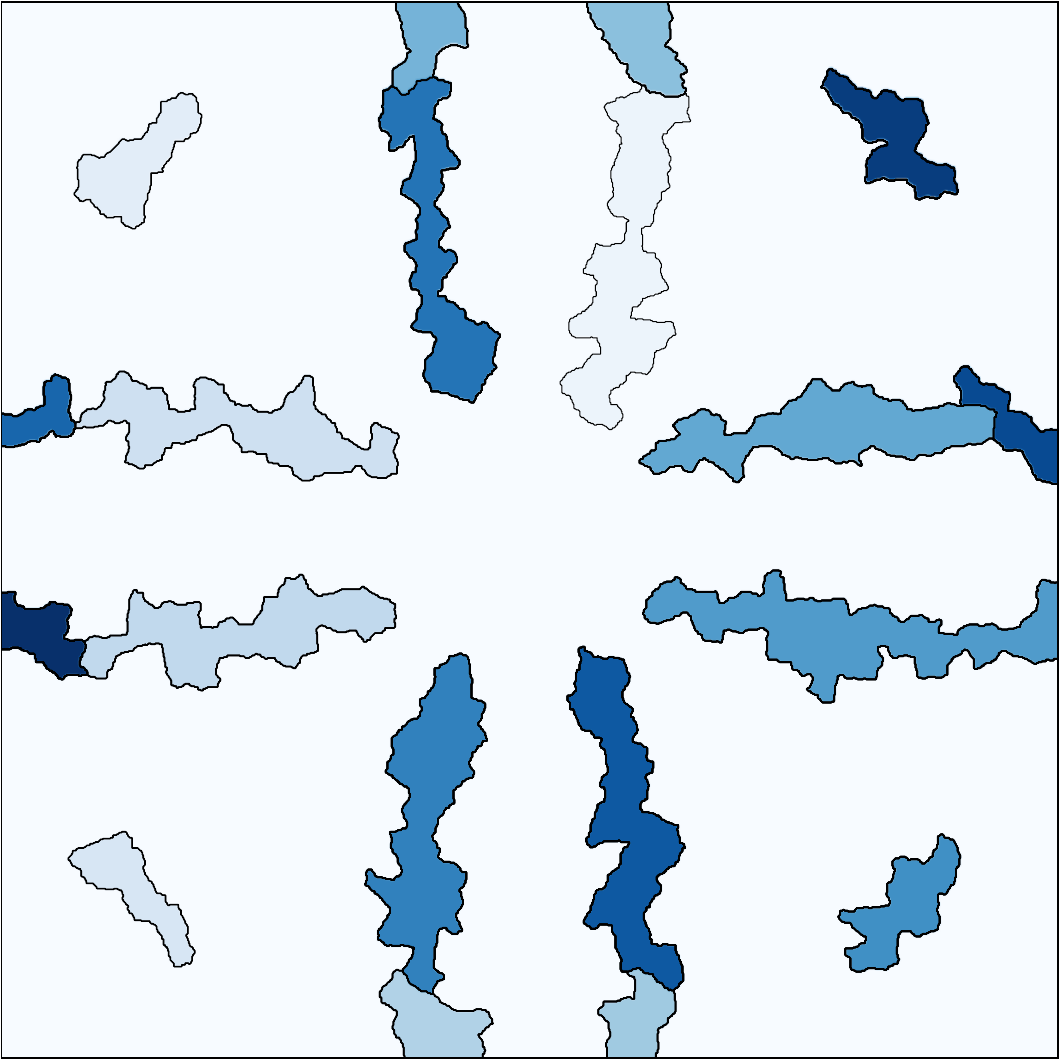}  
		\end{minipage}
	}%

	\subfigure[750d]{\label{fig:ADDM:case1:750d-cs02}  
		\begin{minipage}[t]{0.23\linewidth}
			\centering
			\includegraphics[width=3.3cm]{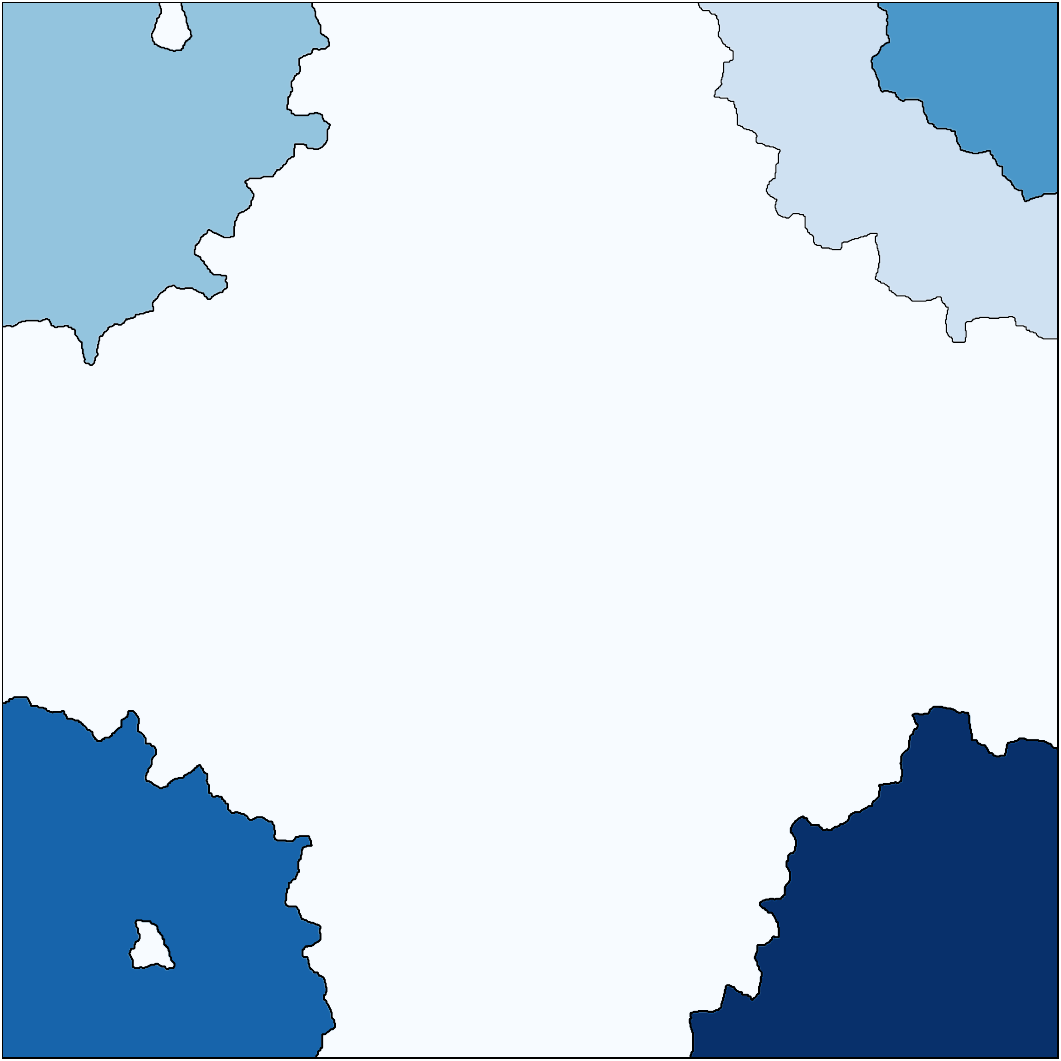}  
		\end{minipage}
	}%
	\subfigure[2000d]{\label{fig:ADDM:case1:2000d-cs02} 
		\begin{minipage}[t]{0.23\linewidth}
			\centering
			\includegraphics[width=3.3cm]{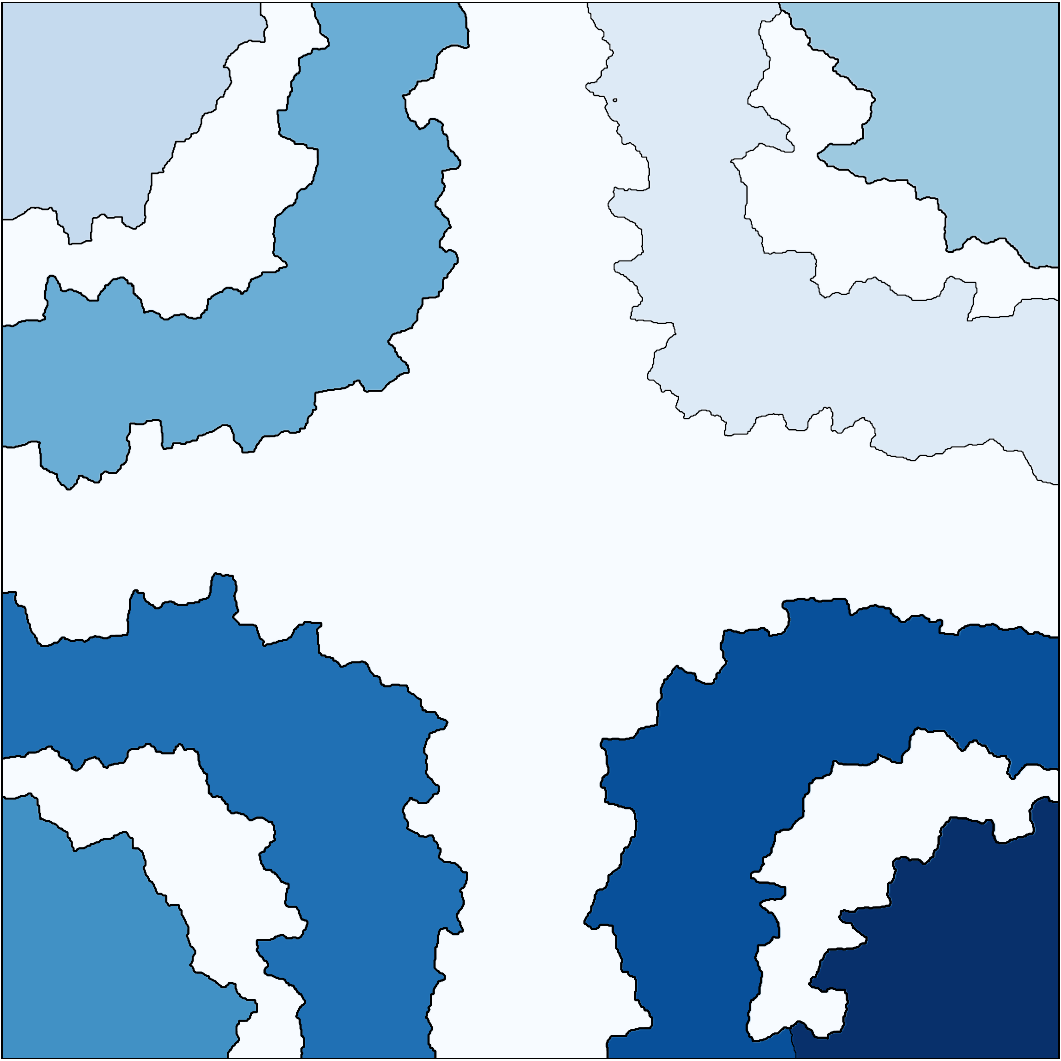}  
		\end{minipage}
	}%
	\subfigure[2820d]{\label{fig:ADDM:case1:2820d-cs02} 
		\begin{minipage}[t]{0.23\linewidth}
			\centering
			\includegraphics[width=3.3cm]{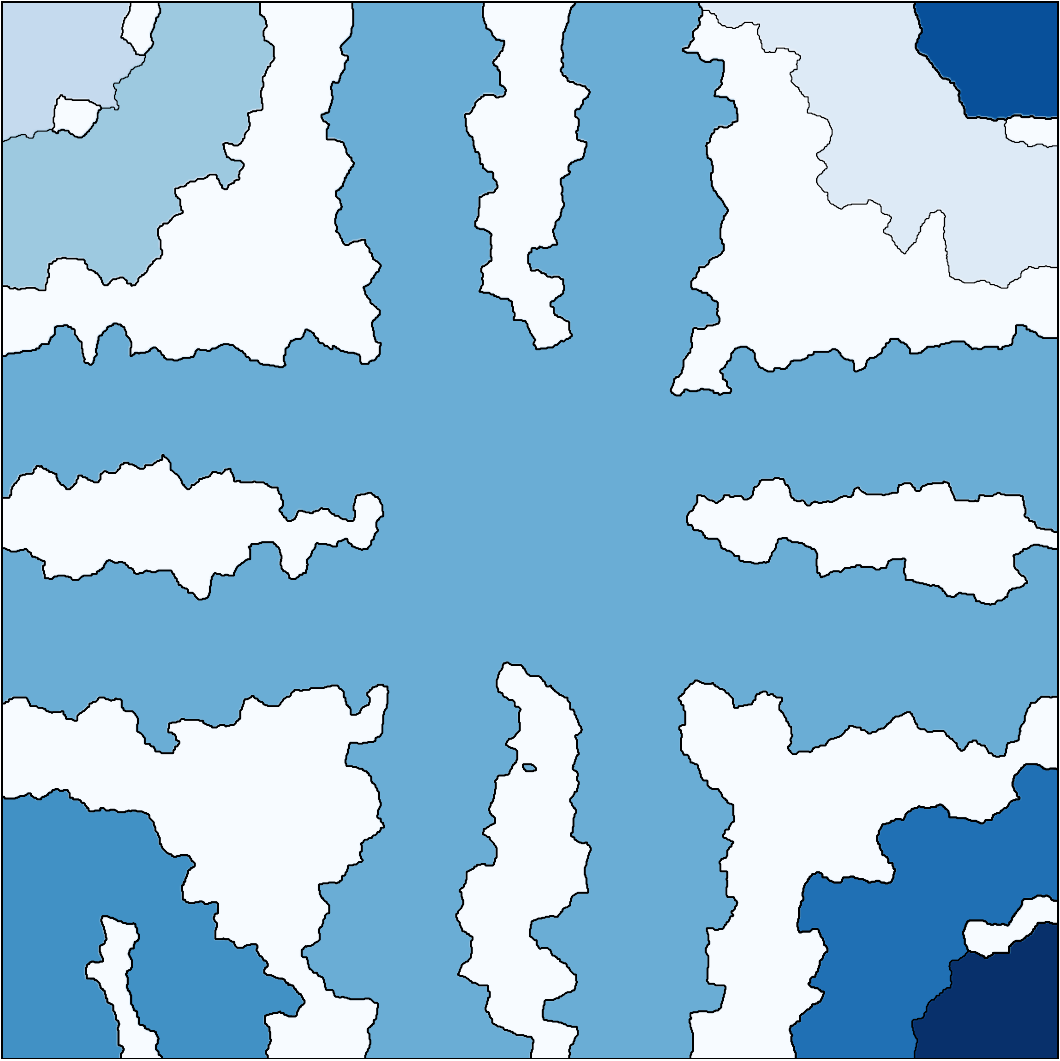}  
		\end{minipage}
	}%
	\subfigure[3000d]{\label{fig:ADDM:case1:3000d-cs02} 
		\begin{minipage}[t]{0.23\linewidth}
			\centering
			\includegraphics[width=3.3cm]{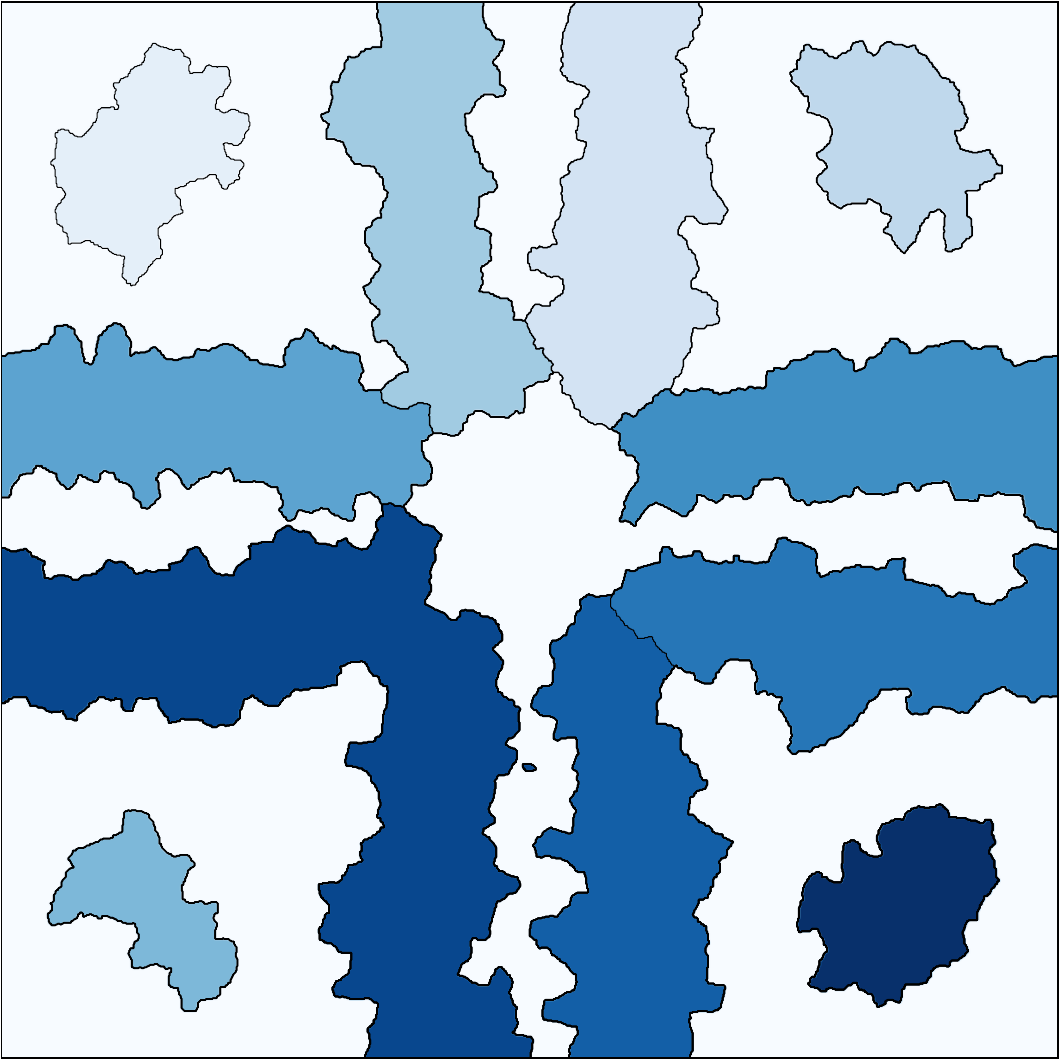}  
		\end{minipage}
	}%

	\subfigure[750d]{\label{fig:ADDM:case1:750d-cs03}  
		\begin{minipage}[t]{0.23\linewidth}
			\centering
			\includegraphics[width=3.3cm]{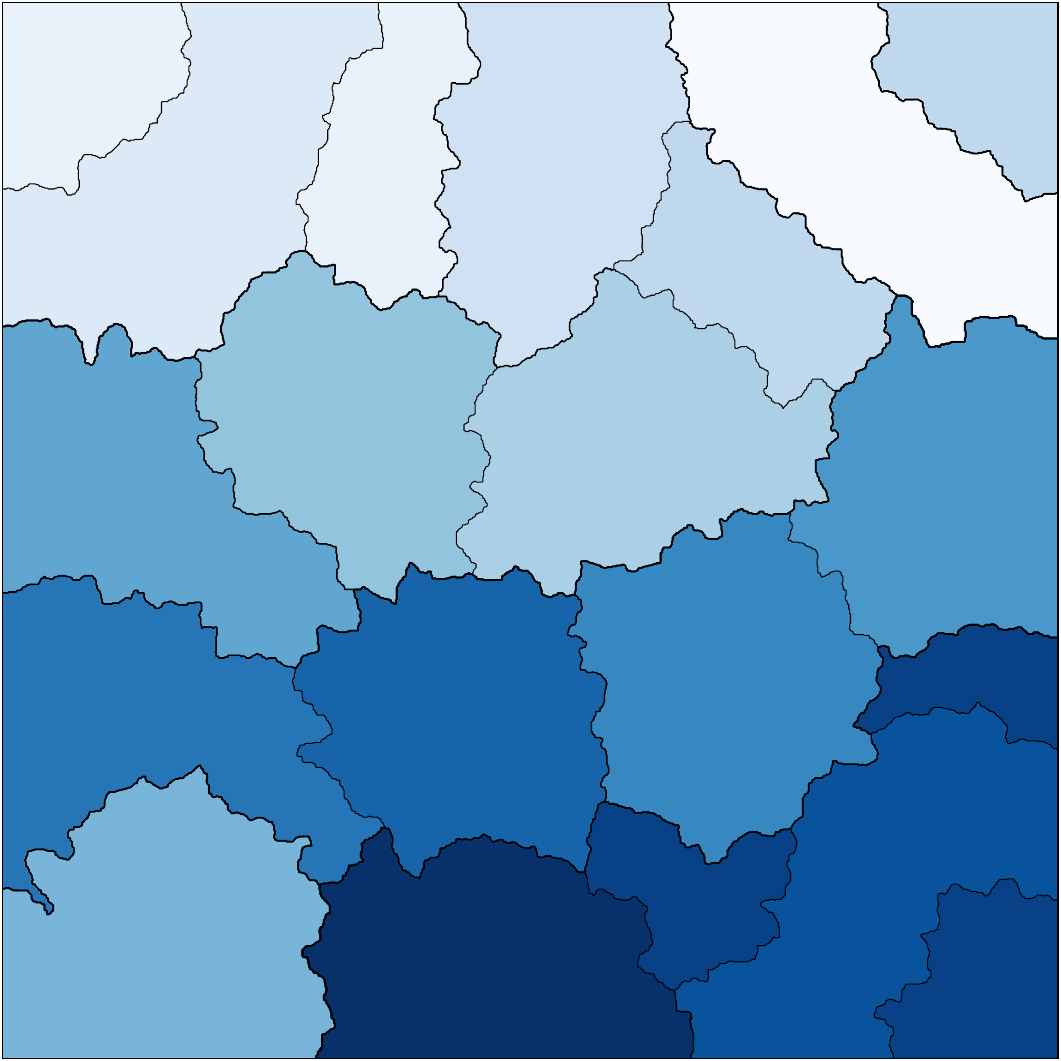}  
		\end{minipage}
	}%
	\subfigure[2000d]{\label{fig:ADDM:case1:2000d-cs03} 
		\begin{minipage}[t]{0.23\linewidth}
			\centering
			\includegraphics[width=3.3cm]{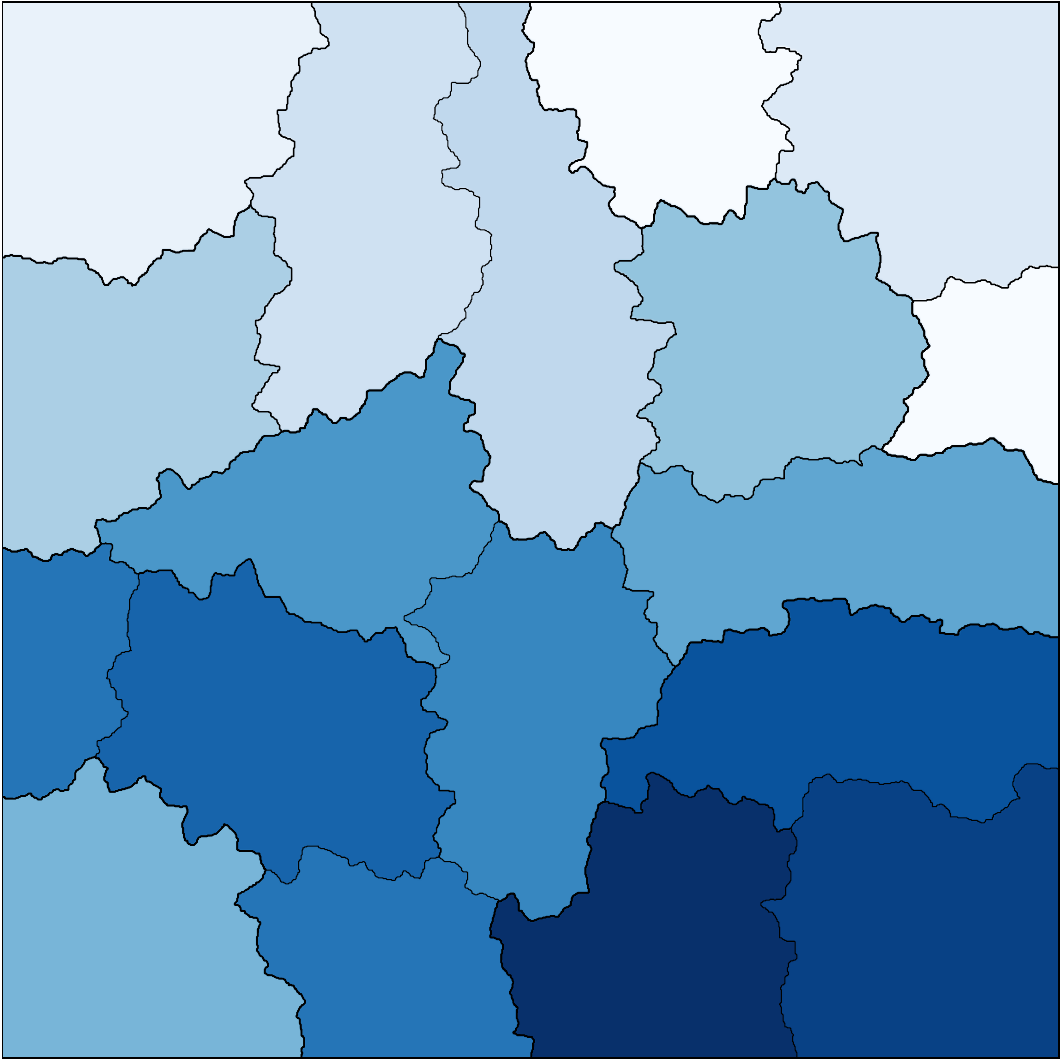}  
		\end{minipage}
	}%
	\subfigure[2820d]{\label{fig:ADDM:case1:2820d-cs03} 
		\begin{minipage}[t]{0.23\linewidth}
			\centering
			\includegraphics[width=3.3cm]{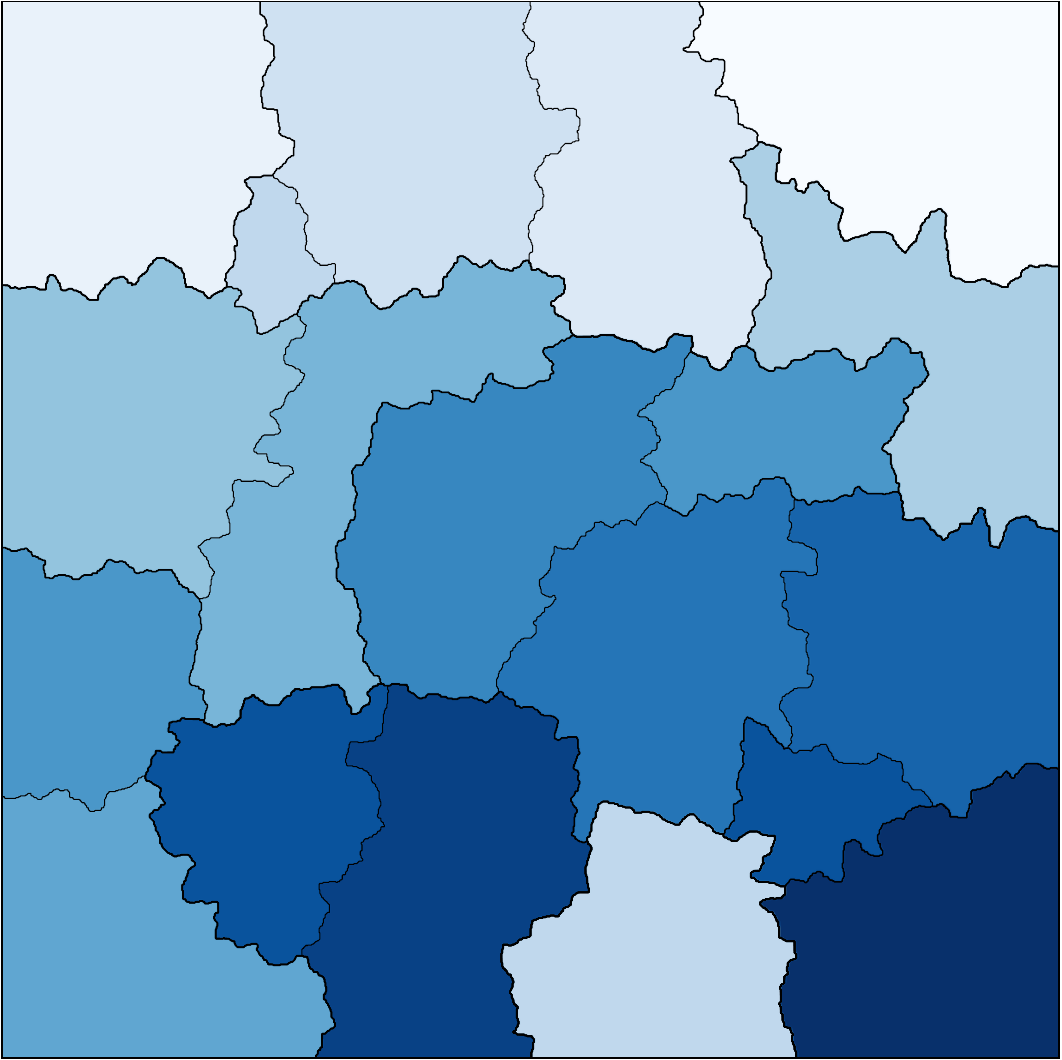}  
		\end{minipage}
	}%
	\subfigure[3000d]{\label{fig:ADDM:case1:3000d-cs03} 
		\begin{minipage}[t]{0.23\linewidth}
			\centering
			\includegraphics[width=3.3cm]{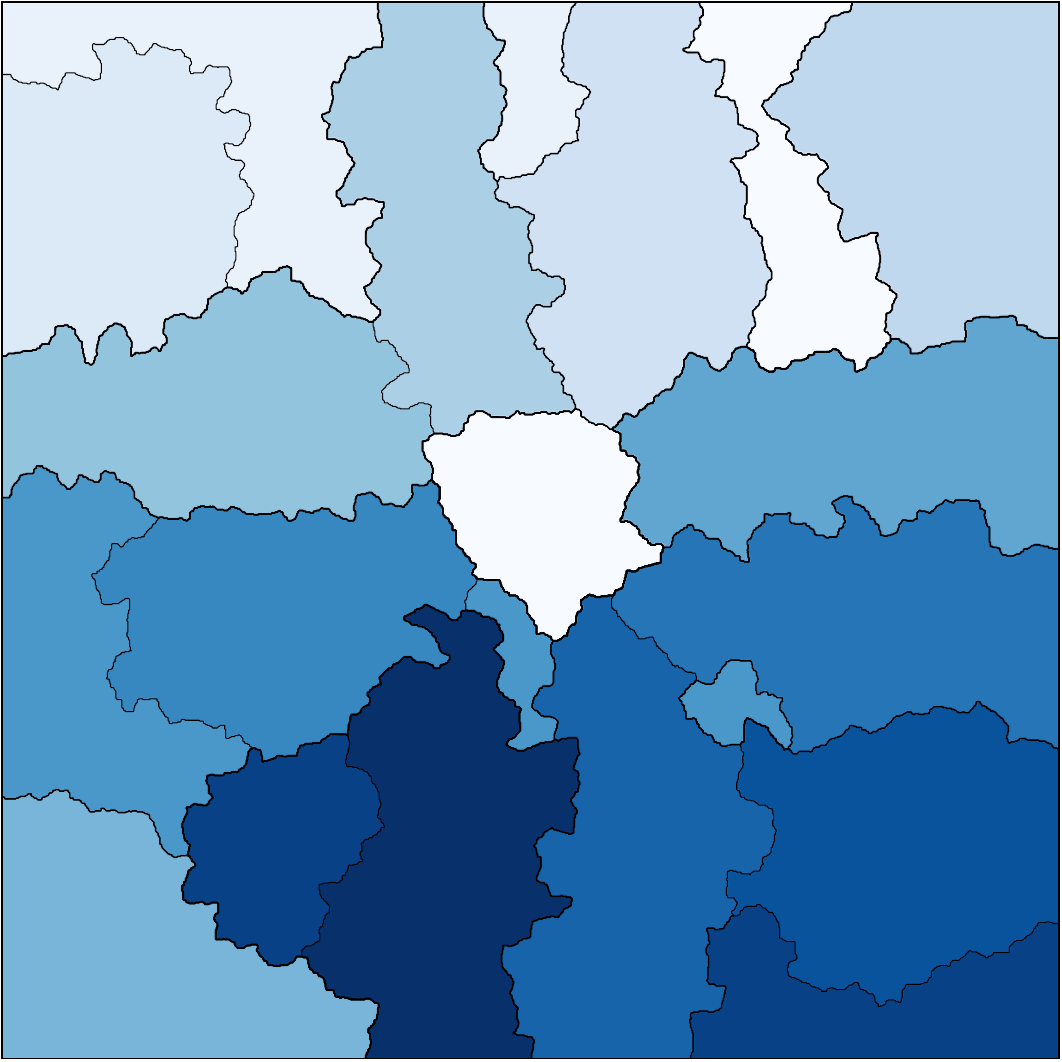}  
		\end{minipage}
	}%
    \caption{Comparison of results for different subdomain coupling strategies: (a)-(d) show the gas phase saturation in the top layer on day 750, day 2000, day 2820, and day 3000, respectively; \new{The closer the color is to red, the closer the saturation is to 1; the closer the color is to blue, the closer the saturation is to 0.} (e)-(h), (i)-(l), and (m)-(p) respectively show the subdomain coupling patterns at corresponding time points when using Strategy A, Strategy B, and Strategy C. \new{Subdomains colored in white will have their subproblems solved independently.} For the remaining categories, subdomains marked with the same color will be solved in a coupled manner.}  
    \label{fig:case1:3000Day:cstype}  
\end{figure}

This subsection examines the impact of three different subdomain adaptive coupling strategies and two boundary conditions on the performance of ADDM. Figure~\ref{fig:case1:3000Day:cstype} shows the gas saturation distribution in the top layer on days 750, 2000, 2820, and 3000, along with the corresponding subdomain coupling pattern for each of the three strategies, where the value of \( c_{_S} \) is set to \( 5 \times 10^{-3} \). For Strategy C, the maximum allowable number of coupled subdomains \( N_{\text{cs}} \) is set to 50. Note that during the simulation, gas primarily accumulates in the top layer due to its lower density; therefore, only the gas saturation distribution in the top layer is presented. Additionally, the gas evolution pattern in the lower layers mirrors that of the top layers, only delayed in time. Therefore, a limited range of subdomain coupling is often observed in the four corners of the subdomain coupling mode diagrams.

Figures~\ref{fig:ADDM:case1:750d-sags}-\ref{fig:ADDM:case1:3000d-sags} illustrate the three stages of fluid injection: pre-gas breakthrough (Figures~\ref{fig:ADDM:case1:750d-sags} and~\ref{fig:ADDM:case1:2000d-sags}), gas breakthrough (Figure~\ref{fig:ADDM:case1:2820d-sags}), and post-gas breakthrough (Figure~\ref{fig:ADDM:case1:3000d-sags}). 
The displacement front is relatively simple during the first stage but becomes increasingly complex in the second and third stages. 
Strategy A captures the displacement front using relatively few subdomains. 
Strategy B considers both the front and its movement, resulting in a larger coupled area and sometimes even in very large coupled regions. Strategy C captures the displacement features while allowing more subdomains to participate in the coupled solution, offering better control over the number and area of coupled subdomains, although it sacrifices some capability in capturing the displacement front.

Next, we compare the impact of different coupling parameters $c_{S}$ and boundary conditions on performance under these three different coupling strategies. For the two types of boundary conditions:
\begin{itemize} 
    \item Constant pressure boundary condition: For subdomain boundaries, the boundary condition uses the pressure values from neighboring cells in the previous time step, denoted as ADDM\_P.

    \item Constant flux boundary condition: For subdomain boundaries, the boundary condition uses the molar flux of components at the interface of neighboring cells in the previous time step, denoted as ADDM\_V.
\end{itemize}  

Table~\ref{tab:case1:3000day-3} presents details on coupling strategies (\texttt{Strategy}), the parameter $c_{_S}$, solution methods (\texttt{Method}), the number of time steps (\texttt{Timestep}), cumulative global Newton--Raphson iterations (\texttt{NRiter}), cumulative global linear iterations (\texttt{LSiter}), the local Newton--Raphson iterations required in the initial solution process
(\texttt{NRiter(DDM)}), and the total simulation runtime (\texttt{Runtime}).

\begin{table}[htpb]  
    \caption{Comparison of performance results for different boundary conditions and coupling strategies. The numbers in parentheses represent the iterations wasted due to solver failures. Bold indicates the best performance results.}  
    \label{tab:case1:3000day-3}  
    \centering
    \footnotesize
	\setlength{\tabcolsep} {1.0pt}  
    \renewcommand{\arraystretch}{1.2}
	\begin{tabular}{cccccccc}  
        \hline
        \texttt{Strategy} & $c_{_S}$ & \texttt{Method} & \texttt{Timestep} & \texttt{NRiter} & \texttt{LSiter} & \texttt{NRiter(DDM)} & \texttt{Runtime}(s)\\
        \hline
        \multirow{4}{*}{A}
        & \multirow{2}{*}{$1 \times 10^{-3}$}
			&ADDM\_P & 753 & 1899(+162) & 9242(+1231)  &2487(+192)  & 2130 \\
        &   &ADDM\_V & 735 & 1652(+260) & 7207(+1797) &2434(+113)  & \textbf{1902} \\
		\cdashline{2-8}
		& \multirow{2}{*}{$5 \times 10^{-3}$}
			&ADDM\_P & 776 & 1918(+202) & 9394(+1438)  &2509(+230)  & 2223 \\
        &   &ADDM\_V & 742 & 1887(+243) & 8595(+1833)  &2471(+116)  & \textbf{2029} \\
        \hline

        \multirow{4}{*}{B}
        & \multirow{2}{*}{$1 \times 10^{-3}$}
			&ADDM\_P & 754 & 1760(+142) & 8529(+1046)  &2475(+170)  & 1998 \\
        &   &ADDM\_V & 737 & 1489(+162) & 6512(+1189)  & 2452(+154)  & \textbf{1805} \\
		\cdashline{2-8}
		& \multirow{2}{*}{$5 \times 10^{-3}$}
			&ADDM\_P & 732 & 1719(+112) & 8431(+894)  & 2461(+135)  & 1910 \\
        &   &ADDM\_V & 721 & 1529(+160) & 6847(+1075)  & 2409(+141)  &  \textbf{1796} \\
        \hline

        \multirow{4}{*}{C}
        & \multirow{2}{*}{$1 \times 10^{-3}$}
			&ADDM\_P & 722 & 1827(+140) & 8973(+1002)  & 2411(+139)  & 2064 \\
        &   &ADDM\_V & 723 & 1852(+190) & 8342(+1224)  &2407(+150)  & \textbf{1991} \\
		\cdashline{2-8}
		& \multirow{2}{*}{$5 \times 10^{-3}$}
			&ADDM\_P & 737 & 1890(+171) & 9371(+1211)  &2451(+114)  & 2045 \\
        &   &ADDM\_V & 736 & 1857(+170) & 8362(+1263)  &2450(+155)  & \textbf{1968} \\
        \hline
    \end{tabular}
\end{table}

As shown in Table~\ref{tab:case1:3000day-3}, ADDM\_V outperforms ADDM\_P in all settings. Specifically, with the same number of \texttt{NRiter(DDM)}, ADDM\_V provides superior initial values, resulting in fewer global Newton--Raphson and linear iterations, and ultimately delivering better overall performance.
This outcome is intuitive, as the flow rate, which is related to the pressure gradient, varies more smoothly than pressure itself. 
Consequently, using the fixed flow rate from the previous time step as the boundary condition is a more effective choice.
In terms of coupling strategy, Strategy B outperforms the other two strategies. Compared to Strategy A, it involves a larger coupled area, which allows for a more complete capture of the nonlinearities near the displacement front and its movement. This is further supported by the observation that, within Strategy A, a smaller value of $c_{_S}$ leads to better performance.
By contrast, compared to Strategy C, it offers greater flexibility in capturing the displacement front, thereby reducing the risk of losing key coupling relationships.
It is noteworthy that when Strategy B is applied (with \( c_{_S} = 5 \times 10^{-3} \)) and a constant flow rate boundary condition is used, ADDM achieves a speedup of up to 577 seconds (24.3\%) compared to Standard.

\subsection{Case 2}\label{subsec:Case 2}
\begin{figure}[H]
	\setlength{\abovecaptionskip}{-0.05cm}
	\centering
	\subfigure[Grid layers 1 to 2]{\label{fig:case2:perm1}
		\begin{minipage}[t]{0.31\linewidth}
			\centering
			\includegraphics[trim=0.4cm 3.0cm 2.5cm 3cm,clip,width=4.4cm]{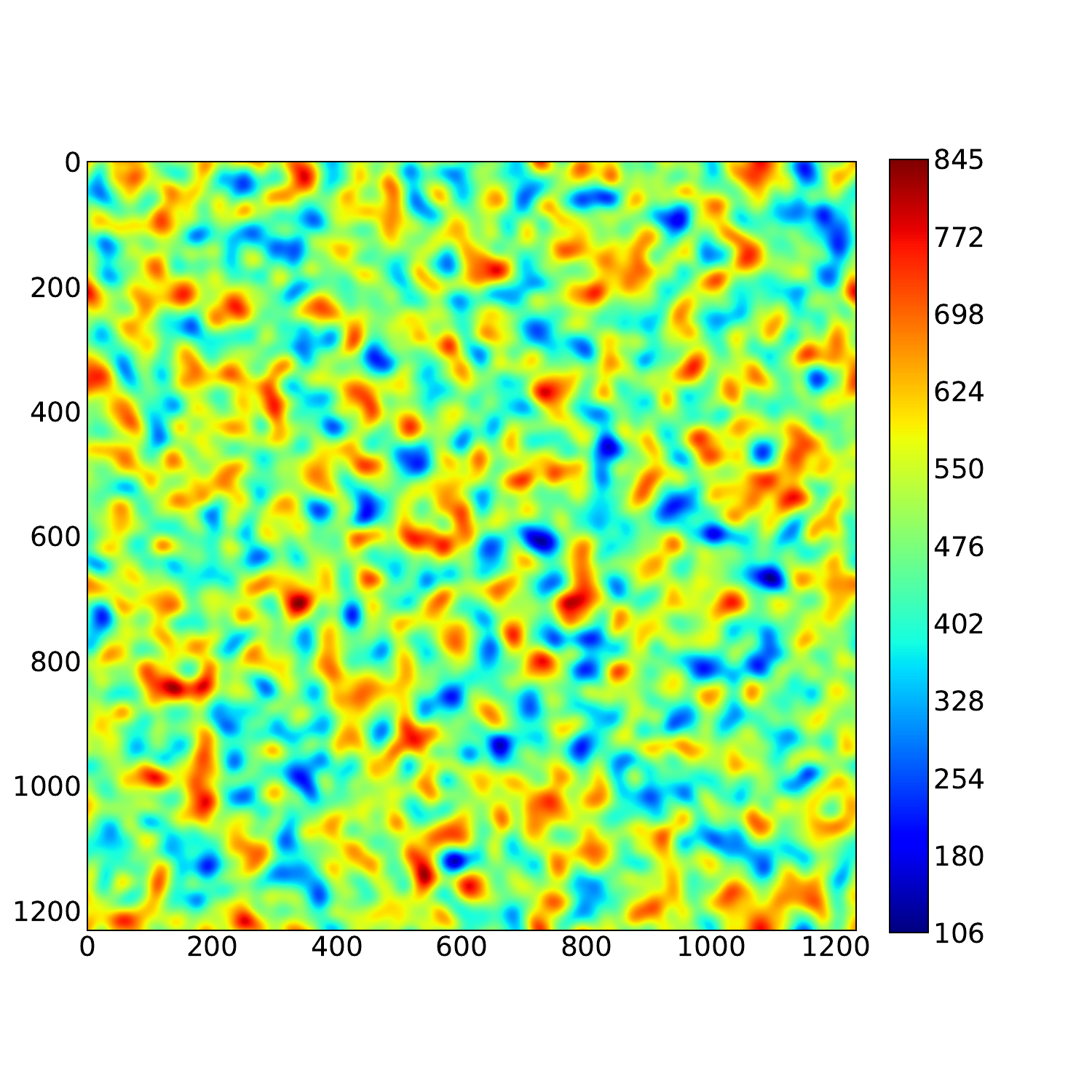}
		\end{minipage}
	}%
	\subfigure[Grid layers 3 to 5]{\label{fig:case2:perm2}
		\begin{minipage}[t]{0.31\linewidth}
			\centering
			\includegraphics[trim=0.4cm 3.0cm 2.5cm 3cm,clip,width=4.4cm]{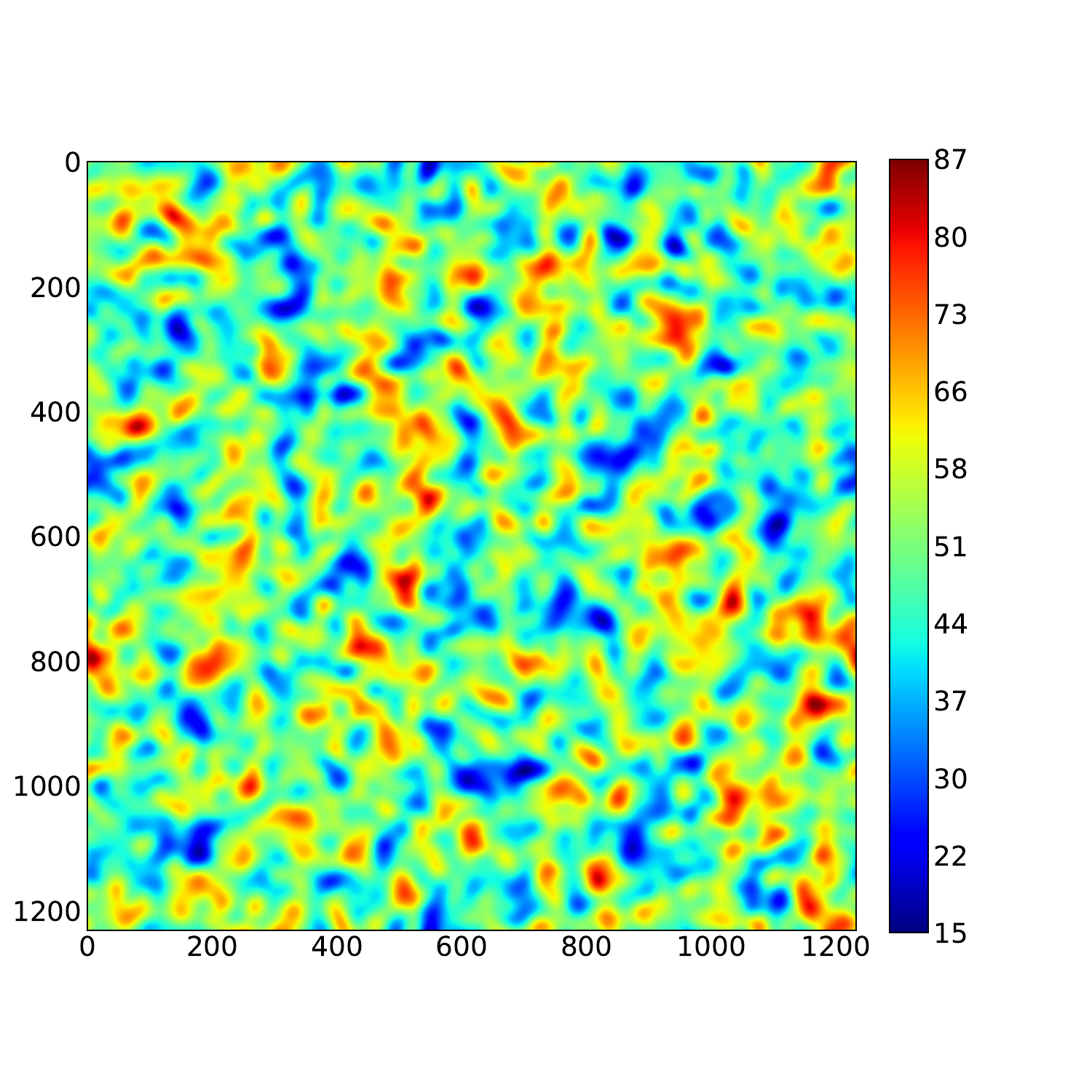}
		\end{minipage}
	}%
	\subfigure[Grid layers 6 to 10]{\label{fig:case2:perm3}
		\begin{minipage}[t]{0.31\linewidth}
			\centering
			\includegraphics[trim=0.4cm 3.0cm 2.5cm 3cm,clip,width=4.4cm]{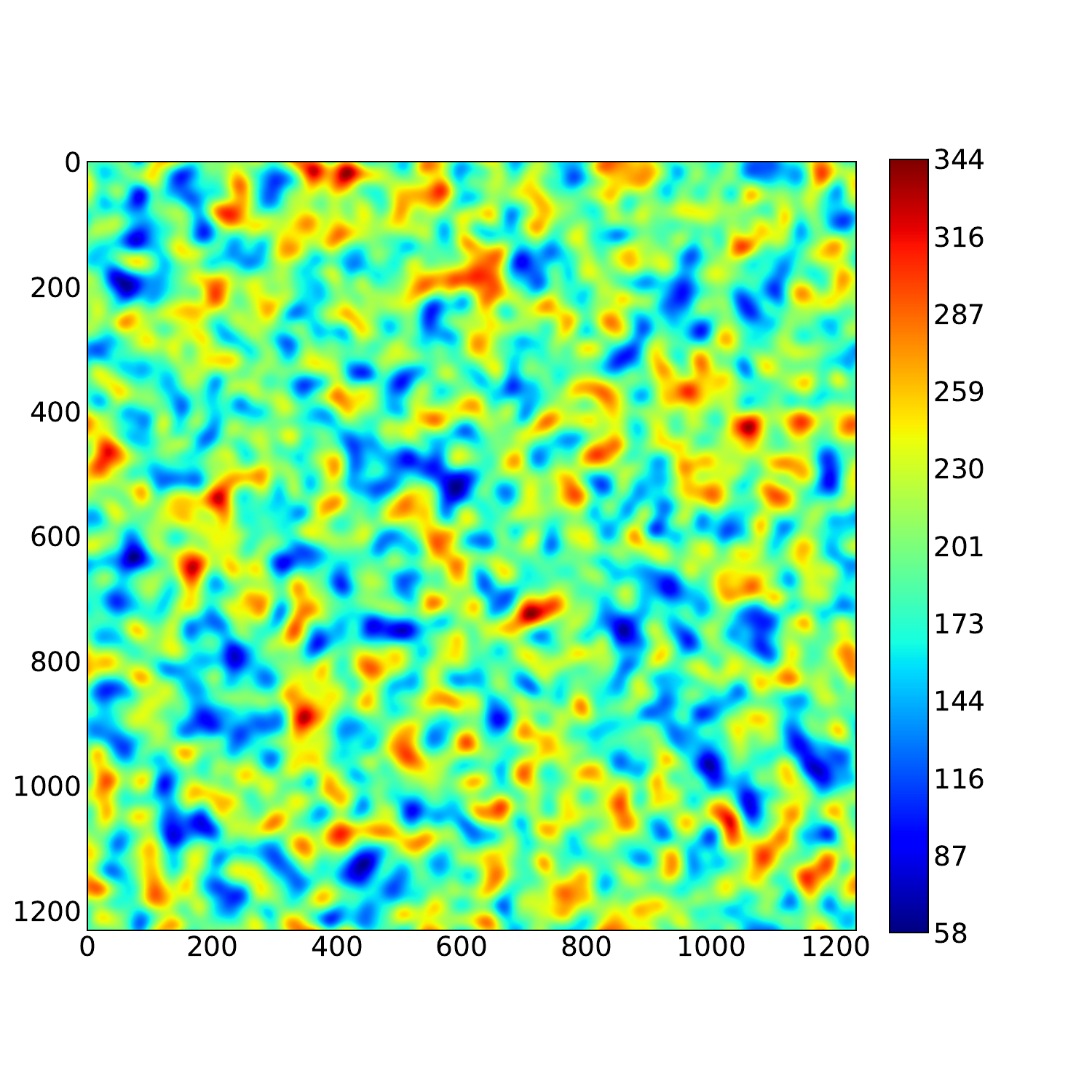}
		\end{minipage}
	}%
    \caption{The horizontal permeability distribution (in mD) for grid layers 1 to 2, 3 to 5, and 6 to 10. }
    \label{fig:case2:perm} 
\end{figure}

We introduced heterogeneity into Case 1 to provide a more comprehensive assessment of ADDM's performance.
The original SPE1 case features a three-layer geological structure with horizontal rock permeability of 500 mD, 50 mD, and 200 mD, respectively. 
In Case 1, following vertical grid refinement, these layers correspond to grid layers 1 to 2, 3 to 5, and 6 to 10, respectively. 
Building on Case 1, we introduced heterogeneity by adjusting the horizontal permeability of these three geological layers.
As shown in Figure~\ref{fig:case2:perm}, the horizontal permeability in each layer follows a Gaussian distribution, with mean values consistent with the original case and a standard deviation equal to 10 times the mean.

\begin{figure}[h]
	\centering
	\subfigure[900d]{\label{fig:ADDM:case2:900d-sags} 
		\begin{minipage}[t]{0.23\linewidth}
			\centering
			\includegraphics[width=3.3cm]{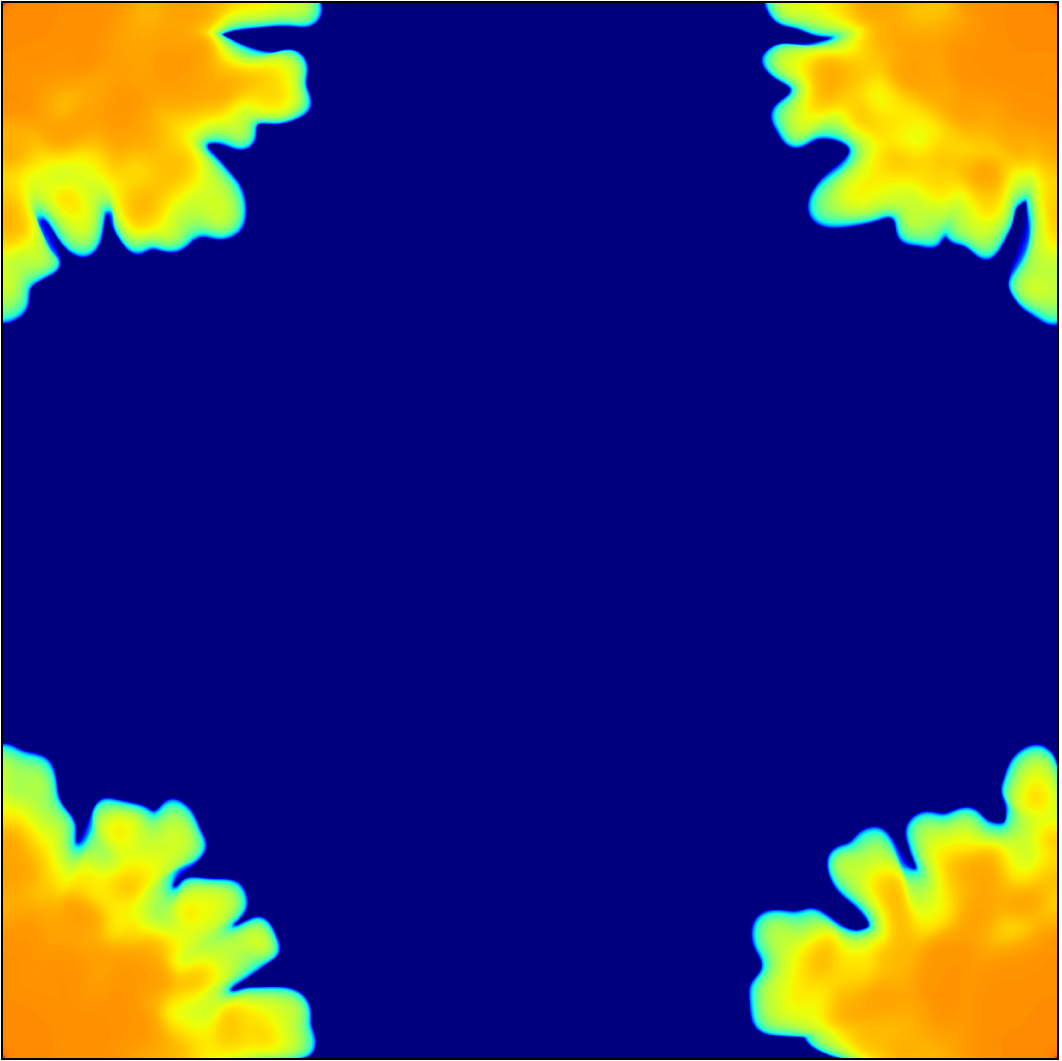}
		\end{minipage}
	}%
	\subfigure[1800d]{\label{fig:ADDM:case2:1800d-sags} 
		\begin{minipage}[t]{0.23\linewidth}
			\centering
			\includegraphics[width=3.3cm]{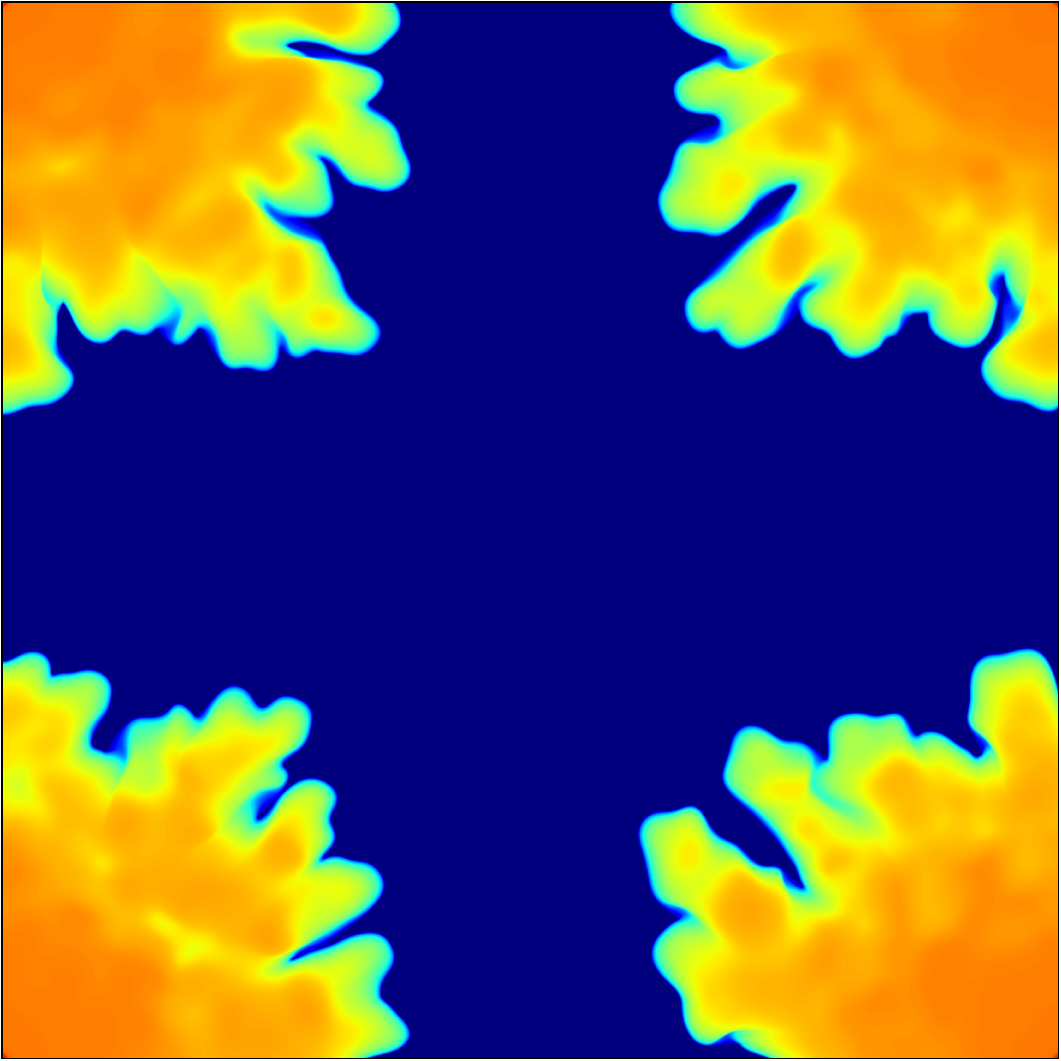}
		\end{minipage}
	}%
	\subfigure[2400d]{\label{fig:ADDM:case2:2400d-sags} 
		\begin{minipage}[t]{0.23\linewidth}
			\centering
			\includegraphics[width=3.3cm]{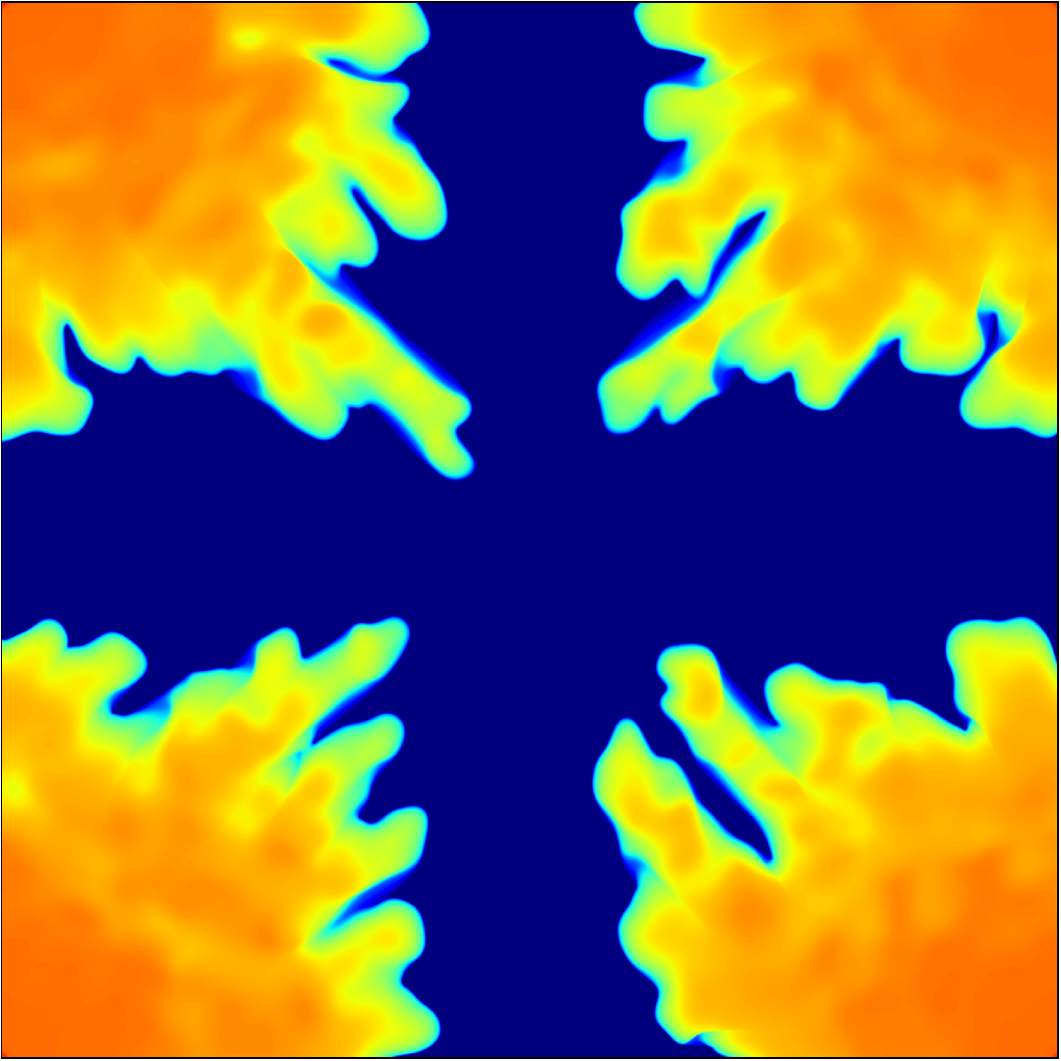}
		\end{minipage}
	}%
	\subfigure[3000d]{\label{fig:ADDM:case2:3000d-sags} 
		\begin{minipage}[t]{0.23\linewidth}
			\centering
			\includegraphics[width=3.3cm]{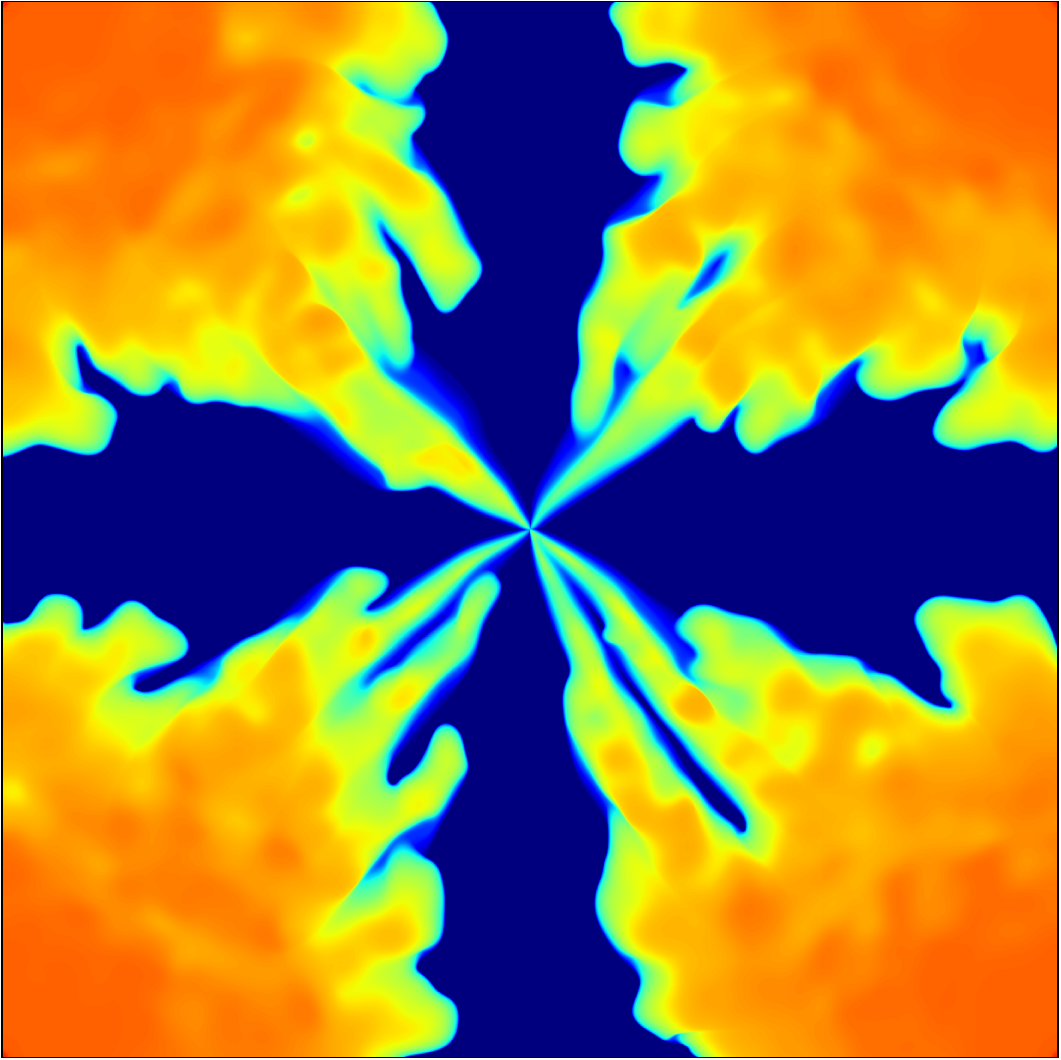}
		\end{minipage}
	}%

	\subfigure[900d]{\label{fig:ADDM:case2:900d-cs02} 
		\begin{minipage}[t]{0.23\linewidth}
			\centering
			\includegraphics[width=3.3cm]{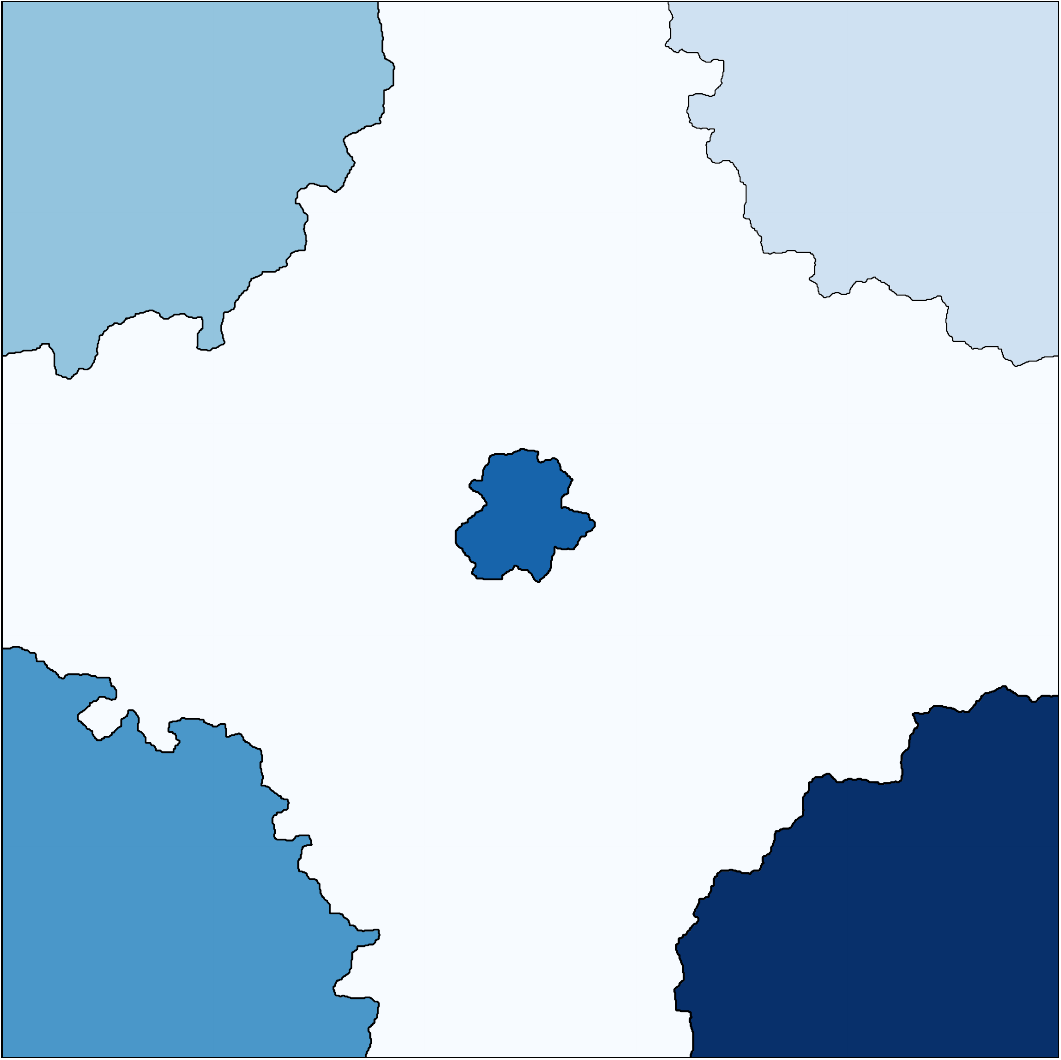}
		\end{minipage}
	}%
	\subfigure[1800d]{\label{fig:ADDM:case2:1800d-cs02} 
		\begin{minipage}[t]{0.23\linewidth}
			\centering
			\includegraphics[width=3.3cm]{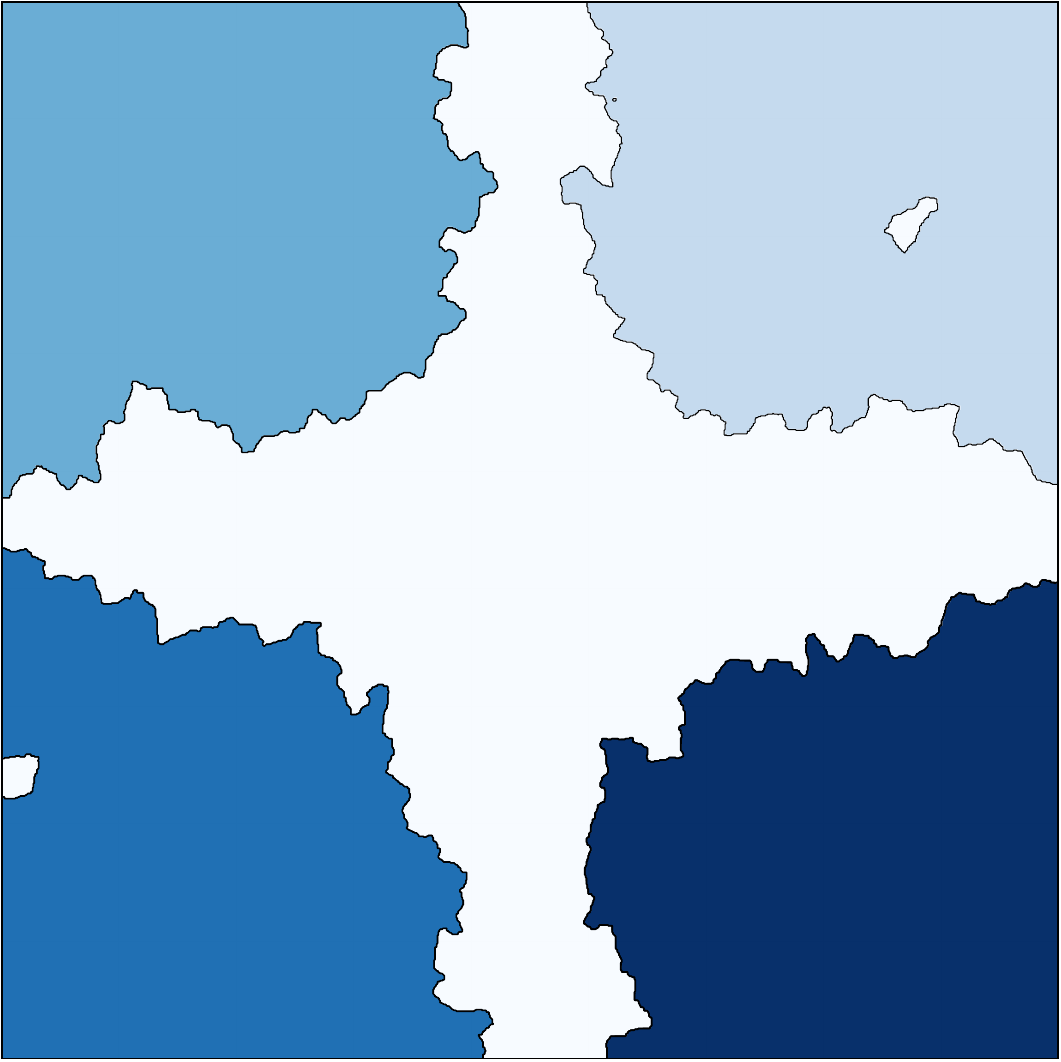}
		\end{minipage}
	}%
	\subfigure[2400d]{\label{fig:ADDM:case2:2400d-cs02} 
		\begin{minipage}[t]{0.23\linewidth}
			\centering
			\includegraphics[width=3.3cm]{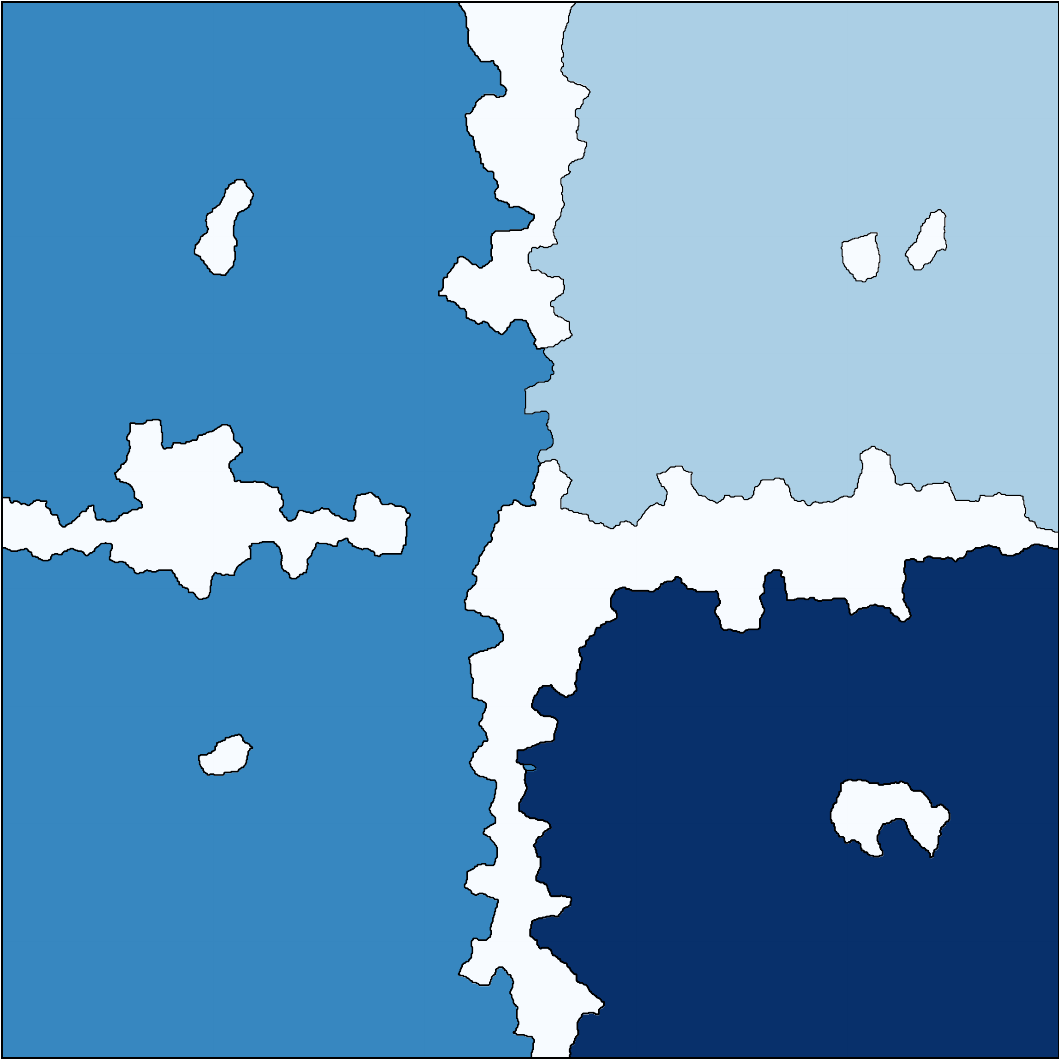}
		\end{minipage}
	}%
	\subfigure[3000d]{\label{fig:ADDM:case2:3000d-cs02} 
		\begin{minipage}[t]{0.23\linewidth}
			\centering
			\includegraphics[width=3.3cm]{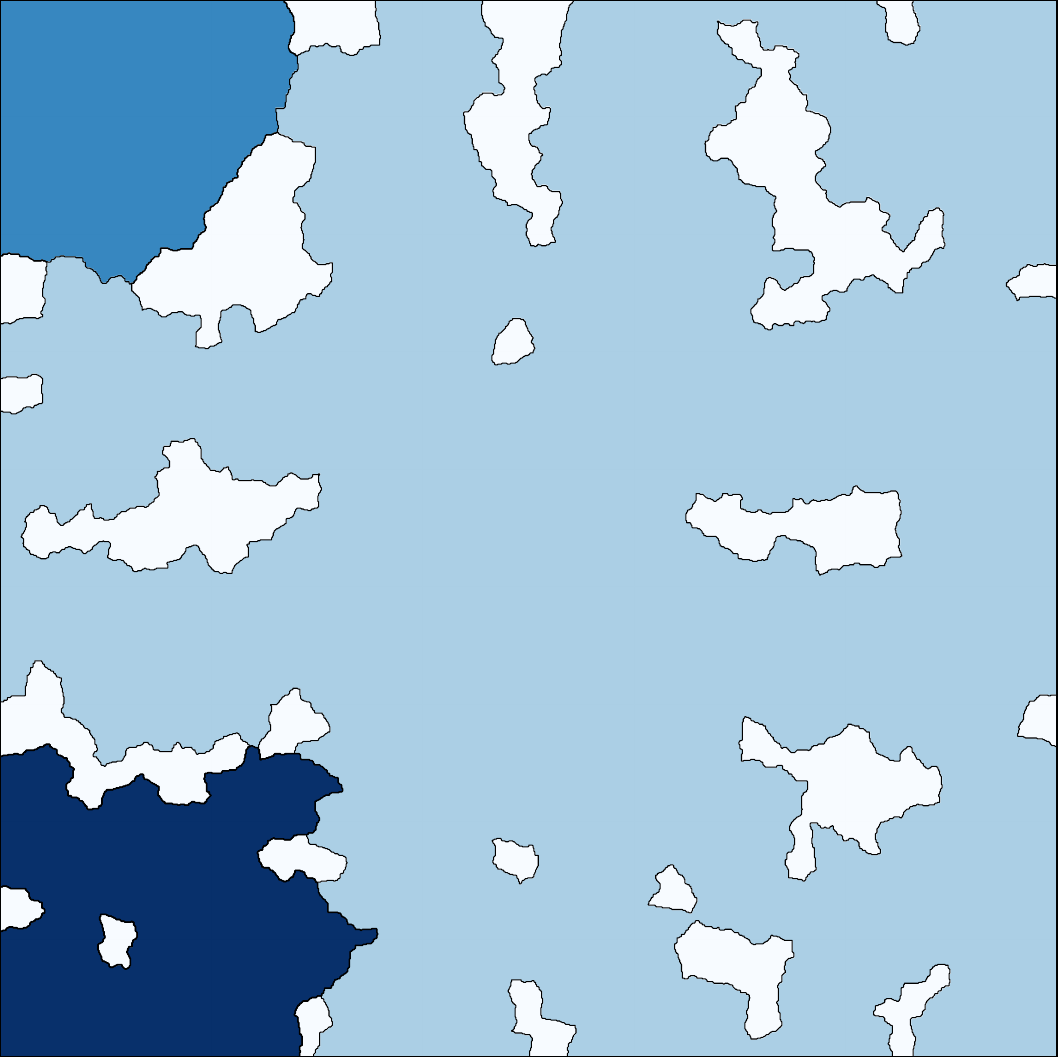}
		\end{minipage}
	}%
    \caption{Comparison of results for subdomain coupling Strategy B: (a)-(d) show the gas phase saturation in the top layer on day 900, day 1800, day 2400, and day 3000, respectively; \new{The closer the color is to red, the closer the saturation is to 1; the closer the color is to blue, the closer the saturation is to 0.} (e)-(h) show the subdomain coupling patterns at corresponding time points. \new{Subdomains colored in white will have their subproblems solved independently.} For the remaining categories, subdomains marked with the same color will be solved in a coupled manner.} 
    \label{fig:case2:cstype}
\end{figure} 

Figure~\ref{fig:case2:cstype} presents the gas phase saturation distribution at the top layer on days 900, 1800, 2400, and 3000, along with the subdomain coupling pattern associated with Strategy B. In this case, the substantial heterogeneity of the medium leads to complex fluid flow, making it challenging to accurately capture the shape and movement of the displacement front. 
The fluid distribution within the displacement front is uneven and evolves over time. 
In this condition, \( c_{_S} \) is set to $ = 10^{-3}$. 
Specifically, for ADDM, the constant flux boundary condition will be used. Compared to Case 1, Strategy B leads to a significant expansion of the coupled subdomain areas, particularly in the later stages of the simulation, which consequently increases the computational cost of solving for the initial values.
Nevertheless, it still demonstrates considerable performance improvements.

Table~\ref{tab:case2:3000day} presents a detailed comparison of the performance results.
The introduction of heterogeneity significantly increased the computational complexity of the problem. 
Compared to Case 1, Standard shows a 38.3\% increase in time steps, a 54.2\% increase in global Newton--Raphson iterations, and a 61.9\% increase in global linear iterations, ultimately resulting in a 60.8\% increase in total runtime.
Similarly, compared to Standard, ASM\_CDDM greatly reduces the required global Newton--Raphson iterations (31.4\%) and global linear iterations (36.4\%). However, this reduction is not sufficient to offset the cost of solving for the initial values, resulting in an overall increase in total runtime.
ASM\_ADDM, which utilizes the subdomain adaptive coupling strategy, further reduces global Newton--Raphson iterations (39.9\%) and global linear iterations (56.2\%), while requiring only 68.9\% of the local Newton--Raphson iterations of ASM\_CDDM. This leads to a performance improvement, resulting in a runtime reduction of 163 seconds.
In ADDM, the global Newton--Raphson iterations and global linear iterations are further reduced by 53.9\% and 59.8\%, respectively, resulting in a runtime reduction of 730 seconds (19.1\%). 
This highlights that, in the proposed algorithm, the constant flux boundary condition is the preferred choice, showcasing superior performance even in the face of complex conditions.

\begin{table} [htpb]  
    \caption{Performance comparison results of four methods for Case 2. The numbers in parentheses represent the iterations wasted due to solver failures. Bold indicates the best performance results.}  
    \label{tab:case2:3000day}  
    \centering
    \small  
    \setlength{\tabcolsep} {3pt}  
    \renewcommand{\arraystretch} {1.2}  
	\begin{tabular} {cccccccc}  
        \hline
        \texttt{Method} & \texttt{Timestep} & \texttt{NRiter} & \texttt{LSiter} & \texttt{NRiter(DDM)}& \texttt{Runtime}(s) \\
        \hline
        Standard       & 1003 & 4983(+1014) & 28574(+7939) & 0 & 3816 \\
        ASM\_CDDM & 1136 & 3363(+753) & 18091(+5144) & 4801(+772) & 4291 \\
        ASM\_ADDM & 1005 & 2875(+730) & 15334(+5178) & 3532(+310) & 3653 \\
        ADDM      & 1006 & 2225(+539) & 10937(+3734) & 3411(+442) & \textbf{3086} \\
        \hline
    \end{tabular} 
\end{table}

\new{
\subsection{Case 3}\label{subsec:Case 3}
This case is a refined version of the SPE5 benchmark~\citep{Killough_SPE5} and represents a compositional reservoir problem involving six components ($\rm C_1, C_3$, $\rm C_6$, $\rm C_{10}, C_{15}$, and $\rm C_{20}$), with both injection and production wells. The reservoir domain measures \(3500 \text{ ft} \times 3500 \text{ ft} \times 100 \text{ ft}\), and the original orthogonal grid consists of \(7 \times 7 \times 3\) cells. To evaluate the performance of the proposed methods for compositional reservoir simulations, the original grid is refined to \(1400 \times 1400 \times 30\), and the system is simulated over a period of 70 days using 2048 processes.

Table~\ref{tab:case3:70day} presents a performance comparison of the Standard, ASM\_CDDM, ASM\_ADDM, and ADDM methods, including the number of time steps (\texttt{Timestep}), cumulative global Newton--Raphson iterations (\texttt{NRiter}), cumulative global linear iterations (\texttt{LSiter}), the local Newton--Raphson iterations required in the initial solution process (\texttt{NRiter(DDM)}), and total simulation runtime (\texttt{Runtime}). 

\begin{table} [htpb]  
\new{
    \caption{\new{Performance comparison results of four methods for Case 3. The numbers in parentheses represent the iterations wasted due to solver failures. Bold indicates the best performance results.}}  
    \label{tab:case3:70day}  
    \centering
    \small  
    \setlength{\tabcolsep}{5.5pt}  
    \renewcommand{\arraystretch} {1.2}  
	\begin{tabular} {cccccccc}  
        \hline
        \texttt{Method} & \texttt{Timestep} & \texttt{NRiter} & \texttt{LSiter} & \texttt{NRiter(DDM)}& \texttt{Runtime}(s) \\
        \hline
        Standard       &56	&378(+86)	&2198(437)	&0			&1167  \\
        ASM\_CDDM &57	&193(+24)	&1094(146)	&531(+243)	&1297  \\
        ASM\_ADDM &48	&112(+24)	&714(+143)	&294(+91)	&814   \\
        ADDM      &47	&107(+5)	&525(+12)	&293(+128)	&\textbf{783}  \\
        \hline
    \end{tabular} 
}
\end{table}

As shown in Table~\ref{tab:case3:70day}, compared with the Standard method, ASM\_CDDM reduces the number of global Newton--Raphson iterations by 53.2\% and global linear iterations by 52.9\%, at the cost of an 11.1\% increase in runtime; ASM\_ADDM achieves reductions of 70.7\% and 67.5\%, respectively, together with a 30.2\% decrease in runtime; and ADDM further improves these reductions to 75.9\% and 79.6\%, respectively, yielding a 32.9\% reduction in runtime.
The results demonstrate that the proposed methods, which rely on a fixed saturation-change threshold, exhibit clear advantages in compositional flow scenarios. However, relying solely on saturation variation may be insufficient to accurately capture the flow front in such complex compositional settings. Future work will focus on developing more comprehensive criteria that incorporate multiple physical variables to further enhance the robustness and accuracy of the front-identification strategy.
}

\subsection{Parallel scalability}\label{subsubsec:Parallel scalability}
The subsection focuses on the parallel scalability analysis of the proposed method, involving numerical simulations with up to 500 million grid elements and over 2 billion degrees of freedom. Strong scalability is assessed to evaluate parallel performance while maintaining a fixed total problem size, focusing on how effectively the method accelerates computations as the number of processes increases. Weak scalability is evaluated to assess parallel performance while keeping the problem size per process constant, examining whether the method can efficiently handle larger overall problems as more processes are added. This analysis provides insights into the method's performance efficiency and resource utilization across different scaling scenarios.

\subsubsection{Strong scalability test} 
\begin{figure} [htpb] 
    \centering
    \includegraphics[width=0.8\linewidth]{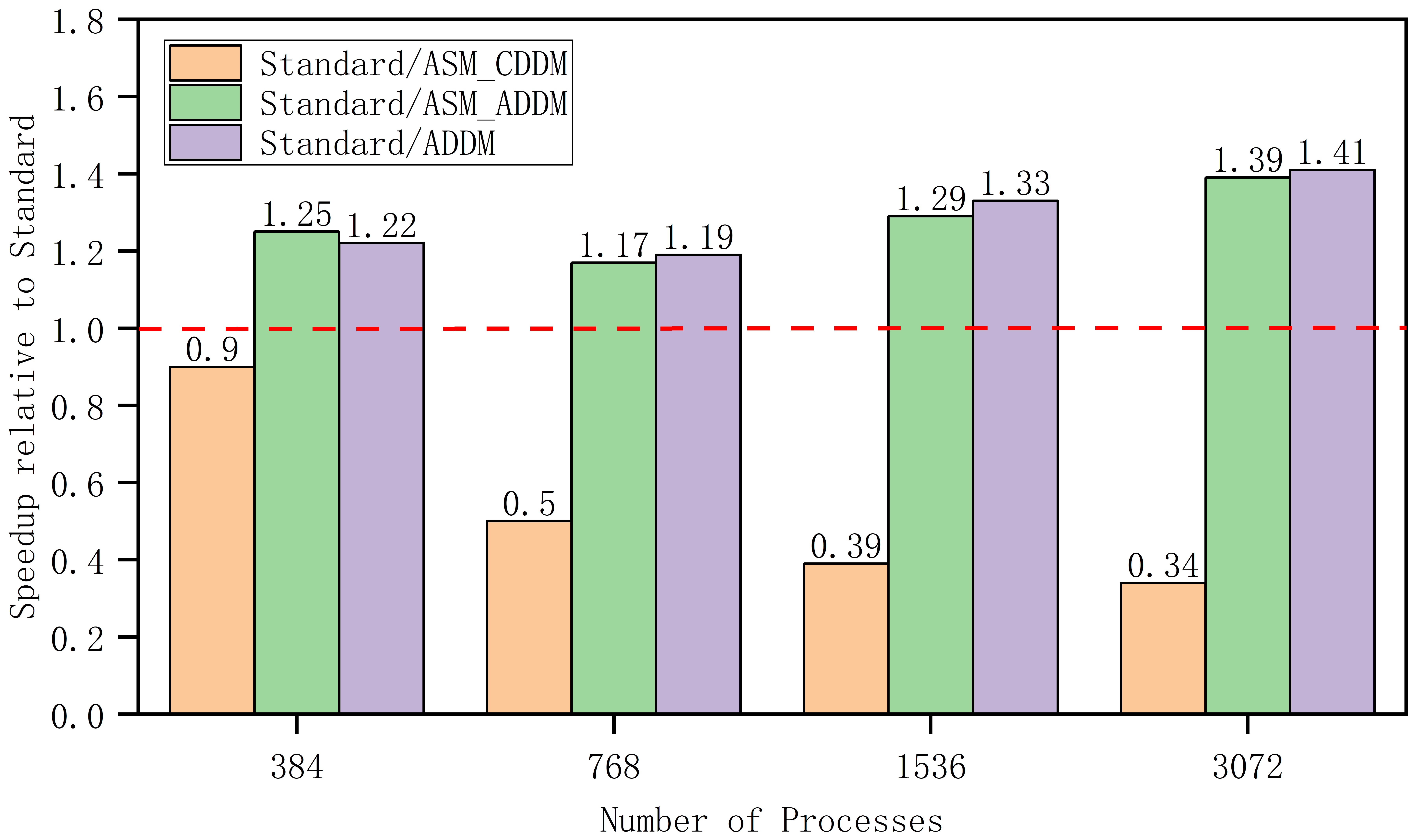} 
    \caption{Speedup of ASM\_CDDM, ASM\_ADDM, and ADDM relative to Standard under different process numbers in the parallel strong scalability test.} 
    \label{fig:Parallel-strong} 
\end{figure}

In Case 1, the grid is extended and refined, resulting in a final grid size of \( 1582 \times 1582 \times 50 \), comprising 125,136,200 grid cells and 500,544,800 degrees of freedom. Each grid cell measures \( 100 \text{ ft} \times 100 \text{ ft} \times 2 \text{ ft} \). The simulations are conducted with 384, 768, 1536, and 3072 processes. The total simulation duration is set to 1000 days.  
Table~\ref{tab:ch5:ADDM:strongP} summarizes the test results for four methods across different process numbers (\texttt{Np}), including the number of time steps (\texttt{Timestep}), cumulative global Newton--Raphson iterations (\texttt{NRiter}), cumulative global linear iterations (\texttt{LSiter}), the proportion of total simulation time spent on global linear solving time (\texttt{LSratio}), total simulation runtime (\texttt{Runtime}), and parallel efficiency (\texttt{PE}). Figure~\ref{fig:Parallel-strong} presents the speedup of ASM\_CDDM, ASM\_ADDM, and ADDM relative to Standard under different numbers of processes.

\begin{table} [h]
    \caption{Parallel strong scalability test results for four methods. Bold indicates the best performance results.}
    \label{tab:ch5:ADDM:strongP} 
    \centering
    \small
    \setlength{\tabcolsep} {3.5pt}
    \renewcommand{\arraystretch} {1.2} 
	\begin{tabular} {crcccccc} 
        \hline
        \texttt{Np} & \texttt{Method} & \texttt{Timestep} & \texttt{NRiter} & \texttt{LSiter}& \texttt{LSratio} & \texttt{Runtime}(s) & \texttt{PE} \\
        \hline
        \multirow{4} {*} {384}  & Standard & 182 & 858 & 3674 &82.3\%   &6957 &100\% \\
        & ASM\_CDDM & 197 & 559 & 2423 &47.8\%   &  7751 &100\%\\
        & ASM\_ADDM & 182 & 348 & 1618 &43.6\%   & \textbf{5558} &100\%\\
        & ADDM      & 182 & \textbf{346} & \textbf{1485} &{42.0\%}   & 5690 &100\%\\
        \hline
        \multirow{4} {*} {768}  & Standard & 182 & 857 & 3828 &82.9\%   & 3407 &{102\%} \\
        & ASM\_CDDM & 322 & 829 & 3490 &38.2\%   & 6836 & 57\% \\
        & ASM\_ADDM & 182 & 351 & 1655 &40.6\%   & 2920 & 95\% \\
        & ADDM      & 182 & \textbf{345} & \textbf{1530} &{39.7\%}   & \textbf{2860} & 99\% \\
        \hline
        \multirow{4} {*} {1536}  & Standard & 182 & 856 & 3977 &84.0\%   & 2060 & 84\%  \\
        & ASM\_CDDM & 415 & 991 & 4255 &36.6\%   &  5223 & 37\%  \\
        & ASM\_ADDM & 182 & 354 & 1662 &44.4\%   &  1603 & 87\% \\
        & ADDM      & 182 & \textbf{346} & \textbf{1514} &{42.2\%}   &  \textbf{1550} & {92\%} \\
        \hline
        \multirow{4} {*} {3072}  & Standard & 182 & 857 & 4052 &88.1\%   & 1500  & 58\% \\
        & ASM\_CDDM & 532 & 1191 & 5212 &37.7\%   & 4401 & 22\% \\
        & ASM\_ADDM & 182 & 355  & 1708 &47.5\%   & 1078 & 64\% \\
        & ADDM      & 182 & \textbf{351}  & \textbf{1536} &{46.0\%}   & \textbf{1064} & {67\%} \\
        \hline
    \end{tabular} 
\end{table}

From Table~\ref{tab:ch5:ADDM:strongP}, several observations can be made.
(1) ASM\_CDDM not only requires a large number of global Newton--Raphson iterations but also suffers from serious non-convergence issues, leading to frequent time step reductions and repeated computations. This problem worsens rapidly as the number of processes increases, making this method considerably slower than the others.
(2) ASM\_ADDM and ADDM, both of which employ the subdomain adaptive coupling strategy, significantly reduce the number of global Newton--Raphson and linear iterations at a low computational cost, thereby accelerating convergence. Moreover, even when the process number increases substantially, the linear iteration numbers remain nearly constant. This demonstrates that an appropriate subdomain adaptive coupling strategy can greatly enhance the robustness of ASM\_ADDM and ADDM as the number of processes grows. The stability of the global Newton--Raphson iteration number further indicates that the initial guesses provided by these methods maintain high quality and are not affected by increasing process numbers.
(3) ASM\_ADDM and ADDM also achieve higher parallel efficiencies. For example, with 3072 processes, the parallel efficiencies of Standard, ASM\_CDDM, ASM\_ADDM, and ADDM are 58\%, 22\%, 64\%, and 67\%, respectively. ADDM attains slightly higher parallel efficiency than ASM\_ADDM, primarily due to its further reduction of global linear iterations, which lowers the time spent on global linear solves. In addition, ADDM avoids boundary information exchange between subdomains when solving local problems.
Furthermore, from Figure~\ref{fig:Parallel-strong}, ASM\_CDDM consistently underperforms, exhibiting even longer runtimes than Standard, and its relative performance deteriorates as the number of processes increases, with the speedup dropping from 0.9 to 0.34. In contrast, ASM\_ADDM and ADDM show clear runtime advantages across all process numbers. Their speedups relative to Standard increase from approximately 1.2 at 384 processes to about 1.4 at 3072 processes, demonstrating improved efficiency with larger parallel configurations. 
Overall, these results demonstrate that the proposed method is both efficient and robust across a wide range of process numbers, and it indicates that appropriate handling of highly nonlinear local subproblems can provide high-quality initial solutions.

\subsubsection{Weak scalability test} 
In the weak scalability tests, we conduct numerical simulations on Case 1 over a 100-day period. Starting with a mesh size of $791 \times 791 \times 25$ (15.6 million cells) using 192 processes, we progressively scale up to a maximum mesh size of $3164 \times 3164 \times 50$ (500.5 million cells), involving 2002.2 million degrees of freedom. This largest case employs 6144 processes, with each process handling 0.3 million degrees of freedom. Given that previous tests reveal substantial declines in convergence and robustness for the ASM\_CDDM method at higher process counts, we exclude it from this comparison. Thus, our analysis focuses solely on the Standard, ASM\_ADDM, and ADDM methods.

\begin{figure}[htpb] 
    \centering
    \includegraphics[width=0.8\linewidth]{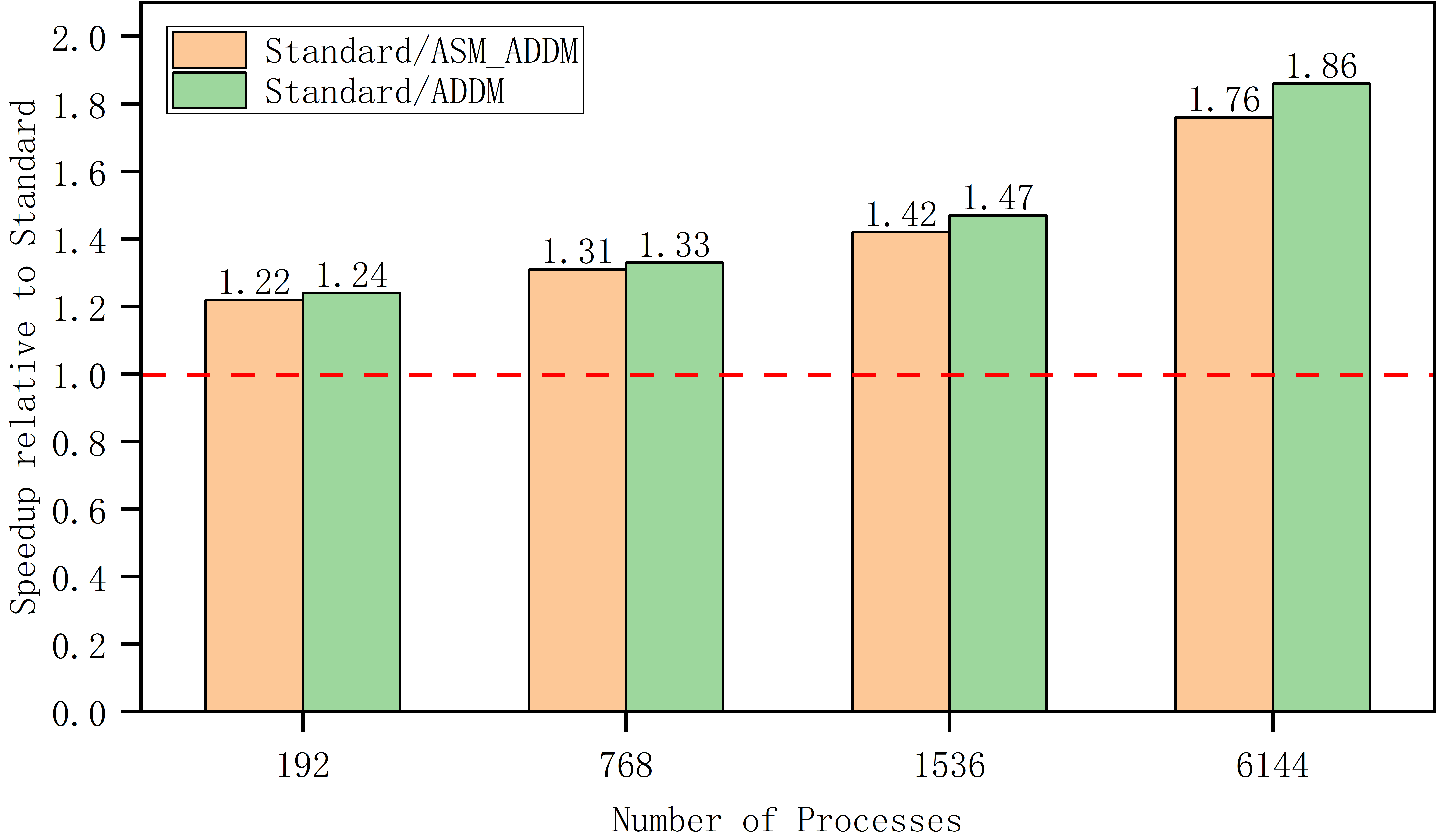} 
    \caption{Speedup of ASM\_ADDM and ADDM relative to Standard under different process numbers in the parallel weak scalability test. } 
    \label{fig:Parallel-weak} 
\end{figure}

Table~\ref{tab:ch5:ADDM:strongW} provides details on the number of mesh cells (\texttt{Nc}), number of processes (\texttt{Np}), solution methods (\texttt{Method}), number of time steps (\texttt{Timestep}), cumulative global Newton--Raphson iterations (\texttt{NRiter}), cumulative global linear iterations (\texttt{LSiter}), the average number of linear iterations per Newton step (\texttt{Avgiter}), the percentage of total simulation time spent to global linear solving (\texttt{LSratio}), and the total simulation runtime (\texttt{Runtime}). Figure~\ref{fig:Parallel-weak} presents the speedup of ASM\_ADDM and ADDM relative to Standard under different mesh cells and numbers of processes.

\begin{table} [htpb]
    \caption{The weak scaling results by using different numbers of processes for the Standard, ASM\_ADDM, and ADDM methods. Bold indicates the best performance results.} 
    \label{tab:ch5:ADDM:strongW} 
    \centering
    \footnotesize
    \setlength{\tabcolsep} {2.0pt} 
    \renewcommand{\arraystretch} {1.15} 
	\begin{tabular} {ccrcccccc} 
        \hline
        \texttt{Nc}(million) & \texttt{Np} & \texttt{Method} & \texttt{Timestep} & \texttt{NRiter} & \texttt{LSiter} & \texttt{Avgiter} & \texttt{LSratio} & \texttt{Runtime}(s)  \\
        \hline
        \multirow{3} {*} {15.6} & \multirow{3} {*} {192}  & Standard & 41 & 137 & 534 &3.9 &74.3\%   & 194  \\
        & & ASM\_ADDM & 41 & 60 & 253 &4.2 &42.3\%   & 159 \\
        & & ADDM      & 41 & 60 & \textbf{227} &\textbf{3.8} &\textbf{40.5\%}   & \textbf{157} \\
        \hline
        \multirow{3} {*} {62.6} & \multirow{3} {*} {768}  & Standard & 42 & 154 & 643 &4.2 &78.0\%   & 261  \\
        & & ASM\_ADDM & 42 & 56 & 269 &4.8 &40.4\%   & 200 \\
        & & ADDM      & 42 & \textbf{54} & \textbf{227} &4.2 &\textbf{36.5\%}   & \textbf{196} \\
        \hline
        \multirow{3} {*} { 125.1 } & \multirow{3} {*} {1536}  & Standard & 41 & 192 & 851 &4.4 &83.2\%   & 435  \\
        & & ASM\_ADDM & 41 & \textbf{67} & 312 &4.7 &43.3\%   &  307 \\
        & & ADDM      & 41 & 68 & \textbf{290} &\textbf{4.3} &\textbf{41.9\%}   &  \textbf{296} \\
        \hline
        \multirow{3} {*} { 500.5 } & \multirow{3} {*} {6144}  & Standard & 47 & 247 & 1215 &4.9 &89.2\%    & 1025  \\
        & & ASM\_ADDM & 44 & 95 & 469 &4.9 &53.4\%   & 581 \\
        & & ADDM      & 44 & \textbf{86} & \textbf{395} &\textbf{4.6} &\textbf{49.2\%}   & \textbf{550} \\
        \hline
    \end{tabular}
\end{table}

According to Table~\ref{tab:ch5:ADDM:strongW}, as the grid is progressively refined, both the problem complexity and the number of required global Newton--Raphson iterations increase. Additionally, the average number of linear iterations per Newton step also grows. For example, using Standard, as the number of mesh cells increases from 15.6 million to 500.5 million, the global Newton--Raphson iterations increase by 110 (an 80.3\% rise), while the average linear iterations per Newton step increase from 3.9 to 4.9 (a 26.2\% rise). Both ASM\_ADDM and ADDM, which utilize a subdomain adaptive coupling strategy, consistently outperform Standard by significantly reducing the number of global Newton--Raphson and linear iterations required. This reduction leads to decreased time spent on global linear solutions. 
In terms of runtime (see Figure~\ref{fig:Parallel-weak}), for 192, 768, 1536, and 6144 processes, ASM\_ADDM achieves speedups over Standard of 1.22, 1.31, 1.42, and 1.76, respectively, while ADDM achieves 1.24, 1.33, 1.47, and 1.86. The slightly higher speedups of ADDM are mainly due to its lower communication and synchronization overhead. These results indicate that the proposed methods gain increasingly significant advantages as the number of processes grows, maintaining high efficiency and robustness compared to Standard, even for complex and refined problems in large-scale parallel environments.

\new{
\subsection{Summary}
Overall, the results demonstrate that the proposed ASM\_ADDM/ADDM methods significantly reduce both nonlinear and linear iteration numbers, thereby yielding clear computational advantages. Compared with the Standard approach, ASM\_ADDM/ADDM reduces the numbers of nonlinear and linear iterations by more than 50\%. Moreover, the total simulation runtime is reduced by 20\%-80\%. It is also observed that these performance gains become more pronounced as the number of processes increases, indicating good parallel scalability of the proposed methods for large-scale simulations.}

\section{Conclusions}\label{sec:Conclusions}
In this work, we propose an adaptively coupled domain decomposition method (ADDM) framework for the fully implicit solution of multiphase and multicomponent flow in porous media. The solution methods developed within this framework effectively capture strong nonlinearities in global problems by defining subproblems in the coupled regions based on fluid flow characteristics, significantly accelerating the convergence of nonlinear solvers. Additionally, we introduce several adaptive coupling strategies and develop a nonlinear problem initialization method within this framework. Numerical experiments confirm the effectiveness of the proposed ADDM framework, using saturation to capture the dynamics of moving interfaces. The methods achieve good parallel performance in both strong and weak scalability, particularly for large-scale parallel applications. \new{In future work, more comprehensive criteria incorporating multiple physical variables will be explored for complex physical scenarios to further improve the robustness and accuracy of the front-identification strategy. Furthermore, the application of the ADDM method to preconditioning techniques will also be investigated.}

%
%




\bibliographystyle{elsarticle-harv} 
\bibliography{ref}

\begin{thebibliography}{51}
\expandafter\ifx\csname natexlab\endcsname\relax\def\natexlab#1{#1}\fi
\providecommand{\url}[1]{\texttt{#1}}
\providecommand{\href}[2]{#2}
\providecommand{\path}[1]{#1}
\providecommand{\DOIprefix}{doi:}
\providecommand{\ArXivprefix}{arXiv:}
\providecommand{\URLprefix}{URL: }
\providecommand{\Pubmedprefix}{pmid:}
\providecommand{\doi}[1]{\href{http://dx.doi.org/#1}{\path{#1}}}
\providecommand{\Pubmed}[1]{\href{pmid:#1}{\path{#1}}}
\providecommand{\bibinfo}[2]{#2}
\ifx\xfnm\relax \def\xfnm[#1]{\unskip,\space#1}\fi
\bibitem[{Aavatsmark(2002)}]{Aavatsmark2002}
\bibinfo{author}{Aavatsmark, I.}, \bibinfo{year}{2002}.
\newblock \bibinfo{title}{An introduction to multipoint flux approximations for
  quadrilateral grids}.
\newblock \bibinfo{journal}{Computat. Geosci.} \bibinfo{volume}{6},
  \bibinfo{pages}{405--432}.
\newblock \DOIprefix\doi{10.1023/A:1021291114475}.
\bibitem[{Aavatsmark et~al.(2008)Aavatsmark, Eigestad, Mallison and
  Nordbotten}]{Aavatsmark2008}
\bibinfo{author}{Aavatsmark, I.}, \bibinfo{author}{Eigestad, G.},
  \bibinfo{author}{Mallison, B.}, \bibinfo{author}{Nordbotten, J.},
  \bibinfo{year}{2008}.
\newblock \bibinfo{title}{A compact multipoint flux approximation method with
  improved robustness}.
\newblock \bibinfo{journal}{Numer. Meth. Part. D. E.} \bibinfo{volume}{24},
  \bibinfo{pages}{1329--1360}.
\newblock \DOIprefix\doi{10.1002/num.20320}.
\bibitem[{Appleyard and Cheshire(1983)}]{Appleyard1983}
\bibinfo{author}{Appleyard, J.R.}, \bibinfo{author}{Cheshire, I.M.},
  \bibinfo{year}{1983}.
\newblock \bibinfo{title}{Nested factorization}, in: \bibinfo{booktitle}{SPE
  Reservoir Simulation Symposium}, \bibinfo{publisher}{SPE}.
\newblock \DOIprefix\doi{10.2118/12264-ms}.
\bibitem[{Aziz(1979)}]{Aziz1979}
\bibinfo{author}{Aziz, K.}, \bibinfo{year}{1979}.
\newblock \bibinfo{title}{Petroleum Reservoir Simulation}.
\newblock \bibinfo{publisher}{Applied Science Publishers}.
\bibitem[{Balay et~al.(2025)Balay, Abhyankar, Adams and
  et~al.}]{petsc_web_page2024}
\bibinfo{author}{Balay, S.}, \bibinfo{author}{Abhyankar, S.},
  \bibinfo{author}{Adams, M.F.}, \bibinfo{author}{et~al.},
  \bibinfo{year}{2025}.
\newblock \bibinfo{title}{{PETS}c {W}eb page}.
\newblock \URLprefix \url{https://petsc.org/}.
\bibitem[{Cai and Keyes(2002)}]{doi:10.1137/S106482750037620X}
\bibinfo{author}{Cai, X.C.}, \bibinfo{author}{Keyes, D.E.},
  \bibinfo{year}{2002}.
\newblock \bibinfo{title}{Nonlinearly preconditioned inexact {Newton}
  algorithms}.
\newblock \bibinfo{journal}{SIAM J. Sci. Comput.} \bibinfo{volume}{24},
  \bibinfo{pages}{183--200}.
\newblock \DOIprefix\doi{10.1137/S106482750037620X}.
\bibitem[{Cai and Li(2011)}]{doi:10.1137/080736272}
\bibinfo{author}{Cai, X.C.}, \bibinfo{author}{Li, X.}, \bibinfo{year}{2011}.
\newblock \bibinfo{title}{Inexact {Newton} methods with restricted additive
  {Schwarz} based nonlinear elimination for problems with high local
  nonlinearity}.
\newblock \bibinfo{journal}{SIAM J. Sci. Comput.} \bibinfo{volume}{33},
  \bibinfo{pages}{746--762}.
\newblock \DOIprefix\doi{10.1137/080736272}.
\bibitem[{Chen et~al.(2006)Chen, Huan and Ma}]{chen_computational_2006}
\bibinfo{author}{Chen, Z.}, \bibinfo{author}{Huan, G.}, \bibinfo{author}{Ma,
  Y.}, \bibinfo{year}{2006}.
\newblock \bibinfo{title}{Computational Methods for Multiphase Flows in Porous
  Media}.
\newblock \bibinfo{publisher}{SIAM}.
\newblock \DOIprefix\doi{10.1137/1.9780898718942}.
\bibitem[{Coats(2003)}]{CoatsCFL2003}
\bibinfo{author}{Coats, K.H.}, \bibinfo{year}{2003}.
\newblock \bibinfo{title}{{IMPES} stability: The {CFL} limit}.
\newblock \bibinfo{journal}{SPE Journal} \bibinfo{volume}{8},
  \bibinfo{pages}{291--297}.
\newblock \DOIprefix\doi{10.2118/85956-PA}.
\bibitem[{Dolean et~al.(2016)Dolean, Gander, Kheriji, Kwok and
  Masson}]{dolean:hal-01171167}
\bibinfo{author}{Dolean, V.}, \bibinfo{author}{Gander, M.J.},
  \bibinfo{author}{Kheriji, W.}, \bibinfo{author}{Kwok, F.},
  \bibinfo{author}{Masson, R.}, \bibinfo{year}{2016}.
\newblock \bibinfo{title}{Nonlinear preconditioning: How to use a nonlinear
  {Schwarz} method to precondition {Newton's} method}.
\newblock \bibinfo{journal}{SIAM J. Sci. Comput.} \bibinfo{volume}{38},
  \bibinfo{pages}{A3357--A3380}.
\newblock \DOIprefix\doi{10.1137/15M102887X}.
\bibitem[{Douglas et~al.(1959)Douglas, Peaceman and Rachford}]{DouglasFIM1959}
\bibinfo{author}{Douglas, Jim, J.}, \bibinfo{author}{Peaceman, D.},
  \bibinfo{author}{Rachford, H.H., J.}, \bibinfo{year}{1959}.
\newblock \bibinfo{title}{A method for calculating multi-dimensional immiscible
  displacement}.
\newblock \bibinfo{journal}{Trans. AIME} \bibinfo{volume}{216},
  \bibinfo{pages}{297--308}.
\newblock \DOIprefix\doi{10.2118/1327-G}.
\bibitem[{Falgout and Yang(2002)}]{Hypre2002}
\bibinfo{author}{Falgout, R.D.}, \bibinfo{author}{Yang, U.M.},
  \bibinfo{year}{2002}.
\newblock \bibinfo{title}{{HYPRE}: A library of high performance
  preconditioners}, in: \bibinfo{booktitle}{Computational Science --- ICCS
  2002}, \bibinfo{publisher}{Springer Berlin Heidelberg},
  \bibinfo{address}{Berlin, Heidelberg}. pp. \bibinfo{pages}{632--641}.
\newblock \DOIprefix\doi{10.1007/3-540-47789-6_66}.
\bibitem[{Feng et~al.(2024)Feng, Li, Liu, Zhang and Zhao}]{Feng2024}
\bibinfo{author}{Feng, C.}, \bibinfo{author}{Li, S.}, \bibinfo{author}{Liu,
  S.}, \bibinfo{author}{Zhang, C.}, \bibinfo{author}{Zhao, L.},
  \bibinfo{year}{2024}.
\newblock \bibinfo{title}{Application-oriented preconditioning of seepage
  mechanics}.
\newblock \bibinfo{journal}{Chinese J. Comput. Phys.} \bibinfo{volume}{41},
  \bibinfo{pages}{98--109}.
\newblock \DOIprefix\doi{10.19596/j.cnki.1001-246x.8791}.
\bibitem[{Feng et~al.(2014)Feng, Shu, Xu and Zhang}]{2014A2}
\bibinfo{author}{Feng, C.}, \bibinfo{author}{Shu, S.}, \bibinfo{author}{Xu,
  J.}, \bibinfo{author}{Zhang, C.}, \bibinfo{year}{2014}.
\newblock \bibinfo{title}{A multi-stage preconditioner for the black oil model
  and its {OpenMP} implementation}.
\newblock \bibinfo{journal}{Lect. Notes Comput. Sci. Eng.}
  \bibinfo{volume}{98}, \bibinfo{pages}{141--153}.
\newblock \DOIprefix\doi{10.1007/978-3-319-05789-7_11}.
\bibitem[{Hwang and Cai(2007)}]{hwang2007class}
\bibinfo{author}{Hwang, F.N.}, \bibinfo{author}{Cai, X.C.},
  \bibinfo{year}{2007}.
\newblock \bibinfo{title}{A class of parallel two-level nonlinear {Schwarz}
  preconditioned inexact {Newton} algorithms}.
\newblock \bibinfo{journal}{Comput. Methods Appl. Mech. Engrg.}
  \bibinfo{volume}{196}, \bibinfo{pages}{1603--1611}.
\newblock \DOIprefix\doi{10.1016/j.cma.2006.03.019}.
\bibitem[{Karypis and Kumar(2009)}]{karypis_metis_2009}
\bibinfo{author}{Karypis, G.}, \bibinfo{author}{Kumar, V.},
  \bibinfo{year}{2009}.
\newblock \bibinfo{title}{{METIS}: Unstructured graph partitioning and sparse
  matrix ordering system}.
\newblock \URLprefix \url{http://www.cs.umn.edu/ metis}.
\bibitem[{Karypis et~al.(2020)Karypis, Schloegel and Kumar}]{ParMETIS2020}
\bibinfo{author}{Karypis, G.}, \bibinfo{author}{Schloegel, K.},
  \bibinfo{author}{Kumar, V.}, \bibinfo{year}{2020}.
\newblock \bibinfo{title}{{ParMETIS}: Parallel graph partitioning and
  fill-reducing matrix ordering}.
\newblock \URLprefix \url{https://github.com/KarypisLab/ParMETIS}.
\bibitem[{Killough and Kossack(1987)}]{Killough_SPE5}
\bibinfo{author}{Killough, J.E.}, \bibinfo{author}{Kossack, C.A.},
  \bibinfo{year}{1987}.
\newblock \bibinfo{title}{Fifth comparative solution project: Evaluation of
  miscible flood simulators}.
\newblock \bibinfo{journal}{SPE Symposium on Reservoir Simulation.}
  \bibinfo{volume}{SPE-16000-MS}.
\newblock \DOIprefix\doi{10.2118/16000-MS}.
\bibitem[{Klemetsdal et~al.(2022)Klemetsdal, Moncorg{\'e}, M{\o}yner and
  Lie}]{klemetsdal2022numerical}
\bibinfo{author}{Klemetsdal, {\O}.}, \bibinfo{author}{Moncorg{\'e}, A.},
  \bibinfo{author}{M{\o}yner, O.}, \bibinfo{author}{Lie, K.A.},
  \bibinfo{year}{2022}.
\newblock \bibinfo{title}{A numerical study of the additive {Schwarz}
  preconditioned exact {Newton} method {(ASPEN)} as a nonlinear preconditioner
  for immiscible and compositional porous media flow}.
\newblock \bibinfo{journal}{Computat. Geosci.} \bibinfo{volume}{26},
  \bibinfo{pages}{1045--1063}.
\newblock \DOIprefix\doi{10.1007/s10596-021-10090-x}.
\bibitem[{Knoll and Keyes(2004)}]{KNOLL2004357}
\bibinfo{author}{Knoll, D.}, \bibinfo{author}{Keyes, D.}, \bibinfo{year}{2004}.
\newblock \bibinfo{title}{Jacobian-free {Newton--Krylov} methods: a survey of
  approaches and applications}.
\newblock \bibinfo{journal}{J. Comput. Phys.} \bibinfo{volume}{193},
  \bibinfo{pages}{357--397}.
\newblock \DOIprefix\doi{10.1016/j.jcp.2003.08.010}.
\bibitem[{Li(2024)}]{ShizhePHD2024}
\bibinfo{author}{Li, S.}, \bibinfo{year}{2024}.
\newblock \bibinfo{title}{Study of parallel numerical methods and software for
  multiphase multicomponent flow in porous media}.
\newblock Ph.D. thesis. University of Chinese Academy of Sciences.
\bibitem[{Li and Zhang(2024)}]{OCP2024}
\bibinfo{author}{Li, S.}, \bibinfo{author}{Zhang, C.S.}, \bibinfo{year}{2024}.
\newblock \bibinfo{title}{{OpenCAEPoro}: A parallel simulation framework for
  multiphase and multicomponent porous media flows}.
\newblock \bibinfo{journal}{Preprint}
  \href{http://arxiv.org/abs/2406.10862}{{\tt arXiv:2406.10862}}.
\bibitem[{Li et~al.(2017)Li, Wu, Zhang, Xu, Feng and Hu}]{LiZheng2017}
\bibinfo{author}{Li, Z.}, \bibinfo{author}{Wu, S.}, \bibinfo{author}{Zhang,
  C.}, \bibinfo{author}{Xu, J.}, \bibinfo{author}{Feng, C.},
  \bibinfo{author}{Hu, X.}, \bibinfo{year}{2017}.
\newblock \bibinfo{title}{Numerical studies of a class of linear solvers for
  fine-scale petroleum reservoir simulation}.
\newblock \bibinfo{journal}{Comput. Visualization Sci.} \bibinfo{volume}{18},
  \bibinfo{pages}{93--102}.
\newblock \DOIprefix\doi{10.1007/s00791-016-0273-3}.
\bibitem[{Liu et~al.(2024)Liu, Gao, Yu and Keyes}]{Liu2023OverlappingMS}
\bibinfo{author}{Liu, L.}, \bibinfo{author}{Gao, W.}, \bibinfo{author}{Yu, H.},
  \bibinfo{author}{Keyes, D.E.}, \bibinfo{year}{2024}.
\newblock \bibinfo{title}{Overlapping multiplicative {Schwarz} preconditioning
  for linear and nonlinear systems}.
\newblock \bibinfo{journal}{J. Comput. Phys.} \bibinfo{volume}{496},
  \bibinfo{pages}{112548}.
\newblock \DOIprefix\doi{10.1016/j.jcp.2023.112548}.
\bibitem[{Liu and Keyes(2015)}]{doi:10.1137/140970379}
\bibinfo{author}{Liu, L.}, \bibinfo{author}{Keyes, D.E.}, \bibinfo{year}{2015}.
\newblock \bibinfo{title}{Field-split preconditioned inexact {Newton}
  algorithms}.
\newblock \bibinfo{journal}{SIAM J. Sci. Comput.} \bibinfo{volume}{37},
  \bibinfo{pages}{A1388--A1409}.
\newblock \DOIprefix\doi{10.1137/140970379}.
\bibitem[{Liu and Keyes(2016)}]{liu2016convergence}
\bibinfo{author}{Liu, L.}, \bibinfo{author}{Keyes, D.E.}, \bibinfo{year}{2016}.
\newblock \bibinfo{title}{Convergence analysis for the multiplicative {Schwarz}
  preconditioned inexact {Newton} algorithm}.
\newblock \bibinfo{journal}{SIAM J. Numer. Anal.} \bibinfo{volume}{54},
  \bibinfo{pages}{3145--3166}.
\newblock \DOIprefix\doi{10.1137/15M1028182}.
\bibitem[{Liu et~al.(2018)Liu, Keyes and Krause}]{liu2018note}
\bibinfo{author}{Liu, L.}, \bibinfo{author}{Keyes, D.E.},
  \bibinfo{author}{Krause, R.}, \bibinfo{year}{2018}.
\newblock \bibinfo{title}{A note on adaptive nonlinear preconditioning
  techniques}.
\newblock \bibinfo{journal}{SIAM J. Sci. Comput.} \bibinfo{volume}{40},
  \bibinfo{pages}{A1171--A1186}.
\newblock \DOIprefix\doi{10.1137/17M1128502}.
\bibitem[{Luo et~al.(2021)Luo, Cai and Keyes}]{luo2021nonlinear}
\bibinfo{author}{Luo, L.}, \bibinfo{author}{Cai, X.C.}, \bibinfo{author}{Keyes,
  D.E.}, \bibinfo{year}{2021}.
\newblock \bibinfo{title}{Nonlinear preconditioning strategies for two-phase
  flows in porous media discretized by a fully implicit discontinuous
  {Galerkin} method}.
\newblock \bibinfo{journal}{SIAM J. Sci. Comput.} \bibinfo{volume}{43},
  \bibinfo{pages}{S317--S344}.
\newblock \DOIprefix\doi{10.1137/20M1344652}.
\bibitem[{Michelsen(1982a)}]{michelsen1982isothermal-psa}
\bibinfo{author}{Michelsen, M.L.}, \bibinfo{year}{1982}a.
\newblock \bibinfo{title}{The isothermal flash problem. part i. stability}.
\newblock \bibinfo{journal}{Fluid phase equilibria} \bibinfo{volume}{9},
  \bibinfo{pages}{1--19}.
\bibitem[{Michelsen(1982b)}]{michelsen1982isothermal-psc}
\bibinfo{author}{Michelsen, M.L.}, \bibinfo{year}{1982}b.
\newblock \bibinfo{title}{The isothermal flash problem. part ii. phase-split
  calculation}.
\newblock \bibinfo{journal}{Fluid Phase Equilibria} \bibinfo{volume}{9},
  \bibinfo{pages}{21--40}.
\newblock \URLprefix \url{10.1016/0378-3812(82)85002-4}.
\bibitem[{Odeh(1981)}]{odeh_comparison_1981}
\bibinfo{author}{Odeh, A.S.}, \bibinfo{year}{1981}.
\newblock \bibinfo{title}{Comparison of solutions to a three-dimensional
  black-oil reservoir simulation problem}.
\newblock \bibinfo{journal}{J. Pet. Technol.} \bibinfo{volume}{33},
  \bibinfo{pages}{13--25}.
\newblock \DOIprefix\doi{10.2118/9723-PA}.
\bibitem[{Peaceman(1978)}]{Peaceman_1978}
\bibinfo{author}{Peaceman, D.}, \bibinfo{year}{1978}.
\newblock \bibinfo{title}{Interpretation of well-block pressures in numerical
  reservoir simulation}.
\newblock \bibinfo{journal}{SPE Journal} \bibinfo{volume}{18},
  \bibinfo{pages}{183--194}.
\newblock \DOIprefix\doi{10.2118/6893-PA}.
\bibitem[{Qiao(2015)}]{Qiao2015}
\bibinfo{author}{Qiao, C.}, \bibinfo{year}{2015}.
\newblock \bibinfo{title}{General purpose compositional simulation for
  multiphase reactive flow with a fast linear solver}.
\newblock Ph.D. thesis. The Pennsylvania State University.
\bibitem[{Rasmussen et~al.(2021)Rasmussen, Sandve, Bao, Lauser, Hove,
  Skaflestad, Kl{\"o}fkorn, Blatt, Rustad, S{\ae}vareid, Lie and
  Thune}]{rasmussen2021open}
\bibinfo{author}{Rasmussen, A.F.}, \bibinfo{author}{Sandve, T.H.},
  \bibinfo{author}{Bao, K.}, \bibinfo{author}{Lauser, A.},
  \bibinfo{author}{Hove, J.}, \bibinfo{author}{Skaflestad, B.},
  \bibinfo{author}{Kl{\"o}fkorn, R.}, \bibinfo{author}{Blatt, M.},
  \bibinfo{author}{Rustad, A.B.}, \bibinfo{author}{S{\ae}vareid, O.},
  \bibinfo{author}{Lie, K.A.}, \bibinfo{author}{Thune, A.},
  \bibinfo{year}{2021}.
\newblock \bibinfo{title}{The open porous media flow reservoir simulator}.
\newblock \bibinfo{journal}{Comput. Math. Appl.} \bibinfo{volume}{81},
  \bibinfo{pages}{159--185}.
\newblock \DOIprefix\doi{10.1016/j.camwa.2020.05.014}.
\bibitem[{Saad(2003)}]{Saad2003}
\bibinfo{author}{Saad, Y.}, \bibinfo{year}{2003}.
\newblock \bibinfo{title}{Iterative Methods for Sparse Linear Systems}.
\newblock \bibinfo{edition}{Second} ed., \bibinfo{publisher}{SIAM}.
\newblock \DOIprefix\doi{10.1137/1.9780898718003}.
\bibitem[{Schlumberger(2021)}]{eclipse_manual_2021}
\bibinfo{author}{Schlumberger}, \bibinfo{year}{2021}.
\newblock \bibinfo{title}{ECLIPSE Technical Description}.
\bibitem[{Siek et~al.(2002)Siek, Lee and Lumsdaine}]{siek_boost_2002}
\bibinfo{author}{Siek, J.}, \bibinfo{author}{Lee, L.Q.},
  \bibinfo{author}{Lumsdaine, A.}, \bibinfo{year}{2002}.
\newblock \bibinfo{title}{The {Boost} {Graph} {Library}: {User} {Guide} and
  {Reference} {Manual}}.
\newblock \bibinfo{publisher}{Addison-Wesley}.
\bibitem[{Skogestad et~al.(2013)Skogestad, Keilegavlen and
  Nordbotten}]{skogestad2013domain}
\bibinfo{author}{Skogestad, J.O.}, \bibinfo{author}{Keilegavlen, E.},
  \bibinfo{author}{Nordbotten, J.M.}, \bibinfo{year}{2013}.
\newblock \bibinfo{title}{Domain decomposition strategies for nonlinear flow
  problems in porous media}.
\newblock \bibinfo{journal}{J. Comput. Phys.} \bibinfo{volume}{234},
  \bibinfo{pages}{439--451}.
\newblock \DOIprefix\doi{10.1016/j.jcp.2012.10.001}.
\bibitem[{Wallis(1983)}]{Wallis1983}
\bibinfo{author}{Wallis, J.}, \bibinfo{year}{1983}.
\newblock \bibinfo{title}{Incomplete {Gaussian} elimination as a
  preconditioning for generalized conjugate gradient acceleration}.
\newblock \bibinfo{journal}{SPE Reservoir Simulation Conference}
  \bibinfo{volume}{SPE-12265}.
\newblock \DOIprefix\doi{10.2118/12265-MS}.
\bibitem[{Wallis et~al.(1985)Wallis, Kendall and Little}]{Wallis1985}
\bibinfo{author}{Wallis, J.}, \bibinfo{author}{Kendall, R.},
  \bibinfo{author}{Little, T.}, \bibinfo{year}{1985}.
\newblock \bibinfo{title}{Constrained residual acceleration of conjugate
  residual methods}.
\newblock \bibinfo{journal}{SPE Reservoir Simulation Conference}
  \bibinfo{volume}{SPE-13536}.
\newblock \DOIprefix\doi{10.2118/13536-MS}.
\bibitem[{Wang et~al.(2018)Wang, Liu, Luo and Chen}]{WANG2018443}
\bibinfo{author}{Wang, K.}, \bibinfo{author}{Liu, H.}, \bibinfo{author}{Luo,
  J.}, \bibinfo{author}{Chen, Z.}, \bibinfo{year}{2018}.
\newblock \bibinfo{title}{Efficient {CPR-type} preconditioner and its adaptive
  strategies for large-scale parallel reservoir simulations}.
\newblock \bibinfo{journal}{J. Comput. Appl. Math.} \bibinfo{volume}{328},
  \bibinfo{pages}{443--468}.
\newblock \DOIprefix\doi{10.1016/j.cam.2017.07.022}.
\bibitem[{Xu and Zou(1998)}]{xu1998some}
\bibinfo{author}{Xu, J.}, \bibinfo{author}{Zou, J.}, \bibinfo{year}{1998}.
\newblock \bibinfo{title}{Some nonoverlapping domain decomposition methods}.
\newblock \bibinfo{journal}{SIAM Review} \bibinfo{volume}{40},
  \bibinfo{pages}{857--914}.
\newblock \DOIprefix\doi{10.1137/S0036144596306800}.
\bibitem[{Yang et~al.(2014)Yang, Moridis and Blasingame}]{YANG2014417}
\bibinfo{author}{Yang, D.}, \bibinfo{author}{Moridis, G.J.},
  \bibinfo{author}{Blasingame, T.A.}, \bibinfo{year}{2014}.
\newblock \bibinfo{title}{A fully coupled multiphase flow and geomechanics
  solver for highly heterogeneous porous media}.
\newblock \bibinfo{journal}{J. Comput. Appl. Math.} \bibinfo{volume}{270},
  \bibinfo{pages}{417--432}.
\newblock \DOIprefix\doi{10.1016/j.cam.2013.12.029}.
\bibitem[{Yang and Hwang(2018)}]{yang2018adaptive}
\bibinfo{author}{Yang, H.}, \bibinfo{author}{Hwang, F.N.},
  \bibinfo{year}{2018}.
\newblock \bibinfo{title}{An adaptive nonlinear elimination preconditioned
  inexact {Newton} algorithm for highly local nonlinear multicomponent {PDE}
  systems}.
\newblock \bibinfo{journal}{Appl. Numer. Math.} \bibinfo{volume}{133},
  \bibinfo{pages}{100--115}.
\newblock \DOIprefix\doi{10.1016/j.apnum.2018.01.008}.
\bibitem[{Yang et~al.(2018)Yang, Sun, Li and Yang}]{YANG2018}
\bibinfo{author}{Yang, H.}, \bibinfo{author}{Sun, S.}, \bibinfo{author}{Li,
  Y.}, \bibinfo{author}{Yang, C.}, \bibinfo{year}{2018}.
\newblock \bibinfo{title}{A scalable fully implicit framework for reservoir
  simulation on parallel computers}.
\newblock \bibinfo{journal}{Comput. Methods Appl. Mech. Engrg.}
  \bibinfo{volume}{330}, \bibinfo{pages}{334--350}.
\newblock \DOIprefix\doi{10.1016/j.cma.2017.10.016}.
\bibitem[{Yang et~al.(2019)Yang, Sun, Li and Yang}]{YANG20192}
\bibinfo{author}{Yang, H.}, \bibinfo{author}{Sun, S.}, \bibinfo{author}{Li,
  Y.}, \bibinfo{author}{Yang, C.}, \bibinfo{year}{2019}.
\newblock \bibinfo{title}{Parallel reservoir simulators for fully implicit
  complementarity formulation of multicomponent compressible flows}.
\newblock \bibinfo{journal}{Comput. Phys. Commun.} \bibinfo{volume}{244},
  \bibinfo{pages}{2--12}.
\newblock \DOIprefix\doi{10.1016/j.cpc.2019.07.011}.
\bibitem[{Yang et~al.(2016)Yang, Yang and Sun}]{yang2016active}
\bibinfo{author}{Yang, H.}, \bibinfo{author}{Yang, C.}, \bibinfo{author}{Sun,
  S.}, \bibinfo{year}{2016}.
\newblock \bibinfo{title}{Active-set reduced-space methods with nonlinear
  elimination for two-phase flow problems in porous media}.
\newblock \bibinfo{journal}{SIAM J. Sci. Comput.} \bibinfo{volume}{38},
  \bibinfo{pages}{B593--B618}.
\newblock \DOIprefix\doi{10.1137/15M1041882}.
\bibitem[{Zhang(2022)}]{Chensong2022}
\bibinfo{author}{Zhang, C.}, \bibinfo{year}{2022}.
\newblock \bibinfo{title}{Linear solvers for petroleum reservoir simulation}.
\newblock \bibinfo{journal}{J. Numer. Methods Comput. Appl.}
  \bibinfo{volume}{43}, \bibinfo{pages}{1--26}.
\newblock \DOIprefix\doi{10.12288/szjs.s2021-0813}.
\bibitem[{Zhang et~al.(2022)Zhang, Yang, Wu and Sun}]{ZHANG2022}
\bibinfo{author}{Zhang, M.}, \bibinfo{author}{Yang, H.}, \bibinfo{author}{Wu,
  S.}, \bibinfo{author}{Sun, S.}, \bibinfo{year}{2022}.
\newblock \bibinfo{title}{Parallel multilevel domain decomposition
  preconditioners for monolithic solution of non-isothermal flow in reservoir
  simulation}.
\newblock \bibinfo{journal}{Comput. Fluids} \bibinfo{volume}{232},
  \bibinfo{pages}{105183}.
\newblock \DOIprefix\doi{10.1016/j.compfluid.2021.105183}.
\bibitem[{Zhao et~al.(2022)Zhao, Feng, Zhang and Shu}]{Zhao2022}
\bibinfo{author}{Zhao, L.}, \bibinfo{author}{Feng, C.}, \bibinfo{author}{Zhang,
  C.S.}, \bibinfo{author}{Shu, S.}, \bibinfo{year}{2022}.
\newblock \bibinfo{title}{Parallel multi-stage preconditioners with adaptive
  setup for the black oil model}.
\newblock \bibinfo{journal}{Comput. Geosci.} \bibinfo{volume}{168},
  \bibinfo{pages}{105230}.
\newblock \DOIprefix\doi{10.1016/j.cageo.2022.105230}.
\bibitem[{Zhao et~al.(2023)Zhao, Li, Zhang, Feng and Shu}]{Zhao2023}
\bibinfo{author}{Zhao, L.}, \bibinfo{author}{Li, S.}, \bibinfo{author}{Zhang,
  C.S.}, \bibinfo{author}{Feng, C.}, \bibinfo{author}{Shu, S.},
  \bibinfo{year}{2023}.
\newblock \bibinfo{title}{An improved multistage preconditioner on {GPUs} for
  compositional reservoir simulation}.
\newblock \bibinfo{journal}{CCF Trans. High Perform. Comput.}
  \bibinfo{volume}{5}, \bibinfo{pages}{144--159}.
\newblock \DOIprefix\doi{10.1007/s42514-023-00136-0}.

\end{thebibliography}

\end{document}